\documentclass[]{siamltex}
\usepackage{./yuyuan}
\usepackage[colorlinks=true, citecolor=blue]{hyperref}
\usepackage{algorithm, algpseudocode, color, tabularx, multirow, graphicx,placeins,caption,subcaption,epstopdf}
\renewcommand{\thefigure}{\arabic{figure}}
\graphicspath{{./figures/OverLASSO/}{./figures/CS/}{./figures/PPI/}}

\renewcommand\epsilon\varepsilon
\newcommand\ABb[1][]{Bw#1 - Kx#1 - b}
\newcommand\BKb[1][]{Bw#1 - Kx#1 - b}
\def\tx{\tilde x}
\def\tw{\tilde w}
\def\ty{\tilde y}

\newcommand\cBx[1][]{\cB_t(x{#1}, x_{[t+1]}, \eta_{[t]})}
\newcommand\cBKx[1][]{\cB_t(Kx{#1}, Kx_{[t+1]}, \theta_{[t]})}

\newcommand\cBw[1][]{\cB_t(w{#1}, w_{[t+1]}, \theta_{[t]})}
\newcommand\cBy[1][]{\cB_t(y{#1}, y_{[t+1]}, \rho_{[t]}^{-1})}

\newcommand\taurho[1][t]{\left(\frac{\tau_{#1}}{\rho_{#1}}-1\right)}
\newcommand\thetarho[1][t]{\left(\frac{\theta_{#1}}{\rho_{#1}}-1\right)}

\title{An Accelerated Linearized Alternating Direction Method of Multipliers}
\author{Yuyuan Ouyang\thanks{Department of Industrial and System Engineering, University of Florida ({\tt ouyang@ufl.edu})}. Part of the research was done while the author was a PhD student at the Department of Mathematics, University of Florida. This author was partially supported by AFRL Mathematical Modeling Optimization Institute. \and Yunmei Chen\thanks{Department of Mathematics, University of Florida ({\tt yun@math.ufl.edu}). This author was partially supported by NSF grants DMS-1115568, IIP-1237814 and DMS-1319050}. \and
Guanghui Lan\thanks{Department of Industrial and System Engineering, University of Florida ({\tt glan@ise.ufl.edu}).
This author was partially supported by NSF grant CMMI-1000347, ONR grant N00014-13-1-0036, NSF DMS-1319050, and NSF CAREER Award CMMI-1254446.}. \and Eduardo Pasiliao Jr. \thanks{Munitions Directorate, Air Force Research Laboratory ({\tt eduardo.pasiliao@eglin.af.mil})}
}

\begin{document}

\maketitle

\begin{abstract}
We present a novel framework, namely AADMM, for acceleration of linearized alternating direction method of multipliers (ADMM). The basic idea of AADMM is to incorporate a multi-step acceleration scheme into linearized ADMM. We demonstrate that for solving a class of convex composite optimization with linear constraints, the rate of convergence of AADMM is better than that of linearized ADMM, in terms of their dependence on the Lipschitz constant of the smooth component. Moreover, AADMM is capable to deal with the situation when the feasible region is unbounded, as long as the corresponding saddle point problem has a solution. A backtracking algorithm is also proposed for practical performance. 
\end{abstract}

\section{Introduction}
	Assume that $\cal W$, $\cal X$ and $\cal Y$ are finite dimensional vectorial spaces equipped with inner product $\langle\cdot,\cdot\rangle$, norm $\|\cdot\|$ and conjugate norm $\|\cdot\|_*$. Our problem of interest is the following affine equality constrained composite optimization (AECCO) problem:
	\begin{align}
		\label{eqnAECCO}
		\min_{x\in X, w\in \cW}G(x) + F(w), \st Bw - Kx = b,
	\end{align}
	where $X\subseteq{\cal X}$ is a closed convex set, $G(\cdot):X\to\R$ and $F(\cdot):{\cal W}\to\R$ are finitely valued, convex and lower semi-continuous functions, and $K:X\to\cal Y$, $B:\cW\to\cal Y$ are bounded linear operators. 
	
	In this paper, we assume that $F(\cdot)$ is simple, in the sense that the optimization problem
		\begin{align}
			\label{eqnF}
			\min_{w\in \cW}\frac{\eta}{2}\|w - c\|^2 + F(w), \text{ where }c\in \cW, \eta\in \R
		\end{align}
		can be solved efficiently. We will use the term ``simple" in this sense throughout this paper, and use the term ``non-simple'' in the opposite sense. We assume that $G(\cdot)$ is non-simple, continuously differentiable, and that there exists $L_G>0$ such that
		\begin{align}
			\label{eqnLG}
			G(x_2) - G(x_1) - \langle\nabla G(x_1), x_2 - x_1\rangle \leq \frac{L_G}{2}\|x_2 - x_1\|^2,\ \forall x_1\in X, x_2\in X.
		\end{align}

	One special case of the AECCO problem in \eqref{eqnAECCO} is when  $B=I$ and $b=0$. Under this situation, problem \eqref{eqnAECCO} is equivalent to the following unconstrained composite optimization (UCO) problem:
	\begin{align}
		\label{eqnUCO}
		\min_{x\in X}f(x): = G(x) + F(Kx).
	\end{align}

	Both AECCO and UCO can be reformulated as saddle point problems. By the method of Lagrangian multipliers, the AECCO problem \eqref{eqnAECCO} is equivalent to the following saddle point problem:
				\begin{align}
					\label{eqnSPP}
					\min_{x\in X, w\in \cW}\max_{y\in \cY}G(x) + F(w) - \langle y, Bw - Kx - b\rangle.
				\end{align}
	The AECCO and UCO problems have found numerous applications in machine learning and image processing. In most application, $G(\cdot)$ is known as the fidelity term and $F(\cdot)$ is the regularization term. For example, consider the following two dimensional total variation (TV) based image reconstruction problem
	\begin{align}
		\label{eqnTV}
		\min_{x\in \F^{n}}\frac 12\|Ax - c\|^2 + \lambda\|D x\|_{2,1}, 
	\end{align}
	where the field $\F$ is either $\R$ or $\C$, $x$ is the n-vector form of a two-dimensional complex or real valued image, $D:\F^n\to \F^{2n}$ is the two-dimensional finite difference operator acting on the image $x$, and
	\begin{align*}
		\|y\|_{2,1}:=\sum_{i=1}^{n}\|(y^{(2i-1)},y^{(2i)})^T\|_2,\ \forall  y\in\F^{2n},
	\end{align*}
	where $\|\cdot\|_2$ is the Euclidean norm in $\R^2$. In \eqref{eqnTV}, the regularization term $\|D x\|_{2,1}$ is the discrete form of TV semi-norm. By setting $G(x):=\|Ax-c\|^2/2$, $F(\cdot):=\|\cdot\|_{2,1}$, $K=\lambda D$, $X=\cX= \F^n$ and $\cW=\F^{2n}$, problem \eqref{eqnTV} becomes a UCO problem in \eqref{eqnUCO}. 
		
\subsection{Notations and terminologies}
\label{secNotation}
In this subsection, we describe some necessary assumptions, notations and terminologies that will be used throughout this paper.

We assume that there exists an optimal solution $(w^*, x^*)$ of \eqref{eqnAECCO} and that there exists $y^*\in\cY$ such that $z^*:=(w^*, x^*, y^*)\in\cal Z$ is a saddle point of \eqref{eqnSPP}, where $\cal Z := \cal W\times\cal X\times\cal Y$. We also use the notation $Z:=\cW\times X\times Y$ if a set $Y\subseteq \cY$ is declared readily. We use $f^*:=G(x^*) + F(w^*)$ to denote the optimal objective value of problem \eqref{eqnAECCO}. Since UCO problems \eqref{eqnUCO} are special cases of AECCO \eqref{eqnAECCO}, we will also use $f^*$ to denote the optimal value $G(x^*) + F(Kx^*)$. 

In view of \eqref{eqnAECCO}, both the objective function value and the feasibility of the constraint should be considered when defining approximate solutions of AECCO, henceforth the following definition comes naturally:
\begin{dfn}
	\label{defSol}
	A pair $(w, x)\in \cW\times X$ is called an $(\epsilon, \delta)$-solution of \eqref{eqnAECCO} if
	\begin{align*}
		G(x) + F(w) - f^* \leq \epsilon,\text{ and } \|Bw - Kx - b\|\leq\delta.
	\end{align*}
	We say that $(w, x)$ has primal residual $\epsilon$ and feasibility residual $\delta$. In particular, if $(w, x)$ is an $(\epsilon, 0)$-solution, then we simply say that it is an $\epsilon$-solution.
\end{dfn}

The feasibility residual $\delta$ in Definition \ref{defSol} measures the violation of the equality constraint, and the primal residual $\epsilon$ measures the gap between the objective value $G(x)+F(w)$ at the approximate solution and the optimal value $f^*$. For an $(\epsilon, \delta)$-solution $(w, x)$ where $\delta>0$, since $(w,x)$ does not satisfy the equality constraint in \eqref{eqnAECCO}, it is possible that $G(x)+F(w)-f^*<0$. However, as pointed out in \cite{lan2013iteration}, a lower bound of $G(x)+F(w)-f^*$ is given by
		\begin{align*}
			G(x)+F(w)-f^*\geq \langle y^*, Bw-Kx-b\rangle\geq -\delta\|y^*\|,
		\end{align*}
where $y^*$ is a component of $z^*=(w^*,x^*,y^*)$, a saddle point of \eqref{eqnSPP}.

In the remainder of this subsection, we introduce some notations that will be used throughout this paper. The following distance constants will be used for simplicity:
\begin{align}
	\label{eqnD}
	\begin{aligned}
	D_{w^*, B}:=\|B(w_1 - w^*)\|, D_{x^*, K}:=\|K(x_1 - x^*)\|, D_{x^*}:=\|x_1 - x^*\|, D_{y^*}:=\|y_1 - y^*\|, 
	\\
	D_{X, K}:=\sup_{x_1, x_2\in X}\|Kx_1 - Kx_2\|, \text{ and }
	 D_S:=\sup_{s_1, s_2\in S}\|s_1 - s_2\|, \text{ for any compact set }S.
	\end{aligned}
\end{align}
For example, for any compact set $Y\subset \cY$, we use $D_Y$ to denote the diameter of $Y$. In addition, we use $x_{[t]}$ to denote sequence $\{x_i \}_{i=1}^t$, where $x_i$'s may either be real numbers, or points in vectorial spaces. We will also equip a few operations on the notation of sequences. Firstly, suppose that $\cV_1$, $\cV_2$ are any vector spaces, $v_{[t+1]}\subset \cV_1$ is any sequence in $\cV_1$ and $\cA:\cV_1\to\cV_2$ is any operator, we use $\cA v_{[t+1]}$ to denote the sequence $\{\cA v_i \}_{i=1}^{t+1}$. Secondly, if $\eta_{[t]}, \tau_{[t]}\subset \R$ are any real valued sequences, and $L\in \R$ is any real number, then  $\eta_{[t]} - L\tau_{[t]}$ denotes $\{\eta_i - L\tau_i \}_{i=1}^t $. Finally, we denote by $\eta_{[t]}^{-1}$ the reciprocal sequence $\{\eta_{i}^{-1}\}_{i=1}^t$ for any non-zero real valued sequence $\eta_{[t]}$.

	\subsection{Augmented Lagrangian and alternating direction method of multipliers}
	\label{secADMMintro}

	In this paper, we study AECCO problems from the aspect of the augmented Lagrangian formulation of \eqref{eqnSPP}:
	\begin{align}
		\label{eqnAL}
		\min_{x\in X, w\in \cW}\max_{y\in \cY}G(x) + F(w) - \langle y, Bw - Kx - b\rangle + \frac{\rho}{2}\|Bw - Kx - b\|^2,
	\end{align}	
		where $\rho$ is a penalty parameter. The idea of analyzing \eqref{eqnAL} in order to solve \eqref{eqnAECCO} is essentially the augmented Lagrangian method (ALM) by Hestenes \cite{hestenes1969multiplier} and Powell \cite{powell1969method} (It is originally called the method of multipliers in \cite{hestenes1969multiplier, powell1969method}; see also the textbooks, e.g., \cite{bertsekas1982constrained, nocedal2006numerical, bertsekas1999nonlinear}). The ALM is a special case of the Douglas-Rachford splitting method \cite{gabay1983applications, douglas1956numerical, lions1979splitting}, which is also an instance of the proximal point algorithm \cite{eckstein1992douglas, rockafellar1976monotone}. The iteration complexity of an inexact version of ALM, where the subproblems are solved iteratively by Nesterov's method, has been studied in \cite{lan2009iteration}. One influential variant of ALM is the ADMM algorithm  \cite{gabay1976dual, glowinski1975approximation}, which is an alternating method for solving \eqref{eqnAL} by minimizing $x$ and $w$ alternatively and then updating the Lagrangian coefficient $y$ (See \cite{boyd2011distributed} for a comprehensive explanation on ALM, ADMM and other algorithms). In compressive sensing and imaging science, the class of Bregman iterative methods is an application of the ALM and the ADMM. In particular, the Bregman iterative method \cite{goldstein2009split} is equivalent to ALM, and the split Bregman method \cite{goldstein2012fast} is equivalent to ADMM. 
	
		We give a brief review on ADMM, and some of its variants. The scheme of ADMM is described in Algorithm \ref{algADMM}.

		\begin{algorithm}
			\caption{\label{algADMM} The alternating direction method of multipliers (ADMM) for solving \eqref{eqnAECCO} }	
			\begin{algorithmic}[]
				\State Choose $x_1\in X$, $w_1\in \cW$ and $y_1\in \cY$.
				\For {$t=1,\ldots,N-1$}
					\begin{align}
						\label{eqnADMMxtn}
						x\tn =&\ \argmin_{x\in X}G(x) - \langle y_t, Bw_t - Kx - b\rangle + \frac{\rho}{2}\|Bw_t - Kx - b\|^2,
						\\
						\label{eqnADMMwtn}
						w\tn =&\ \argmin_{w\in \cW}F(w) - \langle y_t, Bw - Kx\tn - b\rangle + \frac{\rho}{2}\|Bw - Kx\tn - b\|^2,
						\\
						\label{eqnADMMytn}
						y\tn =&\ y_t - \rho(\BKb[\tn]).
					\end{align}
				\EndFor
			\end{algorithmic}
		\end{algorithm}
		
	For non-simple $G$, a linearized ADMM (L-ADMM) scheme generates iterate $x\tn$ in \eqref{eqnADMMxtn} by
			\begin{align}
				\label{eqnLADMMxtn}
				x\tn =&\ \argmin_{x\in X}\langle\nabla G(x_t), x\rangle + \langle y_t, Kx\rangle + \frac{\rho}{2}\|Bw_t - Kx - b\|^2 + \frac{\eta}{2}\|x - x_t\|^2.
			\end{align}
	We may also linearize $\|Bw_t - Kx - b\|^2$, and generate $x\tn$ by
			\begin{align}
				\label{eqnPADMMxtn}
				x\tn =&\ \argmin_{x\in X}G(x) + \langle y_t, Kx\rangle - \rho\langle Bw_t - Kx_t - b, Kx\rangle + \frac{\eta}{2}\|x - x_t\|^2,
			\end{align}
	as discussed in \cite{esser2010general,chambolle2011first}. This variant is called the preconditioned ADMM (P-ADMM). If we linearize both $G(x)$ and $\|Bw_t - Kx - b\|^2$, we have the linearized preconditioned ADMM (LP-ADMM), in which \eqref{eqnADMMxtn} is changed to
			\begin{align}
				\label{eqnLPADMMxtn}
				x\tn =&\ \argmin_{x\in X}\langle\nabla G(x_t), x\rangle +  \langle y_t, Kx\rangle - \rho\langle Bw_t - Kx_t - b, Ax\rangle + \frac{\eta}{2}\|x - x_t\|^2.
			\end{align}

	There has been several works on the convergence analysis and applications of ADMM, L-ADMM, and P-ADMM. It is shown in \cite{chambolle2011first} that P-ADMM (Algorithm 1 with $\theta=1$ in \cite{chambolle2011first}) solves the UCO problem with rate of convergence  
			\[
				\cO\left(\frac{\|K\|D^2}{N}\right),
			\]
			where $N$ is the number of iterations and $D$ depends on the distances $D_{x^*}$ and $D_{y^*}$. There are also several works concerning the tuning of the stepsize $\eta_t$ in L-ADMM, including \cite{ye2011computational, ye2011fast, chen2012fast}.

	For AECCO problems, in \cite{monteiro2013iteration} ADMM is treated as an instance of block-decomposition hybrid proximal extragradient (BD-HPE), and it is proved that the rate of convergence of the primal residual of ADMM for solving AECCO is
			\[
				\cO\left(\frac{D^2}{N}\right),
			\]
			where $D$ depends on $B$, $D_{x^*}$ and $D_{y^*}$. In \cite{he2012convergence}, the convergence analysis of ADMM and P-ADMM is studied based on the variational inequality formulation of \eqref{eqnSPP}, in which similar rate of convergence is achieved under the assumption that both the primal and dual feasible sets in \eqref{eqnSPP} are bounded. In \cite{ouyang2013stochastic}, it is shown that if $X$ is compact, then the rate of convergence of ADMM and L-ADMM for solving the AECCO problem is
				\begin{align}
					\label{eqnRateO}
					G(x^N) + F(w^N) - f^* + \rho\|Bw^N - Kx^N - b\|^2\leq\cO\left(\frac{L_GD_X^2+\rho D_{y^*, B}^2}{N} \right), \forall \rho>0,
				\end{align}
				where $(x^N, w^N)$ is the average of iterates $x_{[N]}$ of the ADMM algorithm. The result in \eqref{eqnRateO} is stronger than the results in \cite{monteiro2013iteration,he2012convergence}, in the sense that both primal and feasibility residuals are included in \eqref{eqnRateO}, while in \cite{monteiro2013iteration,he2012convergence} there is no discussion on the feasibility residual. However, the rate of convergence of the feasibility residual is still not very clear in \eqref{eqnRateO}, considering that $G(x^N) + F(w^N) - f^*$ can be negative.
	
\subsection{Accelerated methods for AECCO and UCO problems}
In a seminal paper \cite{nesterov1983method}, Nesterov introduced a smoothing technique and a fast first-order method that solves a class of composite optimization. When applied to UCO problems, Nesterov's method has optimal rate of convergence
\begin{align}
	\label{eqnUCOOptRate}
	\cO\left(\frac{L_GD_{x^*}^2}{N^2}+\frac{\|K\|D_{x^*}D_Y}{N}\right) \footnotemark,
\end{align}
\footnotetext{It is assumed in \cite{nesterov2005smooth} that $X$ is compact, hence the rate of convergence is dependent on $D_X$. However, the analysis in \cite{nesterov2005smooth} is also applicable for the case when $X$ is unbounded, yielding \eqref{eqnUCOOptRate}.}
\hspace{-.12cm}where $Y$ is the bounded dual space of the UCO problem. Following the breakthrough in \cite{nesterov2005smooth}, much effort has been devoted to the development of more efficient first-order methods for non-smooth optimization (see, e.g., \cite{nesterov2005excessive, auslender2006interior, lan2011primal, daspremont2008smooth, pena2008nash, tseng2008accelerated, becker2011nesta, lan2013bundle}). Although the rate in \eqref{eqnUCOOptRate} is also $\cO(1/N)$, what makes it more attractive is that it allows very large Lipschitz constant $L_G$. In particular, $L_G$ can be as large as $\Omega(N)$, without affecting the rate of convergence (up to a constant factor). However, it should be noted that the boundedness of $Y$ is critical for the convergence analysis of Nesterov's smoothing scheme. Following \cite{nesterov2005smooth}, there has also been several studies on the AECCO and UCO problems, and it has been shown that better acceleration results can be obtained if more assumptions are enforced for the AECCO and UCO problem. We give a list of such assumptions and results.
\begin{enumerate}
	\renewcommand{\theenumi}{\arabic{enumi})}
	\item \emph{Excessive gap technique.} The excessive gap technique is proposed in \cite{nesterov2005excessive} for solving the UCO problem in which $G$ is simple. Comparing to \cite{nesterov2005smooth}, the method in \cite{nesterov2005excessive} does not require the total number of iterations $N$ to be fixed in advance. Furthermore, if $G(\cdot)$ is strongly convex, it is shown that the rate of convergence of the excessive gap technique is $\cO(1/N^2)$.
	
	\item \emph{Special instance.} For the UCO problem, if $K=I$ and $G$ is simple, an accelerated method with skipping steps is proposed in Algorithm 7 of \cite{goldfarb2010fast}, which achieves $\cO(1/N^2)$ rate of convergence. The result is better than \eqref{eqnUCOOptRate}, but with cost of evaluating objective value functions in each iteration. For AECCO problem with compact feasible sets, it is shown in \cite{luo2012linear} that if $G(\cdot)$ is a composition of a strictly convex function and a linear transformation and $F(\cdot)$ is the weighted sum of 1-norm and some 2-norms, the asymptotic rate of convergence of ADMM method and its variants is R-linear.

	\item \emph{Strong convexity.} In \cite{chambolle2011first} for solving the UCO problem in which $G$ is simple, the authors showed that P-ADMM is equivalent to their proposed method, and furthermore, if either $G(\cdot)$ or $F^*(\cdot)$ is uniformly convex, then the rate of convergence of their method can be accelerated to $\cO(1/N^2)$. It is worth noting that this rate of convergence is weaker since it uses a different termination criterion. In addition, if both $G(\cdot)$ and $F^*(\cdot)$ are uniformly convex (hence the objective function in \eqref{eqnUCO} is continuously differentiable), the proposed method in \cite{chambolle2011first} converges linearly. When both $G(x)$ and $F(x)$ are strongly convex in the AECCO problem, an accelerated ADMM method is proposed in \cite{goldstein2012fast}, which achieves the $\cO(1/N^2)$ rate of convergence.
	
\end{enumerate}

It should be noted that all the methods in the above list require more assumptions on the AECCO and UCO 
problems (e.g., simplicity of $G(\cdot)$, strong convexity of $G(\cdot)$ or $F(\cdot)$), in comparison with Nesterov's smoothing scheme. More recently, we proposed an accelerated primal-dual (APD) method for solving the UCO problem \cite{chenoptimal}, which has the same optimal rate of convergence \eqref{eqnUCOOptRate} as that of Nesterov's smoothing scheme in \cite{nesterov2005smooth}. The advantage of the APD method over Nesterov's smoothing scheme is that it does not require boundedness on either $X$ or $Y$. The basic idea of the APD method is to incorporate a multi-step acceleration into LP-ADMM, and this has motivated our studies on accelerating the linearized ADMM method for solving the AECCO and UCO problems.

\subsection{Contribution of the paper}
The main interest of this paper is to develop an accelerated linearized ADMM algorithm for solving AECCO and UCO problems, in which $G$ is a general convex and non-simple function. Our contribution in this paper mainly consists of the following aspects. 

Firstly, we propose an accelerated framework for ADMM (AADMM), which consists two novel accelerated linearized ADMM methods, namely, accelerated L-ADMM (AL-ADMM) and accelerated LP-ADMM (ALP-ADMM). We
prove that AL-ADMM and ALP-ADMM have better rates of convergence than L-ADMM and LP-ADMM in terms of their dependence on $L_G$. In particular, we prove that both accelerated methods can achieve rates  similar to \eqref{eqnUCOOptRate}, hence both of them can efficiently solve problems with large Lipschitz constant $L_G$ (as large as $\Omega(N)$). We show that L-ADMM and LP-ADMM are special instances of AL-ADMM and ALP-ADMM respectively, with rates of convergence $\cO(1/N)$. To improve the performance in practice, we also propose a simple backtracking technique for searching Lipschitz constants $L_G$ and $\|K\|$. 

Secondly, the proposed framework solve both AECCO and UCO problems with unbounded feasible sets, as long as a saddle point of problem \eqref{eqnSPP} exists. Instead of using the perturbation type gap function in \cite{chenoptimal}, our convergence analysis is performed directly on both the primal and feasibility residuals. The estimate of the rate of convergence will depend on the distance from the initial point to the set of optimal solutions. 

\section{An accelerated ADMM framework}
\label{secFramework}

In this section, we propose an accelerated ADMM framework for solving AECCO \eqref{eqnAECCO} and UCO  \eqref{eqnUCO}. The proposed framework, namely AADMM, is presented in Algorithm \ref{algAADMM}.
\begin{algorithm}
	\caption{\label{algAADMM} Accelerated ADMM (AADMM) framework}
	\begin{algorithmic}[]
			\State Choose $x_1\in X$ and $w_1\in \cW$ such that $Bw_1=Kx_1+b$. Choose  Set $x\ag[][1] = x_1$, $w\ag[][1] = w_1$ and $y\ag[][1]=y_1=0$.
			\For {$t=1,\ldots,N-1$}
				\begin{align}
					\label{eqnAADMMxmd}
					x\md = &\ (1-\alpha_t)x\ag + \alpha_t x_t
					,\\
					\nonumber
					x\tn = &\ \argmin_{x\in X}\langle\nabla G(x\md), x\rangle  - \chi\theta_t\langle Bw_t - Kx_t  - b, Kx\rangle 
					\\
					&\ + \frac{(1-\chi)\theta_t}{2}\|Bw_t - Kx  - b\|^2 + \langle y_t, Kx\rangle 
					+ \frac{\eta_t}{2}\|x - x_t\|^2	\label{eqnAADMMxtn}
					,\\
					\label{eqnAADMMxag}
					x\ag[1] = &\ (1-\alpha_t)x\ag + \alpha_tx\tn
					,\\
					\label{eqnAADMMwtn}
					w\tn =&\ \argmin_{w\in \cW}F(w) - \langle y_t, Bw\rangle + \frac{\tau_t}{2}\|Bw - Kx\tn - b\|^2,
					\\
					\label{eqnAADMMwag}
					w\ag[1] = &\ (1-\alpha_t)w\ag + \alpha_tw\tn
					,\\
					\label{eqnAADMMytn}
					y\tn =&\ y_t - \rho_t(Bw\tn - Kx\tn - b),
					\\
					\label{eqnAADMMyag}
					y\ag[1] = &\ (1-\alpha_t)y\ag + \alpha_ty\tn
					.
				\end{align}
				\EndFor
			\State Output $z\ag[][N] = (w\ag[][N], x\ag[][N])$.
	\end{algorithmic}
\end{algorithm}

In AADMM, the binary constant $\chi$ in \eqref{eqnAADMMxtn} is either $0$ or $1$, the superscript ``ag'' stands for ``aggregate'', and ``md'' stands for ``middle''. It can be seen that the middle point $x\md$, and the aggregate points $w\ag[1]$, $x\ag[1]$ and $y\ag[1]$ are weighted sums of all the previous iterates $\{x_i\}_{i=1}^{t}$,  $\{w_i\}_{i=1}^{t+1}$, $\{x_i\}_{i=1}^{t+1}$ and $\{y_i\}_{i=1}^{t+1}$, respectively. If the weights $\alpha_t\equiv 1$, then $x\md=x_t$ and the aggregate points are exactly the current iterates $w\tn$, $x\tn$ and $y\tn$. 
In this case, if $\chi=0$, and $\theta_t = \tau_t = \rho_t \equiv \rho$, then AADMM becomes L-ADMM, and if in addition $G$ is simple, then AADMM becomes ADMM. On the other hand, if $\chi=1$, then AADMM becomes LP-ADMM, and if in addition $G$ is simple, AADMM becomes P-ADMM. 

In this work, we will show that if $G$ is non-simple, by properly specifying the parameter $\alpha_t$, we can significantly improve the rate of convergence of Algorithm \ref{algAADMM} in terms of its dependence on $L_G$, with about the same iteration cost. We call the acceleration for $\chi=0$ the accelerated L-ADMM (AL-ADMM), and call that for $\chi=1$ the accelerated LP-ADMM (ALP-ADMM).

Next, we define certain appropriate gap functions.

\subsection{Gap functions}
For any $\tilde z = (\tilde w,\tilde x,\tilde y)\in \cZ$ and $z = (w,x,y)\in \cZ$, we define
\begin{align}
	\label{eqnQ}
	Q(\tilde{w}, \tilde{x}, \tilde{y}; w, x, y):= [G(x) + F(w) -\langle \tilde{y}, Bw-Kx-b\rangle] - [G(\tilde x) + F(\tilde w) -\langle {y}, B\tilde w-K\tilde x-b\rangle].
\end{align}
For simplicity, we use the notation $Q(\tilde{z}; z):=Q(\tilde{w}, \tilde{x}, \tilde{y}; w, x, y)$, and under different situations, we may use notations $Q(\tilde{z}; w, x, y)$ or $Q(\tilde{w}, \tilde{x}, \tilde{y}; z)$ for the same meaning. We can see that $Q(z^*, z)\geq 0$ and $Q(z, z^*)\leq 0$ for all $z\in \cZ$, where $z^*$ is a saddle point of \eqref{eqnSPP}, as defined in Section \ref{secNotation}. For compact sets $W\subset\cW, X\subset\cX, Y\subset\cY$, the duality gap function
\begin{align}
	\label{eqnDGap}
	\sup_{\tw\in W, \tx\in X, \ty\in Y}Q(\tw, \tx, \ty; w, x, y)
\end{align}
measures the accuracy of an approximate solution $(w,x,y)$ to the saddle point problem
	\begin{align*}
		\min_{x\in X, w\in W}\max_{y\in Y}G(x) + F(w) - \langle y, Bw - Kx - b\rangle.
	\end{align*}
However, our problem of interest \eqref{eqnAECCO} has a saddle point formulation \eqref{eqnSPP}, in which the feasible set $(\cW, X, \cY)$ may be unbounded. Recently, a perturbation-based termination criterion is employed by Monteiro and Svaiter \cite{monteiro2010complexity, monteiro2011complexity, monteiro2013iteration} for solving variational inequalities and saddle point problems. This termination criterion is based on the enlargement of a maximal monotone operator, which is first 
introduced in \cite{burachik1997enlargement}. One advantage of using this termination criterion is that its definition does not 
depend on the boundedness of the domain of the operator. We modify this termination criterion and propose a modified version of the gap function in \eqref{eqnDGap}. More specifically, we define
\begin{align}
	\label{eqnGap}
	g_Y(v,z):=\sup_{\tilde y\in Y}Q(w^*, x^*, \tilde y; z) + \langle v, \tilde y\rangle
\end{align}
for any closed set $Y\subseteq \cY$, and for any $z\in Z$ and $v\in Y$. In addition, we denote
\begin{align}
	\label{eqnGap0}
	\bar g_Y(z) := g_Y(0, z) =\sup_{\tilde y\in Y}Q(w^*, x^*, \tilde y; z).
\end{align}
If $Y=\cY$, we will omit the subscript $Y$ and simply use notations $g(v, z)$ and $\bar g(z)$.

\vgap 

In Propositions \ref{proGap} and \ref{proGap0} below, we describe the relationship between the gap functions \eqref{eqnGap}--\eqref{eqnGap0} and the approximate solutions to problems \eqref{eqnAECCO} and \eqref{eqnUCO}.

\begin{pro}
	\label{proGap}
	For any $Y\subseteq\cY$, if $g_Y(Bw-Kx-b, z)\leq\varepsilon<\infty$  and $\|Bw-Kx-b\|\leq\delta$ where $z=(w,x,y)\in\cZ$, then $(w, x)$ is an $(\epsilon, \delta)$-solution of \eqref{eqnAECCO}. In particular, when $Y=\cY$, for any $v$ such that $g(v,z)\leq\varepsilon<\infty$ and $\|v\|\leq\delta$, we always have $v=Bw-Kx-b$.
	\begin{proof}
		By \eqref{eqnQ} and \eqref{eqnGap}, for all $v\in\cY$ and $Y\subseteq \cY$, we have
		\begin{align*}
			\begin{aligned}
			g_Y(v, z) = &\ \sup_{\ty\in Y}[G(x) + F(w) - \langle \ty, \BKb\rangle] - [G(x^*) + F(w^*) ] + \langle v, \ty\rangle
			\\
			= &\ G(x) + F(w) - f^* + \sup_{\ty\in Y}\langle -\ty, \BKb-v\rangle.
			\end{aligned}
		\end{align*}
		From the above we see that if $g_Y(Bw-Kx-b,z) = G(x)+F(w)-f^*\le\epsilon$ and $\|Bw-Kx-b\|\leq\delta$, then $(w,z)$ is an  $(\epsilon, \delta)$-solution. In addition, if $Y=\cY$, we can also see that $g(v,z)=\infty$ if $v\not=\BKb$, hence $g(v,z)<\infty$ implies that $v=Bw-Kx-b$.
	\end{proof}
\end{pro}

\vgap

From Proposition \ref{proGap} we can see that when $Y=\cY$ and $g(v,z)\le\epsilon$, $\|v\|$ is always the feasibility residual of the approximate solution $(w,x)$. Proposition \ref{proGap0} below shows that in some special cases, there exists an approximate solution to problem \eqref{eqnAECCO} that has zero feasibility residual.

\vgap

\begin{pro}
	\label{proGap0}
	Assume that $B$ is an one-to-one linear operator such that $B\cW = \cY$, and $F(\cdot)$ is Lipschitz continuous, then the set $Y:=(B^*)^{-1}\dom F^*$ is bounded. Moreover, if $\bar g_Y(z)\leq \epsilon$, then the pair $(\tw, x)$ is an $\epsilon$-solution of \eqref{eqnAECCO}, where $\tw = (B^*)^{-1}(Kx+b)$. 
	\begin{proof}
		We can see that $\tw$ is well-defined since $B\cW = \cY$. Also, using the fact that $F(\cdot)$ is finite valued, by Corollary 13.3.3 in \cite{rockafellar1970convex} we know that $\dom F^*$ is bounded, hence $Y$ is bounded. In addition, as $B\tw - Kx - b=0$, 
		\begin{align*}
			\bar g_Y(z) = &\ \sup_{\ty\in Y}[G(x) + F(w) - \langle \ty, \BKb\rangle] - [G(x^*) + F(w^*) ] 
			\\
			= &\ G(x) + F(w) - f^* + \sup_{\ty\in Y}\langle -\ty, Bw-B\tilde{w}\rangle
			\\
			= &\ G(x) + F(\tilde w) - f^* + \sup_{\tilde y\in Y}[F(w) - F(\tilde w) - \langle B^*\tilde y, w - \tilde w\rangle].
		\end{align*}
		If $B^*Y\cap\partial F(\tw)\not=\emptyset$, then from the convexity of $F(\cdot)$ we have
		\begin{align*}
			\bar g_Y(z)\geq G(x) + F(\tw) - f^*,
		\end{align*}
		thus $(\tw, x)$ is an $\epsilon$-solution. To finish the proof it suffices to show that $B^*Y\cap\partial F(\tw)\not=\emptyset$. Observing that 
		\begin{align*}
			\sup_{\bar w\in B^*Y}\langle \tw, \bar w\rangle - F^*(\bar w) = \sup_{\bar w\in \dom F^*}\langle \tw, \bar w\rangle - F^*(\bar w) =  \sup_{\bar w\in \cW}\langle \tw, \bar w\rangle - F^*(\bar w),
		\end{align*}
		and using the fact that $Y$ is closed, we can conclude that there exists $B^*\ty\in B^*Y$ such that $B^*\ty$ attains the supremum of the function $\langle \tw, \bar w\rangle - F^*(\bar w)$ with respect to $\bar w$. By Theorem 23.5 in \cite{rockafellar1970convex}, we have $B^*\ty\in\partial F(\tw)$, and hence $\partial F(\tw)\cap B^*Y\not = \emptyset$. 
	\end{proof}
\end{pro}

A direct consequence of the above proposition is that for the UCO problem, if $F(\cdot)$ is Lipschitz continuous and $\bar g_Y(z)\leq \epsilon$, then $(x, Kx)$ is an $\varepsilon$-solution.

\subsection{Main estimations}

In this subsection, we present the main estimates that will be used to prove the rate of convergence for AADMM.

\begin{lem}
	\label{lemQBound}
	Let
		\begin{align}
			\label{eqnGamma}
			\Gamma_t = 
			\begin{cases}
				\Gamma_1 & \text{ when }\alpha_t=1,
				\\ 
				(1-\alpha_t)\Gamma_{t-1} & \text{ when }t>1.
			\end{cases}
		\end{align}
	For all $y\in\cY$, the iterates $\{z\ag\}_{t\geq 1}:=\{(w\ag, x\ag, y\ag)\}_{t\geq 1}$ of Algorithm \ref{algAADMM} satisfy
		\begin{align}
		\label{eqnQBound}
			\begin{aligned}
				&\ \frac{1}{\Gamma_t}Q(w^*, x^*, y;z\ag[1]) - \sum_{i=2}^t\left(\frac{1-\alpha_i}{\Gamma_i} - \frac{1}{\Gamma_{i-1}} \right)Q(w^*, x^*, y;z\ag[][i])
				\\
				\leq &\ \cBx[^*] + \cBy + \cB_t(Bw^*, Bw_{[t+1]}, \theta_{[t]}) - \chi\cB_t(Kx^*, Kx_{[t+1]}, \theta_{[t]})
				\\
				&\  - \sumt\frac{\alpha_i(\tau_i - \theta_i)}{2\Gamma_i}\|Bw\tn[i] - Kx^*-b \|^2  + \sumt\frac{\alpha_i(\tau_i - \theta_i)}{2\Gamma_i}\|K(x\tn[i] - x^*)\|^2 
				\\
				&\ - \sumt \frac{\alpha_i(\tau_i-\rho_i)}{2\Gamma_i\rho_i^2}\|y_i - y\tn[i]\|^2 - \sumt \frac{\alpha_i}{2\Gamma_i}\left({\eta_i} - L_G\alpha_i - \chi\theta_i\|K\|^2 \right)\|x_i - x\tn[i]\|^2.
			\end{aligned}
		\end{align}
			where the term $\cB_t(\cdot,\cdot,\cdot)$ is defined as follows: for any point $v$ and any sequence $v_{[t+1]}$ in any vectorial space $\cV$, and any real valued sequence $\gamma_{[t]}$,
				\begin{align}
					\label{eqnB}
					\cB_t(v, v_{[t+1]}, \gamma_{[t]}): = \sum_{i=1}^t\frac{\alpha_i}{2\Gamma_i}\gamma_i\left(\|v_i - v\|^2 - \|v\tn[i] - v\|^2 \right).
				\end{align}
		\begin{proof}
		To start with, we prove an important property of the function $Q(\cdot,\cdot)$ under Algorithm \ref{algAADMM}. By convexity of $G(\cdot)$ we have
		\begin{align}
			\label{eqnLGinProof}
			\begin{aligned}
				&\ G(x\ag[1]) 
				\leq \  G(x\md)+\langle \nabla G(x\md), x\ag[1]-x\md\rangle + \frac{L_G}{2}\|x\ag[1]-x\md\|^2.
			\end{aligned}
		\end{align}
		Moreover, by equations \eqref{eqnAADMMxmd} and \eqref{eqnAADMMxag}, $x\ag[1] - x\md[1] = \alpha_t(x\tn - x_t)$. Using this observation, equation \eqref{eqnLGinProof} and the convexity of $G(\cdot)$, we have
		\begin{align}
			\label{eqnGRecur}
			\begin{aligned}
				&\ G(x\ag[1]) 
				\leq \  G(x\md) + \langle\nabla G(x\md), x\ag[1]-x\md\rangle + \frac{L_G\alpha_t^2}{2}\|x_{t+1}-x_t\|^2
				\\
				= &\ G(x\md) + (1-\alpha_t)\langle\nabla G(x\md), x\ag - x\md\rangle + \alpha_t\langle\nabla G(x\md), x_{t+1} - x\md\rangle + \frac{L_G\alpha_t^2}{2}\|x_{t+1}-x_t\|^2
				\\
				= &\ (1-\alpha_t)\left[G(x\md)+ \langle\nabla G(x\md), x\ag - x\md\rangle\right]  + \alpha_t\left[G(x\md)+\langle\nabla G(x\md), x_{t+1} - x\md\rangle\right] 
				\\
				&\ + \frac{L_G\alpha_t^2}{2}\|x_{t+1}-x_t\|^2
				\\
				= &\ (1-\alpha_t)\left[G(x\md)+ \langle\nabla G(x\md), x\ag - x\md\rangle\right]  +  \alpha_t\left[G(x\md)+\langle\nabla G(x\md),  x - x\md\rangle\right] 
				\\
				&\ + \alpha_t\langle\nabla G(x\md), x_{t+1} -  x\rangle + \frac{L_G\alpha_t^2}{2}\|x_{t+1}-x_t\|^2
				\\
				\leq&\ (1-\alpha_t)G(x\ag) + \alpha_tG( x) +  \alpha_t\langle\nabla G(x\md), x_{t+1} -  x\rangle + \frac{L_G\alpha_t^2}{2}\|x_{t+1}-x_t\|^2, \ \forall x\in X.
			\end{aligned}
		\end{align}
		By \eqref{eqnADMMwtn}, \eqref{eqnADMMytn}, \eqref{eqnQ}, \eqref{eqnGRecur} and the convexity of $F(\cdot)$, we conclude that		
		\begin{align}
			\label{eqnQRecur}
			\begin{aligned}
			&\ Q(z;z\ag[1]) - (1-\alpha_t)Q(z;z\ag)
			\\
			= &\ [G(x\ag[1]) + F(w\ag[1]) -\langle y, \BKb[{\ag[1]}]\rangle] 
			- [G(x) + F(w) -\langle {y\ag[1]}, \BKb\rangle]
			\\
			&\ - (1-\alpha_t)[G(x\ag) + F(w\ag) -\langle y, \BKb[{\ag}]\rangle]  + (1-\alpha_t)[G(x) + F(w) -\langle {y\ag}, \BKb\rangle]
			\\
			= &\ \left[G(x\ag[1]) - (1-\alpha_t)G(x\ag) - \alpha_tG(x) \right] + \left[F(w\ag[1]) - (1-\alpha_t)F(w\ag) - \alpha_tF(w) \right] 
			\\
			&\ -\alpha_t\langle y, \BKb[\tn]\rangle + \alpha_t\langle y\tn, \BKb\rangle
			\\
			\leq &\ \alpha_t\left\{\vphantom{\frac12}\langle\nabla G(x\md), x_{t+1} -  x\rangle  + \left[F(w\tn) - F(w) \right]+ \frac{L_G\alpha_t}{2}\|x_{t+1}-x_t\|^2 \right.
			\\
			&\ \left.- \langle y, \BKb[\tn]\rangle + \langle y\tn, \BKb\rangle  \frac{}{}\right\}
			.\\
			\end{aligned}
		\end{align}			
		Next, we examine the optimality conditions in  \eqref{eqnAADMMxtn} and \eqref{eqnAADMMwtn}. for all $x\in X$ and $w\in \cW$, we have
			\begin{align*}
				& \langle\nabla G(x\md) + \eta_t(x\tn - x_t), x\tn - x\rangle  - \langle \theta_t (Bw_t - K\tx_t - b) - y_t, K(x\tn - x)\rangle \leq 0, \text{ and}
				\\
				& F(w\tn) - F(w) + \langle\tau_t(\BKb[\tn]) - y_t, B(w\tn - w)\rangle \leq 0,
			\end{align*}
		where 
		\begin{align}
			\label{eqntx}
			\tilde x_t:=\chi x_t + (1-\chi)x\tn.
		\end{align}
			Observing from \eqref{eqnAADMMytn} that $\BKb[\tn] = (y_t - y\tn)/\rho_t$ and $Bw_t - K\tilde{x}_t- b = (y_t - y\tn)/\rho_t - K(\tilde x_t - x\tn) + B(w_t - w\tn)$, the optimality conditions become
			\begin{align*}
				&\ \langle\nabla G(x\md) + {\eta_t}(x\tn - x_t), x\tn - x\rangle + \langle \left(\frac{\theta_t}{\rho_t}-1\right)(y_t - y\tn) - y\tn, -K(x\tn - x)\rangle
				\\
				&\ \ \ + \theta_t\langle K(\tilde x_t - x\tn), K(x\tn - x)\rangle + \theta_t\langle B(w_t - w\tn), -K(x\tn - x)\rangle \leq 0, \text{ and}
				\\
				&\ F(w\tn) - F(w) +\langle \left(\frac{\tau_t}{\rho_t}-1\right)(y_t - y\tn) - y\tn, B(w\tn - w)\rangle\leq 0.
			\end{align*}
			Therefore,
			\begin{align}
				\label{eqnOpt}
				\begin{aligned}
					&\ \langle\nabla G(x\md), x\tn - x\rangle +  F(w\tn) - F(w)  - \langle y, \BKb[\tn]\rangle + \langle y\tn, \BKb\rangle
					\\
					\leq &\ \langle{\eta_t}(x_t - x\tn), x\tn - x\rangle + \langle y\tn - y, \ABb[\tn]\rangle
					\\
					&\ - \langle \thetarho(y_t - y\tn), -K(x\tn - x)\rangle - \langle\taurho(y_t -y\tn), B(w\tn - w)\rangle
					\\
					&\ + \theta_t\langle K(x\tn - \tilde x_t), K(x\tn - x)\rangle +  \theta_t\langle B(w\tn - w_t), -K(x\tn - x)\rangle.
				\end{aligned}
			\end{align}
			Three observations on the right hand side of \eqref{eqnOpt} are in place. Firstly, by \eqref{eqnAADMMytn} we have
			\begin{align}
				\label{eqnRHS1}
				\begin{aligned}
				&\ \langle {\eta_t}(x_t - x\tn), x\tn - x\rangle + \langle y\tn - y, \ABb[\tn]\rangle
				\\
				= &\ {\eta_t}\langle x_t - x\tn, x\tn - x\rangle + \frac{1}{\rho_t}\langle y\tn - y, y_t - y\tn\rangle
				\\
				= &\ \frac{\eta_t}2(\|x_t - x\|^2 - \|x\tn - x\|^2) - \frac{\eta_t}{2}(\|x_t - x\tn\|^2 ) 
				 + \frac{1}{2\rho_t}\left(\|y_t - y\|^2 - \|y\tn - y\|^2 - \|y_t - y\tn\|^2 \right),
				\end{aligned}
			\end{align}
			and secondly, by \eqref{eqnAADMMytn} we can see that 
			\begin{align}
				\label{eqnBwAx}
				B(w\tn - w) = \frac{1}{\rho_t}(y_t - y\tn) + (Kx\tn - Kx) - (\BKb),
			\end{align}
			and
			\begin{align}
				\label{eqnRHS2}
				\begin{aligned}
				 &\ \langle \thetarho(y_t - y\tn), K(x\tn - x)\rangle - \langle\taurho(y_t -y\tn), \frac{1}{\rho_t}(y_t - y\tn) + (Kx\tn - Kx)\rangle
				\\
				= &\ \frac{\tau_t - \theta_t}{\rho_t}\langle y_t - y\tn, -K(x\tn - x)\rangle - \frac{\tau_t - \rho_t}{\rho_t^2}\|y_t -y\tn\|^2 
				\\
				= &\ \frac{\tau_t - \theta_t}2\left[\frac{1}{\rho_t^2}\|y_t - y\tn\|^2 + \|K(x\tn - x)\|^2 - \|\frac{1}{\rho_t}(y_t - y\tn) + K(x\tn - x)\|^2 \right] - \frac{\tau_t - \rho_t}{\rho_t^2}\|y_t -y\tn\|^2
				\\
				= &\ \frac{\tau_t - \theta_t}{2}\left[\frac{1}{\rho_t^2}\|y_t - y\tn\|^2 + \|K(x\tn - x)\|^2 - \|Bw\tn - Kx - b\|^2 \right] - \frac{\tau_t - \rho_t}{\rho_t^2}\|y_t -y\tn\|^2
				.
				\end{aligned}
			\end{align}
			Thirdly, from \eqref{eqntx} we have
			\begin{align}
				\label{eqnRHS3}
				\begin{aligned}
				&\ \theta_t\langle K(x\tn - \tilde x_t), K(x\tn - x)\rangle+ \theta_t\langle B(w\tn - w_t), -K(x\tn - x)\rangle
				\\
				= &\ -\frac{\chi\theta_t}{2}\left(\|K(x_t - x)\|^2 - \|K(x\tn - x)\|^2 - \|K(x_t - x\tn)\|^2 \right)
				\\
				&\ + \frac{\theta_t}{2}\left(\|Bw_t - Kx - b\|^2 - \|Bw\tn - Kx - b\|^2 + \|\ABb[\tn]\|^2 - \|Bw_t - Kx\tn - b\|^2 \right)
				\\
				\leq &\ -\frac{\chi\theta_t}{2}\left(\|K(x_t - x)\|^2 - \|K(x\tn - x)\|^2\right) + \frac{\chi\theta_t\|K\|^2}{2}\|x_t - x\tn\|^2
				\\
				&\ +\frac{\theta_t}{2}\left(\|Bw_t - Kx - b\|^2 - \| Bw\tn - Kx - b\|^2 \right) + \frac{\theta_t}{2\rho_t^2}\|y_t - y\tn\|^2 - \frac{\theta_t}{2}\|Bw_t - Kx\tn - b\|^2,
				\end{aligned}
			\end{align}
		where the last inequality results from the fact that
		\begin{align}
			\label{eqnLKinProof}
			\chi\|K(x_t-x\tn)\|\le \chi\|K\|\|x_t-x\tn\|.
		\end{align}
		Applying \eqref{eqnOpt} -- \eqref{eqnRHS3} to 	\eqref{eqnQRecur}, we have
		\begin{align}
			\label{tmp}
			\begin{aligned}
			&\ \frac{1}{\Gamma_t}Q(z;z\ag[1]) - \frac{1-\alpha_t}{\Gamma_t}Q(z;z\ag)
			\\
			\leq&\ \frac{\alpha_t}{\Gamma_t} \left\{\vphantom{\frac 12} \right.
			\frac{\eta_t}2(\|x_t - x\|^2 - \|x\tn - x\|^2)  + \frac{1}{2\rho_t}(\|y_t - y\|^2 - \|y\tn - y\|^2) - \frac{\tau_t- \rho_t}{2\rho_t^2}\|y_t - y\tn\|^2
			\\
			&\ + \frac{\theta_t}{2}\|Bw_t - Kx - b\|^2 - \frac{\tau_t}{2}\|Bw\tn - Kx - b\|^2 - \frac{\chi\theta_t}{2}(\|K(x_t - x)\|^2 - \|K(x\tn - x)\|^2) 
			\\
			&\ + \langle \taurho(y_t - y\tn), \ABb\rangle + \frac{\tau_t - \theta_t}{2}\|K(x\tn - x)\|^2  -\frac{\theta_t}{2}\|Bw_t - Kx\tn - b\|^2
			\\
			&\ - \frac 12\left({\eta_t} - L_G\alpha_t - \chi\theta_t\|K\|^2\right)\|x_t - x\tn\|^2 
			\left.\vphantom{\frac12}\right\}.
			\end{aligned}
		\end{align}
		Letting $w=w^*$ and $x=x^*$ in the above, observing from \eqref{eqnGamma} that $\Gamma_{t-1} = (1-\alpha_t)/\Gamma_t$, in view of \eqref{eqnB} and applying the above inequality inductively, we conclude \eqref{eqnQBound}.
	\end{proof}
	\end{lem}

	There are two major consequences of Lemma \ref{lemQBound}. If $\alpha_t\equiv 1$ for all $t$, then the left hand side of \eqref{eqnQBound} becomes $\frac{1}{\Gamma_1}\sum_{i=2}^t Q(z; z\ag[][i])$. On the other hand, if $\alpha_t\in [0, 1)$ for all $t$, then in view of \eqref{eqnGamma}, the left hand side of \eqref{eqnQBound} is $Q(z;z\ag[1])/\Gamma_t$. This difference is the main reason why we can accelerate the rate of convergence of AADMM in terms of $L_G$.

	\vgap 
	
	In the next lemma, we provide possible bounds of $\cB(\cdot,\cdot,\cdot)$ in Lemma \ref{lemQBound}.
\begin{lem}
	\label{lemB}
	Suppose that $\cV$ is any vector space and $V\subset\cV$ is any convex set. For any $v\in V$, $v_{[t+1]}\subset \cV$ and $\gamma_{[t]}\subset \R$, we have the following:
	\begin{enumerate}
		\renewcommand{\theenumi}{\alph{enumi})}
		\item If the sequence $\{\alpha_i\gamma_i/\Gamma_i \} $ is decreasing, then
		\begin{align}
			\label{eqnBUB}
			\cB_t(v, v_{[t+1]}, \gamma_{[t]}) \leq \frac{\alpha_1\gamma_1}{2\Gamma_1}\|v_1 - v\|^2 - \frac{\alpha_t\gamma_t}{2\Gamma_t}\|v\tn - v\|^2.
		\end{align}
		\item If the sequence $\{\alpha_i\gamma_i/\Gamma_i \} $ is increasing, $V$ is bounded and $v_{[t+1]}\subset V$, then
		\begin{align}
			\label{eqnBBD}
			\cB_t(v, v_{[t+1]}, \gamma_{[t]}) \leq \frac{\alpha_t\gamma_t}{2\Gamma_t}D_V^2 - \frac{\alpha_t\gamma_t}{2\Gamma_t}\|v\tn - v\|^2.
		\end{align}
	\end{enumerate}
	\begin{proof}
		By \eqref{eqnB} we have
		\begin{align*}
			\cB_t(v, v_{[t+1]}, \gamma_{[t]}) = \frac{\alpha_1\gamma_1}{2\Gamma_1}\|v_1 - v\|^2 - \sumt[t-1]\left(\frac{\alpha_i\gamma_i}{2\Gamma_i} - \frac{\alpha\tn[i]\gamma\tn[i]}{2\Gamma\tn[i]} \right)\|v\tn[i] - v\|^2 - \frac{\alpha_t\gamma_t}{2\Gamma_t}\|v\tn - v\|^2.
		\end{align*}
		If the sequence $\{\alpha_i\gamma_i/\Gamma_i \} $ is decreasing, then the above equation implies \eqref{eqnBUB}. If the sequence $\{\alpha_i\gamma_i/\Gamma_i \} $ is increasing, $V$ is bounded and $v_{[t+1]}\subset V$, then from the above equation we have
		\begin{align*}
			\cB_t(v, v_{[t+1]}, \gamma_{[t]}) \leq &\  \frac{\alpha_1\gamma_1}{2\Gamma_1}D_V^2 - \sumt[t-1]\left(\frac{\alpha_i\gamma_i}{2\Gamma_i} - \frac{\alpha\tn[i]\gamma\tn[i]}{2\Gamma\tn[i]} \right)D_V^2 - \frac{\alpha_t\gamma_t}{2\Gamma_t}\|v\tn - v\|^2
			\\
			= &\ \frac{\alpha_t\gamma_t}{2\Gamma_t}D_V^2 - \frac{\alpha_t\gamma_t}{2\Gamma_t}\|v\tn - v\|^2,
		\end{align*}
		hence \eqref{eqnBBD} holds.
	\end{proof}
\end{lem}

	\subsection{Convergence results on solving UCO problems in bounded domain}
\label{secUCO}

We study UCO problems with bounded feasible sets in this subsection. In particular, throughout this subsection we assume that
\begin{align}
	\label{BD}
	\text{\emph{Both $X$ and $Y:=\dom F^*$ are compact, and $B=I$, $b=0$.}}
\end{align}

It should be noted that the boundedness of $Y$ above is equivalent to the Lipschitz continuity of $F(\cdot)$ (see, e.g, Corollary 13.3.3 in \cite{rockafellar1970convex}).

The following Theorem \ref{thmADMMI} generalizes the convergence properties of ADMM algorithms. Although the convergence analysis of ADMM, L-ADMM and P-ADMM has already been done in several literatures (e.g., \cite{monteiro2013iteration,he2012convergence,chambolle2011first,ouyang2013stochastic}), Theorem \ref{thmADMMI} gives a unified view of the convergence properties of all ADMM algorithms.
\begin{thm}
	\label{thmADMMI}
	In AADMM, if the parameters of are set to $\alpha_t\equiv 1$, $\theta_t\equiv\tau_t\equiv\rho_t\equiv\rho$ and $\eta_t\equiv L_G + \chi\rho\|K\|^2$, then 
	\begin{align}
		\label{eqnADMMI}
		G(x^{t+1}) + F(Kx^{t+1}) - f^* \leq \frac{L_G}{2t}D_{X}^2 + \frac{\chi\rho}{2t}\|K\|^2D_{X}^2 +  \frac{(1-\chi)\rho}{2t} D_{X, K}^2 + \frac{D_Y^2}{2\rho t},
	\end{align}
	where $x^{t+1}:=\ds\frac{1}{t}\sum_{i=2}^{t+1}x_i$. Specially, if $\rho$ is given by
	\begin{align}
		\label{eqnRhoD}
		\rho = \frac{D_Y}{\chi\|K\|D_X + (1-\chi)D_{X,K}},
	\end{align}
	then
	\begin{align}
		\label{eqnADMMISimple}
		G(x^{t+1}) + F(\tw^{t+1}) - f^* \leq \frac{L_GD_X^2}{2t} + \frac{\chi\|K\|D_XD_Y + (1-\chi)D_{X,K}D_Y}{t}.
	\end{align}
	\begin{proof}
		Since $\alpha_t\equiv 1$, By \eqref{eqnAADMMxag}, \eqref{eqnAADMMwag} and \eqref{eqnAADMMyag} we have $x\ag = x_t$, $w\ag = w_t$ and $y\ag = y_t$, and we can see that $\Gamma_t \equiv 1$ satisfies \eqref{eqnGamma} . Applying the parameter settings to RHS of \eqref{eqnQBound} in Lemma \ref{lemQBound}, we have
		\begin{align*}
			\cBx[^*] = &\ \frac{\eta}{2}(\|x_1 - x^*\|^2 - \|x\tn - x^*\|^2) 
			\\
			\le&\ \frac{L_G}{2}D_{x^*}^2 + \frac{\chi\rho}{2}\|K\|^2D_{x^*}^2 - \frac{\chi\rho}{2}\|K\|^2\|x\tn - x^*\|^2
			,\\
			\cBw[^*] = &\ \frac{\rho}{2}(\|w_1 - w^*\|^2 - \|w\tn - w^*\|^2) \leq \frac{\rho D_{w^*}^2}{2} = \frac{\rho D_{x^*, K}^2}{2}
			,\\
			-\chi\cBKx[^*] = &\ -\frac{\chi\rho}{2}(\|Kx_1 - Kx^*\|^2 - \|Kx\tn - Kx^*\|^2) 
			\\
			\leq &\ -\frac{\chi\rho}{2} D_{x^*, K}^2 + \frac{\chi\rho}{2}\|K\|^2\|x\tn - x^*\|^2,
			\\
			\cBy \leq&\ \frac{1}{2\rho}(\|y_1 - y\|^2 - \|y\tn - y\|^2) \leq \frac{D_Y^2}{2\rho }, \ \forall y\in Y.
		\end{align*}
		Therefore, by Lemma \ref{lemQBound} we have
		\begin{align*}
			&\ \sum_{i=2}^{t+1}Q(w^*, x^*, y;z_i)
			\le \frac{L_G}{2}D_{x^*}^2 + \frac{\chi\rho}{2}\|K\|^2D_{x^*}^2 +  \frac{(1-\chi)\rho}{2} D_{x^*, K}^2 + \frac{D_Y^2}{2\rho }
			\\
			\le &\ \frac{L_G}{2}D_{X}^2 + \frac{\chi\rho}{2}\|K\|^2D_{X}^2 +  \frac{(1-\chi)\rho}{2} D_{X, K}^2 + \frac{D_Y^2}{2\rho },\ \forall y\in Y.
		\end{align*}
		Furthermore, noticing that for all $y\in Y$, by the convexity of $Q(x^*, w^*, y; \cdot)$,
		\begin{align*}
			Q(w^*, x^*, y; z^{t+1}) \leq \frac{1}{t}\sum_{i=2}^{t+1}Q(w^*, x^*, y;z_i), \text{ where }z^{t+1}:=\frac{1}{t}\sum_{i=2}^{t+1}z_i.
		\end{align*}
		Applying the two inequalities above to \eqref{eqnGap0} and Proposition \ref{proGap0}, we conclude \eqref{eqnADMMI}, and \eqref{eqnADMMISimple} follows immediately.
	\end{proof}

	\end{thm}

\vgap 

Although AADMM unifies all ADMM algorithms, what makes it most special is the variable weighting sequence $\{\alpha_t \}_{t\geq 1}$ (rather than $\alpha_t=1$) that accelerates its convergence rate with respect to its dependence on $L_G$, as shown in 
Theorem \ref{thmAPD} below. 

\begin{thm}
	\label{thmAPD}
	In AADMM, if the parameters are set to
	\begin{align}
		\label{condAPD}
		\alpha_{t} = \frac{2}{t+1}, 
		\tau_t = \rho_t \equiv \rho,\ \theta_t = \frac{(t-1)\rho}{t},\ \text{and }\eta_t = \frac{2L_G + \chi\rho t\|K\|^2}{t},
	\end{align}
	then 
		\begin{align}
			\label{eqnAADMMIrho}
			G(x\ag[1]) + F(Kx\ag[1]) - f^*
			\leq &\ \frac{2L_GD_X^2}{t(t+1)} + \frac{1}{t+1}\left[{\chi{\rho}\|K\|^2D_X^2}  + (1-\chi){\rho}D_{X, K}^2+ \frac{D_Y^2}{{\rho}}\right].
		\end{align}	
	In particular, if $\rho$ is given by \eqref{eqnRhoD}, then
		\begin{align}
			\label{eqnAADMMI}
			&\ G(x\ag[1]) + F(Kx\ag[1]) - f^*
			\leq \ \frac{2L_GD_X^2}{t(t+1)} + \frac{2}{t+1}\left[{\chi\|K\|D_XD_Y}  + (1-\chi)D_{X, K}D_Y\right].
		\end{align}

\begin{proof}
	It is clear that
	\begin{align}
		\label{eqnAlphaGammaSimple}
		\alpha_t=\frac{2}{t+1}\text{ and }\Gamma_t=\frac{2}{t(t+1)}\text{ satisfies \eqref{eqnGamma}, and }\frac{\alpha_t}{\Gamma_t}=t.
	\end{align}
	By the parameter setting \eqref{condAPD} and the definition of $\cB(\cdot,\cdot,\cdot)$ in \eqref{eqnB}, it is easy to calculate that
		\begin{align*}
			&\ {\eta_t} - L_G\alpha_t - \chi\theta_t\|K\|^2  \ge 0,\ \tau_t\geq \theta_t,
			\\
			&\ \cB_t(w^*, w_{[t+1]}, \theta_t) - \sumt\frac{\alpha_i(\tau_i - \theta_i)}{2\Gamma_i}\|w\tn[i] - w^*\|^2 
			= -\frac{\rho t}{2}\|w\tn - w^*\|^2 \leq 0, 
			\\
			&\ -\chi\cB_t(Kx^*, Kx_{[t+1]}, \theta_t) + \sumt\frac{\alpha_i(\tau_i - \theta_i)}{2\Gamma_i}\|Kx\tn[i] - Kx^*\|^2 
			\\
			= &\ \frac{\chi\rho t}{2}\|Kx\tn - Kx^*\|^2 + \frac{(1-\chi)\rho}{2}\sumt\|Kx\tn[i] - Kx^*\|^2
			\leq \frac{\chi\rho t}{2}\|K\|^2\|x\tn - x^*\|^2 + \frac{(1-\chi)\rho t}{2}D_{X,K}^2.
		\end{align*}
		Moreover, by \eqref{eqnAADMMwtn}, \eqref{eqnAADMMytn} and Moreau's decomposition theorem (see, e.g., \cite{moreau1962decomposition,combettes2005signal,esser2010general}), we have
		\begin{align}
			\label{eqnMoreau}
			\begin{aligned}
			y\tn =&\ y_t - \rho(w\tn - Kx\tn - b)
			\\
			= &\ (y_t + \rho Kx\tn + \rho b) - \rho\argmin_{w\in \cW}F(w) + \frac{\rho}{2}\|w - \frac{1}{\rho}(y_t - Kx\tn - b)\|^2
			\\
			= &\ \argmin_{y\in \cY} F^*(y) + \frac{1}{2\rho}\|y - \frac{1}{\rho}(y_t - Kx\tn - b)\|^2,
			\end{aligned}
		\end{align}
		which  implies that $y_{[t+1]}\subset Y$. Using this observation together with the fact that $\alpha_t/(\Gamma_t\rho_t) = t/\rho$, and applying \eqref{eqnBBD} in Lemma \ref{lemB}, we obtain
		\begin{align*}
			\cBy \leq &\ \frac{t}{2{\rho}}D_Y^2, \ \forall y\in Y.
		\end{align*}
		Finally, noting that $\alpha_t\eta_t/\Gamma_t = 2L_G + \chi\rho t\|K\|^2$, by \eqref{eqnBBD} in Lemma \ref{lemB} we have
		\begin{align*}
			&\ \cBx[^*] \leq \frac{\alpha_t\eta_t}{2\Gamma_t}D_X^2 - \frac{\alpha_t\eta_t}{2\Gamma_t}\|x\tn - x^*\|^2 
			\le L_GD_X^2 + \frac{\chi\rho t}{2}\|K\|^2D_X^2 - \frac{\chi\rho t}{2}\|K\|^2\|x\tn - x^*\|^2.
		\end{align*}
		Applying all above inequalities to \eqref{eqnQBound} in Lemma \ref{lemQBound}, we have
		\begin{align*}
			&\ \frac{1}{\Gamma_t}Q(w^*, x^*, y; z\ag[1])
			\leq  L_GD_X^2 + \frac{\chi\rho t}{2}\|K\|^2D_X^2 + \frac{(1-\chi)\rho t}{2}D_{X,K}^2 + \frac{t}{2\rho }D_Y^2,\ \forall y\in Y.
		\end{align*}
		Using \eqref{eqnAlphaGammaSimple} and applying Proposition \ref{proGap0}, we conclude \eqref{eqnAADMMIrho}, and \eqref{eqnAADMMI} comes from \eqref{eqnRhoD} and \eqref{eqnAADMMIrho}.
	\end{proof}

	\end{thm}

\vgap

In view of Theorems \ref{thmADMMI} and \ref{thmAPD}, several remarks on the AADMM algorithms are in place. 
Firstly, Theorem \ref{thmAPD} provides an example of choosing stepsizes in AL-ADMM and ALP-ADMM, that leads to better convergence properties w.r.t the dependence on $L_G$ than L-ADMM and LP-ADMM respectively. In particular, AL-ADMM and ALP-ADMM allow $L_G$ to be as large as $\Omega(N)$ without affecting the rate of convergence (up to a constant factor). The comparison of these AADMM algorithms in terms of their rates of convergence is shown in Table \ref{tabRateI}. 
Secondly, ALP-ADMM has the same rate of convergence as Nesterov's smoothing scheme \cite{nesterov2005smooth}, and achieves optimal rate of convergence \eqref{eqnUCOOptRate}. Moreover, we can see from \eqref{eqnMoreau} that the APD method in \cite{chenoptimal} is equivalent to ALP-ADMM. Nonetheless, AL-ADMM has better constant in the estimation of rate of convergence than both ALP-ADMM and Nesterov's smoothing scheme, since $D_{X,K}\leq \|K\|D_X$. However, the computational time for solving problem \eqref{eqnAADMMxtn} with $\chi=0$ is usually higher than that for $\chi=1$, hence AL-ADMM has higher iteration cost than that of ALP-ADMM. The trade-off between better rate constants and cheaper iteration costs has to be considered in practice.
Thirdly, while Theorem \ref{thmADMMI} describes only the ergodic convergence of the ADMM algorithms, Theorem \ref{thmAPD} describes the convergence of aggregate sequences $\{z\ag[1]\}_{t\ge 1}$, which are exactly the outputs of the accelerated schemes.
Finally, in ADMM methods we have $\tau_t=\rho_t=\theta_t$, while in Theorem \ref{thmAPD} we only have $\tau_t=\rho_t$, although $\theta_t\to\rho_t$ when $t\to\infty$. In fact, if the total number of iterations is given, it is possible to choose a set of equal stepsize parameters, as described by Theorem \ref{thmAADMMIN} below.

\begin{thm}
	\label{thmAADMMIN}
	In AADMM, if the total number of iterations $N$ is chosen, and the parameters are set to
	\begin{align*}
		\alpha_t=\frac{2}{t+1},\ \theta_t = \tau_t = \rho_t = \frac{\rho N}{t},\text{ and } \eta_t = \frac{2L_G + \chi\rho N\|K\|^2}{t},
	\end{align*}
	where $\rho$ is given by \eqref{eqnRhoD}, then
	\begin{align}
		\label{eqnAADMMINSimple}
		G(x\ag[][N]) + F(Kx\ag[][N]) - f^* \leq \frac{2L_GD_X^2}{N(N-1)} + \frac{2}{N-1}\left[\chi\|K\|D_XD_Y + (1-\chi)D_{X,K}D_Y\right].
	\end{align}
	\begin{proof}
		Using equation \eqref{eqnAlphaGammaSimple} as well as the definition of $\cB(\cdot, \cdot, \cdot)$ in \eqref{eqnB}, it is easy to calculate that
		\begin{align*}
			\eta_t - L_G\alpha_t - \chi\theta_t\|K\|^2\geq&\  0,
			\\
			\cBx[^*] = &\ \frac{2L_G+\chi\rho N\|K\|^2}{2}(\|x_1 - x^*\|^2 -\|x\tn - x^*\|^2)
			\\
			\leq &\ \frac{2L_G+\chi\rho N\|K\|^2}{2}(D_{x^*}^2 -\|x\tn - x^*\|^2)
			,\\
			\cBw[^*] = &\ \frac{\rho N}{2}(\|w_1 - w^*\|^2 - \|w\tn - w^*\|^2) \leq \frac{\rho N D_{w^*}^2}{2} = \frac{\rho N D_{x^*, K}^2}{2}
			,\text{ and }\\
			-\chi\cBKx[^*] = &\ -\frac{\chi\rho N}{2}(\|Kx_1 - Kx^*\|^2 - \|Kx\tn - Kx^*\|^2)
			\\
			\leq &\ -\frac{\chi\rho N}{2}D_{x^*,K}^2 + \frac{\chi\rho N\|K\|^2}{2}\|x\tn - x^*\|^2.
		\end{align*}
		On the other hand, noting that $\alpha_t/(\Gamma_t\rho_t) = t^2/(\rho N)$, by \eqref{eqnBBD} in Lemma \ref{lemB} we have
		\begin{align*}
			\cBy \leq \frac{t^2}{2\rho N}(D_Y^2 - \|y\tn - y^*\|^2) \leq \frac{t^2D_Y^2}{2\rho N} \le \frac{N}{2\rho}D_Y^2, \ \forall y\in Y,\forall t\leq N.
		\end{align*}
		Applying all the above inequalities to \eqref{eqnQBound} in Lemma \ref{lemQBound}, we conclude
\begin{align*}
			&\ \frac{1}{\Gamma_t} Q(w^*, x^*, y;z\ag[1])
			\le L_GD_{x^*}^2 + \frac{\chi\rho N}{2}\|K\|^2D_{x^*}^2 + \frac{(1-\chi)\rho N}{2}D_{x^*,K}^2 + \frac{N}{2\rho}D_Y^2,
			\\
			\le &\ L_GD_{X}^2 + \frac{\chi\rho N}{2}\|K\|^2D_{X}^2 + \frac{(1-\chi)\rho N}{2}D_{X,K}^2 + \frac{N}{2\rho}D_Y^2.
\end{align*}
		Setting $t=N-1$, and applying \eqref{eqnAlphaGammaSimple}, \eqref{eqnRhoD} and the above inequality to Proposition \ref{proGap0}, we obtain \eqref{eqnAADMMINSimple}.
	\end{proof}
	\end{thm}

\begin{table}[h]
	\caption{\label{tabRateI} Rates of convergence of instances of AADMM for solving UCO with bounded feasible set}
	\ \\
	\centering
	\renewcommand*{\arraystretch}{2.5}
	\begin{tabular}{cll}
		\hline
		& No preconditioning ($\chi=0$) & Preconditioned ($\chi=1$)
		\\\hline
		ADMM & $\ds\cO\left(\frac{D_{X, K}D_Y}{t} \right)$ & $\ds\cO\left(\frac{\|K\|D_XD_Y}{t} \right)$
		\\
		Linearized ADMM & $\ds\cO\left(\frac{L_GD_X^2}{t} + \frac{D_{X, K}D_Y}{t} \right)$ & $\ds\cO\left(\frac{L_GD_X^2}{t} + \frac{\|K\|D_XD_Y}{t} \right)$
		\\
		Accelerated & $\ds\cO\left(\frac{L_GD_X^2}{t^2} + \frac{D_{X, K}D_Y}{t} \right)$ & $\ds\cO\left(\frac{L_GD_X^2}{t^2} + \frac{\|K\|D_XD_Y}{t} \right)$
		\\\hline
	\end{tabular}
\end{table}

\subsection{Convergence results on solving AECCO problems}
\label{secAECCO}
In this section, we study the rate of convergence of AADMM for solving general AECCO problems without boundedness assumption for either $X$ or $Y$, in terms of both primal and feasibility residuals. We start with the convergence analysis of ADMM algorithms as a special case of AADMM where $\alpha_t=1$, $\theta_t=\tau_t=\rho_t=\rho$.
\begin{thm}
	\label{thmADMMII}
	In AADMM, if $\alpha_t\equiv 1$, $\theta_t\equiv\tau_t\equiv\rho_t\equiv\rho$ and $\eta_t\equiv\eta\geq L_G + \chi\rho\|K\|^2$, then
	\begin{align}	
		\label{eqnADMMIIprimal}
		G(x^{t+1}) + F(w^{t+1}) - f^* \leq &\ \frac{1}{2t}\left(\eta D_{x^*}^2 + \rho (1 - \chi) D_{x^*, K}^2\right)
	\end{align}
	and
	\begin{align}
		\label{eqnADMMIIfeas}
		\|Bw^{t+1} - Kx^{t+1} - b\|^2 \leq &\ \frac{2}{t^2}\left(\frac{2D_{y^*}^2}{\rho^2} + \frac{\eta D_{x^*}^2}{\rho} + (1 - \chi) D_{x^*, K}^2 \right),
	\end{align}
	where $x^{t+1}:=\ds\frac{1}{t}\sum_{i=2}^{t+1}x_i$ and $w^{t+1}:=\ds\frac{1}{t}\ds\sum_{i=2}^{t+1}w_i$. Specially, if $\rho = 1$ and $\eta = L_G + \chi\|K\|^2$, then
		\begin{align}
			\label{eqnADMMIIprimSimple}
			G(x^{t+1}) + F(w^{t+1}) - f^* \leq &\ \frac{1}{2t}( L_GD_{x^*}^2 + \chi \|K\|^2D_{x^*}^2 + (1-\chi)D_{x^*, K}^2)
		\end{align}
		and
		\begin{align}
			\label{eqnADMMIIfeasSimple}
			\|Bw^{t+1} - Kx^{t+1} - b\| \leq &\ \frac{2\sqrt{L_G}D_{x^*}}{t} + \frac{\chi\sqrt{2}\|K\|D_{x^*}}{t} + \frac{(1-\chi)\sqrt{2}D_{x^*,K}}{t} + \frac{2D_{y^*}}{t}.
		\end{align}
\begin{proof}
	Similar as the proof of Theorem \ref{thmADMMI}, we have
	\begin{align}
		\nonumber
		&\ Q(w^*, x^*, y;z^{t+1}) 
		\\
		\label{tmp1}
		\leq&\ \frac{1}{2t}\left[L_GD_{x^*}^2 + \chi\rho\|K\|^2D_{x^*}^2 +  (1-\chi)\rho D_{x^*, K}^2 + \frac{1}{\rho}(\|y_1 - y\|^2 - \|y\tn - y\|^2)\right]
		\\
		\label{tmp2}
		\leq&\ \frac{1}{2t}\left[L_GD_{x^*}^2 + \chi\rho\|K\|^2D_{x^*}^2 +  (1-\chi)\rho D_{x^*, K}^2 \right] - \langle\frac{1}{\rho t}(y_1 - y\tn), y\rangle,
	\end{align}
	where $z^{t+1} = \sum_{t=2}^{t+1}z_i$. 
	Noting that $Q(z^*, z^{t+1}) \geq 0$, by \eqref{tmp1} we have
	\[
		\|y\tn - y^*\|^2 \leq \rho L_GD_{x^*}^2 + \chi\rho^2\|K\|^2D_{x^*}^2 +  (1-\chi)\rho^2 D_{x^*, K}^2 + D_{y^*}^2,
	\]
	hence if we let $v\tn = (y_1 - y\tn)/(\rho t)$, then we have
	\begin{align*}
		&\ \|v\tn\|^2 \leq \frac{2}{\rho^2t^2}(\|y_1 - y^*\|^2 + \|y\tn - y^*\|^2 ) 
		\\
		\leq&\ \frac{2}{t^2}(\frac{L_G D_{x^*}^2}{\rho} +\chi\|K\|^2D_{x^*}^2 +  (1-\chi) D_{x^*, K}^2 + \frac{2}{\rho^2}D_{y^*}^2 ).
	\end{align*}
	Furthermore, by \eqref{tmp2} we have
	\begin{align*}
		g(v\tn, z^{t+1}) \leq \frac{1}{2t}\left[L_GD_{x^*}^2 + \chi\rho\|K\|^2D_{x^*}^2 +  (1-\chi)\rho D_{x^*, K}^2 \right].
	\end{align*}
	Applying the two inequalities above to Proposition \ref{proGap} we obtain \eqref{eqnADMMIIprimal} and \eqref{eqnADMMIIfeas}. The results in \eqref{eqnAADMMIIprimSimple} and \eqref{eqnAADMMIIfeasSimple} then follows immediately.
\end{proof}

\end{thm}

\vgap

From Theorem \ref{thmADMMII} we see that the for ADMM algorithms, the rate of convergence of both primal and feasibility residuals are of order $\cO(1/t)$.
The detailed rate of convergence of each algorithm is listed in Tables \ref{tabRateIIprim} and \ref{tabRateIIfeas}. We observe that a larger value of $\rho$ will increase the right side of \eqref{eqnADMMIIprimal}, but decrease that of \eqref{eqnADMMIIfeas}. Hence, an ``optimal" selection of $\rho$ will be determined by considering both primal and feasibility residuals together. For the sake of simplicity, we set $\rho=1$.

\vgap

In Theorem \ref{thmAADMMII} below, we show that there exists a weighting sequence $\{\alpha_t\}_{t\geq 1}$ that improves the rate of convergence of Algorithm \ref{algAADMM} in terms of its dependence on $L_G$.

\begin{thm}
	\label{thmAADMMII}
	In AADMM, if the total number of iterations is set to $N$, and the parameters are set to
	\begin{align}
		\label{condAADMMII}
		\alpha_t = \frac{2}{t+1},\ \theta_t = \tau_t = \frac{N}{t},\ \rho_t = \frac{t}{N},\ \text{and }
		\eta_t = \frac{2L_G + \chi N\|K\|^2}{t},
	\end{align}
	then 
	\begin{align}
		\label{eqnAADMMIIprimSimple}&
\begin{aligned}
		G(x\ag[][N]) + F(w\ag[][N]) - f^*
		\leq &\ \frac{2L_GD_{x^*}^2}{N(N-1)} + \frac{1}{2(N-1)}\left[\chi\|K\|^2D_{x^*}^2 + (1-\chi)D_{x^*, K}^2\right],
\end{aligned}
	\end{align}
	and
	\begin{align}
		\label{eqnAADMMIIfeasSimple}&
\begin{aligned}
		\|\BKb[{\ag[][N]}]\|
		\leq &\  \frac{4\sqrt{L_G}D_{x^*}}{(N-1)\sqrt{N}} +  \frac{2\sqrt{2}\chi\|K\|D_{x^*}}{N-1} + \frac{2\sqrt{2}(1-\chi)D_{x^*, K}}{N-1} + \frac{4D_{y^*}}{N-1}.
\end{aligned}
	\end{align}
\begin{proof}
	Using equations \eqref{condAADMMII}, \eqref{eqnAlphaGammaSimple} and \eqref{eqnB}, we can calculate that
		\begin{align*}
			\eta_t - L_G\alpha_t - \chi\theta_t\|K\|^2\ge&\ 0,\ \tau_t\geq \rho_t \text{ for all }t\le N,
			\\
			\cBx[^*] =&\ \frac{2L_G+\chi N\|K\|^2}{2}(D_{x^*}^2 - \|x\tn - x^*\|^2),
			\\
			\cBy = &\ \frac{N}{2}(\|y_1 - y\|^2 - \|y\tn - y\|^2), \ \forall y\in \cY
			,\\
			\cB_t(Bw^*, Bw_{[t+1]}, \theta_{[t]}) = &\ \frac{ N}{2}(\|Bw_1 - Bw^*\|^2 - \|B\tn - Bw^*\|^2) \leq \frac{ N}{2}D_{w^*, B}^2 = \frac{ N}{2}D_{x^*,K}^2,
			\\
			-\chi\cB_t(Kx^*, Kx_{[t+1]}, \theta_{[t]}) =&\ -\frac{\chi N}{2}(\|Kx_1 - Kx^*\|^2 - \|Kx\tn - Kx^*\|^2)
			,\\
			\leq &\ -\frac{\chi N}{2}(D_{x^*, K} - \|K\|^2\|x\tn - x^*\|^2)
			.
		\end{align*}
		Applying all the above calculations to \eqref{eqnQBound} in Lemma \ref{lemQBound}, we have
		\begin{align*}
			&\ \frac{1}{\Gamma_t}Q(w^*, x^*, y; z\ag[1]) 
			\\
			\leq&\ L_GD_{x^*}^2 + \frac{\chi N}{2}\|K\|^2D_{x^*}^2+\frac{(1-\chi) N}{2}D_{x^*,K}^2 + \frac{N}{2}(\|y_1 - y\|^2 - \|y\tn - y\|^2), \ \forall y\in \cY.
		\end{align*}
		Two consequences to the above estimation can be derived.
		Firstly, since $Q(z^*;z\ag[1])\geq 0$, we have
		\begin{align*}
			\|y\tn - y^*\|^2 \leq \frac{2 L_G}{N}D_{x^*}^2 + \chi\|K\|^2D_{x^*}^2 + (1-\chi) D_{x^*,K}^2 + D_{y^*}^2,
		\end{align*}
		and
		\begin{align*}
			&\ \|y_1 - y\tn\|^2 \leq 2(\|y_1 - y^*\|^2 + \|y\tn - y^*\|^2) \leq \frac{4 L_G}{N}D_{x^*}^2 + 2\chi\|K\|^2D_{x^*}^2 + 2(1-\chi) D_{x^*,K}^2 + 4D_{y^*}^2.
		\end{align*}
		Secondly, since $\|y_1 - y\|^2 - \|y\tn - y\|^2 = \|y_1\|^2 - \|y\tn\|^2 - 2\langle y_1 - y\tn, y\rangle \leq  - 2\langle y_1 - y\tn, y\rangle$,
		\begin{align*}
			\frac{1}{\Gamma_t}Q(w^*, x^*, y; z\ag[1]) + N\langle y_1 - y\tn, y\rangle \leq L_GD_{x^*}^2 + \frac{\chi N}{2}\|K\|^2D_{x^*}^2 + \frac{(1-\chi) N}{2}D_{x^*, K}^2, \ \forall y\in \cY.
		\end{align*}
		Letting $t=N-1$ and $v_N:=2(y_1 - y\tn)/(N-1)$, and applying \eqref{eqnAlphaGammaSimple} and the two above inequalities to Proposition \ref{proGap}, we obtain \eqref{eqnAADMMIIprimSimple} and \eqref{eqnAADMMIIfeasSimple}.
\end{proof}
\end{thm}

\vgap 

Comparing \eqref{eqnADMMIIprimSimple} and \eqref{eqnADMMIIfeasSimple} with \eqref{eqnAADMMIIprimSimple} and \eqref{eqnAADMMIIfeasSimple} respectively, AL-ADMM and ALP-ADMM are better than both L-ADMM and LP-ADMM respectively, in terms of their rates of convergence of both primal and feasibility residuals. The rates of convergence of AADMM algorithms are outlined in Tables \ref{tabRateIIprim} and \ref{tabRateIIfeas}.

\begin{table}[t]
	\caption{\label{tabRateIIprim} Rates of convergence of the primal residuals of AADMM instances for solving general AECCO}
	\ \\
	\centering
	\renewcommand*{\arraystretch}{2.5}
	\begin{tabular}{cll}
		\hline
		& No preconditioning ($\chi=0$) & Preconditioned ($\chi=1$)
		\\\hline
		ADMM & $\ds\cO\left(\frac{D_{x^*, K}^2}{N} \right)$ & $\ds\cO\left(\frac{\|K\|D_{x^*}^2}{N} \right)$
		\\
		Linearized ADMM & $\ds\cO\left(\frac{L_GD_{x^*}^2 + D_{x^*, K}^2}{N} \right)$ & $\ds\cO\left(\frac{L_GD_{x^*}^2 + \|K\|D_{x^*}^2}{N} \right)$
		\\
		Accelerated & $\ds\cO\left(\frac{L_GD_{x^*}^2}{N^2} + \frac{D_{x^*, K}^2}{N} \right)$ & $\ds\cO\left(\frac{L_GD_{x^*}^2}{N^2} + \frac{\|K\|D_{x^*}^2}{N} \right)$
		\\\hline
	\end{tabular}
\end{table}

\begin{table}
	\caption{\label{tabRateIIfeas} Rates of convergence of the feasibility residuals of AADMM instances for solving general AECCO}
	\ \\
	\centering
	\renewcommand*{\arraystretch}{2.5}
	\begin{tabular}{cll}
		\hline
		& No preconditioning ($\chi=0$) & Preconditioned ($\chi=1$)
		\\\hline
		ADMM & $\ds\cO\left(\frac{D_{x^*, K} + D_{y^*}}{N} \right)$ & $\ds\cO\left(\frac{\|K\|D_{x^*} + D_{y^*}}{N} \right)$
		\\
		Linearized ADMM & $\ds\cO\left(\frac{\sqrt{L_G}D_{x^*} + D_{x^*, K} + D_{y^*}}{N} \right)$ & $\ds\cO\left(\frac{\sqrt{L_G}D_{x^*} + \|K\|D_{x^*} + D_{y^*}}{N} \right)$
		\\
		Accelerated  & $\ds\cO\left(\frac{\sqrt{L_G}D_{x^*}}{N^{3/2}} + \frac{D_{x^*, K}+ D_{y^*}}{N} \right)$ & $\ds\cO\left(\frac{\sqrt{L_G}D_{x^*}}{N^{3/2}} + \frac{\|K\|D_{x^*}+ D_{y^*}}{N} \right)$
		\\\hline
	\end{tabular}
\end{table}

\subsection{A simple backtracking scheme}
\label{secLineSearch}

We have discussed the rate of convergence of Algorithm \ref{algAADMM}, with the assumption that both $L_G$ and $\|K\|$ are given. In practice, we may need backtracking techniques to estimate both constants. In this subsection, we propose a simple backtracking technique for AL-ADMM and ALP-ADMM.

From the proof of Lemma \ref{lemQBound}, we can see that if $L_G$ and $\|K\|$ in \eqref{eqnLGinProof} and \eqref{eqnLKinProof} are replaced by $L_t$ and $M_t$ respectively, i.e.,
\begin{align}
	\label{eqnLt}
	&\ G(x\ag[1][t])\le G(x\md)+\langle \nabla G(x\md), x\ag[1]-x\md\rangle + \frac{L_t}{2}\|x\ag[1]-x\md\|^2 \text{ and }
	\\
	\label{eqnMt}
	&\ \chi \|K(x_t - x\tn)\|\le \chi M_t\|x_t - x\tn\|,
\end{align}
then Lemma \ref{lemQBound} still holds. On the other hand, to prove Theorems \ref{thmADMMI} through \ref{thmAADMMII}, in addition to Lemma \ref{lemQBound}, we require monotonicity of the sequences $\alpha_{[t]}\eta_{[t]}/\Gamma_{[t]}$, $\alpha_{[t]}\tau_{[t]}/\Gamma_{[t]}$, $\alpha_{[t]}\theta_{[t]}/\Gamma_t$ and $\alpha_t/(\Gamma_t\rho_t)$, and
\begin{align}
	\label{eqnetat}
	&\ \eta_t - L_t\alpha_t - \chi\theta_tM_t^2\ge 0,
\end{align}
The monotonicity of these sequences is also used in Lemma \ref{lemB}, which helps to prove the boundedness of distances $\cB(\cdot,\cdot,\cdot)$ at the RHS of \eqref{eqnQBound} in Lemma \ref{lemQBound}. From these observations, we can simply use the following choice of parameters:
\begin{align*}
	\theta_t = \tau_t = \frac{\nu_t \alpha_t}{\Gamma_t},\ \rho_t = \frac{\alpha_t}{\nu_t\Gamma_t},\ \eta_t = L_t\alpha_t + \chi\theta_tM_t^2,
\end{align*}
where we assume that $\nu_{[t]}$, $M_{[t]}$ are both monotone. It should be noted that the monotonicity of $\alpha_{[t]}\eta_{[t]}/\Gamma_{[t]}$ relies on $\{L_t\alpha_t^2/\Gamma_t \}_{t\ge 1}$, which is trivial if we simply set $L_t\alpha_t^2 = \Gamma_t$. In addition, in view of the RHS of \eqref{eqnQBound}, we require $\tau_t\ge \rho_t$, i.e., $\nu_t\geq \alpha_t/\Gamma_t$. We summarize all the discussions above to a simple backtracking procedure below.

{	
	\renewcommand\thealgorithm{1}
	\addtocounter{algorithm}{-1}
	\floatname{algorithm}{Procedure}
	\begin{algorithm}[H]
	\caption{\label{funBT}Backtracking procedure for AL-ADMM and ALP-ADMM at the $t$-th iteration}
	\begin{algorithmic}[1]
		\makeatletter
		\setcounter{ALG@line}{-1}
		\makeatother
		\Procedure{Backtracking}{$L_{t-1}$, $M_{t-1}$, $\Gamma_{t-1}$, $\nu_{t-1}$, $x_t$, $x\ag$, $L_{min}$}
			\State 
			\label{codeGoto}
			$L_t\leftarrow \max\{L_{min}, L_{t-1}/2\}$,  $M_t=M_{t-1}$ and $v_t = v_{t-1}$.
			\Comment Initialization
			\State
			\label{codeBacktracking}
			Estimate ${\alpha_t}\in [0,1]$ by solving the quadratic equation
				\begin{align}
					\label{eqnAlphaL}
					L_t\alpha_t^2 = \Gamma_{t-1}(1-{\alpha}_t),
				\end{align}
			\indent and set ${\Gamma}_t \leftarrow \Gamma_{t-1}(1-{\alpha}_t)$, ${\nu}_t = \max\{\nu_{t-1}, {\alpha}_t/{\Gamma}_t\}$.
			\State
			Choose stepsize parameters as
			\begin{align}
				\label{eqnStepL}
				{\theta}_t={\tau}_t=\frac{\rho{\nu}_t{\Gamma}_t}{{\alpha}_t},\ {\rho}_t=\frac{\rho{\alpha}_t}{{\Gamma}_t{\nu}_t}, \text{ and }{\eta}_t = \frac{{\Gamma}_t}{{\alpha}_t} + \chi{\theta}_t {M}_t^2,
			\end{align}
			\indent and calculate iterates \eqref{eqnAADMMxmd} -- \eqref{eqnAADMMxag}.
			\If 
			{$G(x\ag[1]) - G(x\md) - \langle \nabla G(x\md), x\ag[1] - x\md\rangle > \frac{{L}_t}{2}\|x\ag[1] - x\md\|^2$}
			\label{codeLt}
			\State
				Set ${L}_t\leftarrow 2{L}_t$. Go to \ref{codeBacktracking}.
				\Comment  Backtracking $L_G$
		    \ElsIf 
		    {$\chi\|Kx\tn - Kx_t\|> \chi {M}_t\|x\tn - x_t\|$}
		    \label{codeMt}
		    \State
		    	Set ${M}_t\leftarrow 2{M}_t$. Go to \ref{codeBacktracking}.
		    	\Comment  Backtracking $\|K\|$
		    \EndIf\label{codeBacktrackingEnd}
		    \State 
		    \label{codeBTreturn}
		    \textbf{return} $L_t$, $M_t$, $\Gamma_t$, $\nu_t$, $x\tn$, $x\ag[1]$, $\tau_t$, $\rho_t$, $\alpha_t$
		\EndProcedure
	\end{algorithmic}
	\end{algorithm}
}	
A few remarks are in place for the above backtracking procedure. Firstly, steps \ref{codeBacktracking} through \ref{codeBacktrackingEnd} are the backtracking steps, which terminates only when the conditions in steps \ref{codeLt} and \ref{codeMt} are both satisfied. Clearly, in each call to the backtracking procedure, steps \ref{codeLt} and \ref{codeMt} will only be performed finitely many times, and the returned values $L_t$ and $M_t$ satisfies $L_{min}\le L_t\le 2L_G$ and $M_t\le 2\|K\|$, respectively. Secondly, while $M_t\ge M_{t-1}$ and $\nu_{t}\ge \nu_{t-1}$, the value of $L_{t}$ in step \ref{codeBTreturn} is not necessarily greater than $L_{t-1}$. Finally, the multiplier for increasing or decreasing ${L}_t$ and ${M}_t$ is 2, which can be replaced by any number that is greater than 1.

The scheme of AADMM with backtracking is presented in Algorithm \ref{algAPDL}.
\begin{algorithm}[h]
	\caption{\label{algAPDL} AADMM with backtracking}
	\begin{algorithmic}[]
\makeatletter
\setcounter{ALG@line}{-1}
\makeatother
			\State Choose $x_1\in X$ and $w_1\in \cW$ such that $Bw_1=Kx_1+b$, $L_0\ge L_{min}>0$ and $M_0, \nu_0, \rho>0$. Set $x\ag[][1] \leftarrow x_1$, $w\ag[][1] \leftarrow w_1$, $y\ag[][1]\leftarrow y_1=0$,  $\Gamma_0\leftarrow L_0$, $t\leftarrow 1$.
			\For{$t=1,\cdots,N-1$}
				\State ($L_t$, $M_t$, $\Gamma_t$, $\nu_t$, $x\tn$, $x\ag[1]$, $\tau_t$, $\rho_t$, $\alpha_t$)$\leftarrow$BACKTRACKING($L_{t-1}$, $M_{t-1}$, $\Gamma_{t-1}$, $\nu_{t-1}$, $x_t$, $x\ag$, $L_{min}$)
				\State Calculate iterates \eqref{eqnAADMMwtn} -- \eqref{eqnAADMMyag}.
			\EndFor
	\end{algorithmic}
\end{algorithm}

We start by considering UCO problems with bounded feasible sets $X$ and $Y$. Theorem \ref{thmAADMML} below summarizes the convergence properties of Algorithm \ref{algAPDL} for solving bounded UCO problems.

\begin{thm}
	\label{thmAPDL}
	If we set $\nu_0 = -\infty$ and apply Algorithm \ref{algAPDL} to the UCO problem \eqref{eqnUCO} under assumption \eqref{BD}, then
		\begin{align}
			\label{eqnAPDLrho}
			\begin{aligned}
			&\ G(x\ag[1]) + F(Kx\ag[1]) - f^*
			\\
			\leq &\ \frac{4L_GD_X^2}{t^2} + \frac{4L_G}{L_{min}(t-1)}\left[6{\chi{\rho}\max\{4M_0^2, \|K\|^2\}D_X^2}  + (1-\chi){\rho}D_{X, K}^2+ \frac{D_Y^2}{{\rho}}\right].
			\end{aligned}
		\end{align}	
	In particular, if $\rho=\frac{D_Y}{\sqrt{6}\chi\max\{2M_0, \|K\|\}D_X + (1+\chi)D_{X,K}}$, then
		\begin{align}
			\label{eqnAPDL}
			\begin{aligned}
			&\ G(x\ag[1]) + F(Kx\ag[1]) - f^*
			\\
			\leq &\ \frac{4L_GD_X^2}{t^2} + \frac{4L_G}{L_{min}(t-1)}\left[\sqrt{6}{\chi\max\{2M_0, \|K\|\}D_XD_Y}  + (1-\chi)D_{X, K}D_Y\right].
			\end{aligned}
		\end{align}		
	\begin{proof}
		As discussed after Procedure \ref{funBT}, we have
		\begin{align}
			\label{eqnLMUB}
			L_{min}\le L_t\le 2L_G \text{ and } 0\le M_t\le 2\|K\|.
		\end{align}
		We can now estimate the bounds of $\alpha_t$ and $\Gamma_t$. 
		 By \eqref{eqnGamma} we have $1/\Gamma_t = 1/\Gamma_{t-1} + \alpha_t/\Gamma_t$, hence
		\begin{align*}
			\sqrt{\frac{1}{\Gamma_t}} - \sqrt{\frac{1}{\Gamma_{t-1}}} = \frac{1/\Gamma_t - 1/\Gamma_{t-1}}{\sqrt{1/\Gamma_t} + \sqrt{1/\Gamma_{t-1}}} = \frac{\alpha_t/\Gamma_t}{\sqrt{1/\Gamma_t} + \sqrt{1/\Gamma_{t-1}}}.
		\end{align*}
		Observing from equations \eqref{eqnGamma},  \eqref{eqnAlphaL} and \eqref{eqnLMUB} that
		\begin{align}
			\label{eqnAlphaGammaBound}
			1/(2L_G)\le \alpha_t^2/\Gamma_t\le 1/L_{min}, 
		\end{align}
		we have
		\begin{align*}
			\sqrt{\frac{1}{\Gamma_t}} - \sqrt{\frac{1}{\Gamma_{t-1}}} \ge&\  \frac{\alpha_t/\Gamma_t}{2\sqrt{1/\Gamma_t}} = \frac{\alpha_t}{2\sqrt{\Gamma_t}} \ge \frac{1}{2\sqrt{2L_G}}, \text{ and }
			\\
			\sqrt{\frac{1}{\Gamma_t}} - \sqrt{\frac{1}{\Gamma_{t-1}}} \le&\   \frac{\sqrt{1/(L_{min}\Gamma_t)}}{\sqrt{1/\Gamma_t} + \sqrt{1/\Gamma_{t-1}}} \le \frac{1}{\sqrt{L_{min}}}.
		\end{align*}
		Therefore, by induction we conclude that
		\begin{align}
			\label{eqnGammaBound}
			\frac{t}{2\sqrt{2L_G}}\le \sqrt{\frac{1}{\Gamma_t}} \le \frac{t}{\sqrt{L_{min}}} + \frac{1}{\sqrt{L_0}}\le \frac{t+1}{\sqrt{L_{min}}}, \text{ or }\frac{L_{min}}{(t+1)^2}\le \Gamma_t \le \frac{8L_G}{t^2}.
		\end{align}
		Now let us examine the RHS of \eqref{eqnQBound} in Lemma \ref{lemQBound}. Without loss of generality, we assume that $2M_0\le \|K\|$. Indeed, if $2M_0>\|K\|$, then $M_t\equiv 2M_0$ for all $t\ge 1$.
		Since $\nu_{[t]}$ and $M_{[t]}$ are monotonically increasing, by \eqref{eqnStepL} and \eqref{eqnBBD} in Lemma \ref{lemB}, we have
		\begin{align*}
			\cBx[^*] \le&\ \frac{1+\chi\nu_t\rho M_t^2}{2}(D_X^2 - \|x\tn - x^*\|^2) \le \frac{1+\chi\nu_t\rho M_t^2}{2}D_X^2
			\\
			\le&\ \frac{1}{2}D_X^2 + 2{\chi\nu_t\rho }\|K\|^2D_X^2
			,
			\\
			\cBy = &\ \frac{\nu_t}{2\rho}(\|y_1 - y\|^2 - \|y\tn - y\|^2)\le \frac{\nu_t}{2\rho}D_Y^2, \ \forall y\in \cY
			,\\
			\cB_t(w^*, w_{[t+1]}, \theta_{[t]}) = &\ \frac{\nu_t\rho}{2}(\|w_1 - w^*\|^2 - \|w\tn - w^*\|^2) \leq \frac{\nu_t\rho}{2}D_{X, K}^2.
		\end{align*}
		On the other hand, by \eqref{eqnBUB} in Lemma \ref{lemB} we have
		\begin{align*}
			-\chi\cB_t(Kx^*, Kx_{[t+1]}, \theta_{[t]})\le &\ -\frac{\chi\nu_1\rho}{2}(\|Kx_1 - Kx^*\|^2 - \|Kx\tn - Kx^*\|^2) \le
			\frac{\chi\nu_1\rho}{2}\|K\|^2\|x\tn - x^*\|^2 
			\\
			\le&\ \frac{\chi\nu_t\rho}{2}\|K\|^2D_X^2
			.
		\end{align*}
		Applying the above calculations on $\cB(\cdot,\cdot,\cdot)$ to Lemma \ref{lemQBound}, we have
		\begin{align*}
			&\ \frac{1}{\Gamma_t}Q(w^*, x^*, y; z\ag[1]) \le \frac{D_X^2}{2} + \frac{\nu_t\rho }{2}D_{X, K}^2 + \frac{5\chi\nu_t\rho}{2}\|K\|^2D_{X}^2 + \frac{\nu_t}{2\rho}D_Y^2
			\\
			\le &\ \frac{D_X^2}{2} + {3\chi\nu_t\rho}\|K\|^2D_{X}^2  + \frac{(1-\chi)\nu_t\rho }{2}D_{X, K}^2 + \frac{\nu_t}{2\rho}D_Y^2
			.
		\end{align*}
		Observe that by \eqref{eqnAlphaGammaBound} and \eqref{eqnGammaBound}, $\alpha_t/\Gamma_t\le (t+1)/L_{min}$, and that
		\begin{align}
		\label{eqnvtbound}
		\nu_t\le \max_{i=1,\ldots t}\alpha_i/\Gamma_i \le (t+1)/L_{min}. 
		\end{align}
		Using the previous two inequalities and \eqref{eqnGammaBound}, we have
		\begin{align*}
			&\ \bar g_Y(z\ag[1]) 
			\le \frac{4L_GD_X^2}{t^2}+ \frac{24\chi\rho L_G(t+1)}{t^2L_{min}}\|K\|^2D_X^2 + \frac{4(1-\chi)\rho L_G(t+1)}{t^2L_{min}}D_{X,K}^2 + \frac{4L_G(t+1)}{t^2L_{min}\rho}D_Y^2
			\\
			\le &\ \frac{4L_GD_X^2}{t^2}+ \frac{24\chi\rho L_G}{L_{min}(t-1)}\|K\|^2D_X^2   + \frac{4L_G}{L_{min}\rho(t-1)}D_Y^2.
		\end{align*}
		The above inequality, in view of Proposition \ref{proGap0}, then implies  \eqref{eqnAPDLrho} and \eqref{eqnAPDL}.
	\end{proof}
\end{thm}

\vgap

For AECCO problems when both $X$ and $Y$ are bounded, we can also apply Algorithm \ref{algAPDL} with $\chi=0$, as long as the maximum number of iterations $N$ is given. Theorem \ref{thmAADMML} below describes the convergence properties of AL-ADMM with backtracking for solving general AECCO problems.
\begin{thm}
	\label{thmAADMML}
		If we choose $\chi=0$, $\rho=1$, and $\nu_0= N/L_{min}$ in Algorithm \ref{algAPDL}, then
	\begin{align}
		\label{eqnAPDLUprim}&
\begin{aligned}
		G(x\ag[][N]) + F(w\ag[][N]) - f^*
		\leq &\ \frac{4L_GD_{x^*}^2}{(N-1)^2} + \frac{4L_GD_{x^*,K}^2}{L_{min}(N-1)},\text{ and }
\end{aligned}
		\\
		\label{eqnAPDLUfeas}&
\begin{aligned}
		\|\BKb[{\ag[][N]}]\|
		\leq &\  \frac{16\sqrt{L_G}D_{x^*}}{\sqrt{L_{min}}(N-1)^{3/2}} +  \frac{16\sqrt{2}\sqrt{L_G}D_{x^*,K}}{\sqrt{L_{min}}(N-1)} + \frac{32L_GD_{y^*}}{L_{min}(N-1)}.
\end{aligned}
	\end{align}
	\begin{proof}
		In view of step \ref{codeBacktracking} in Procedure \ref{funBT}, equation \eqref{eqnvtbound} and the choice of $\nu_0$, we can see that $\nu_t\equiv N/L_{min}$. By \eqref{eqnGamma}, \eqref{eqnAlphaL}, \eqref{eqnB} and \eqref{eqnStepL}, we have
		\begin{align*}
			\cBx[^*] =&\ \frac{1}{2}(D_{x^*}^2 - \|x\tn - x^*\|^2) \le \frac{1}{2}D_{x^*}^2,
			\\
			\cBy = &\ \frac{N}{2L_{min}}(\|y_1 - y\|^2 - \|y\tn - y\|^2), \ \forall y\in \cY
			,\\
			\cB_t(Bw^*, Bw_{[t+1]}, \theta_{[t]}) = &\ \frac{N}{2L_{min}}(\|Bw_1 - Bw^*\|^2 - \|B\tn - Bw^*\|^2) \leq \frac{N}{2L_{min}}D_{x^*, K}^2
			.
		\end{align*}
		Using the fact that $\tau_t\ge \rho_t$ and $\chi=0$, and applying the above calculations to Lemma \ref{lemQBound}, we have
		\begin{align*}
			\frac{1}{\Gamma_t}Q(w^*, x^*, y; z\ag[1]) 
			\leq&\ \frac{1}{2}D_{x^*}^2 + \frac{N}{2L_{min}}D_{x^*,K}^2 + \frac{N}{2L_{min}}(\|y_1 - y\|^2 - \|y\tn - y\|^2), \ \forall y\in \cY.
		\end{align*}	
		Similarly to the proof of Theorem \ref{thmAADMMII}, we have
		\begin{align*}
			&\|y\tn - y^*\|^2 \le\ \frac{L_{min}}{N}D_{x^*}^2 + D_{x^*,K}^2 + D_{y^*}^2,\ 
			\|y_1 - y\tn\|^2 \le\ \frac{2L_{min}}{N}D_{x^*}^2 + 2D_{x^*,K}^2 + 4D_{y^*}^2,\text{ and }
			\\
			&\frac{1}{\Gamma_t}Q(w^*, x^*, y; z\ag[1]) + \frac{N}{L_{min}}\langle y_1 - y\tn, y\rangle 
			\leq\ \frac{1}{2}D_{x^*}^2 + \frac{N}{2L_{min}}D_{x^*,K}^2,\ \forall y\in\cY.	
		\end{align*}
		Setting $v\tn = {\Gamma_tN}(y_1 - y\tn)/L_{min}$, $t=N-1$ and applying \eqref{eqnGammaBound}, we have 
	\begin{align}
		&\begin{aligned}
				Q(w^*, x^*, y; z\ag[][N]) + \langle v_{N}, y\rangle 
				\leq &\ \frac{4L_GD_{x^*}^2}{(N-1)^2} + \frac{4L_GD_{x^*,K}^2}{L_{min}(N-1)},
		\end{aligned}
				\\
		&\begin{aligned}
				\|v_N\|
				\leq &\  \frac{8\sqrt{2}\sqrt{L_GN}D_{x^*}}{\sqrt{L_{min}}(N-1)^{2}} +  \frac{8\sqrt{2}N\sqrt{L_G}D_{x^*,K}}{\sqrt{L_{min}}(N-1)^2} + \frac{16NL_GD_{y^*}}{L_{min}(N-1)^2}
				\\
				\le &\ \frac{16\sqrt{L_G}D_{x^*}}{\sqrt{L_{min}}(N-1)^{3/2}} +  \frac{16\sqrt{2}\sqrt{L_G}D_{x^*,K}}{\sqrt{L_{min}}(N-1)} + \frac{32L_GD_{y^*}}{L_{min}(N-1)}.
		\end{aligned}
	\end{align}
	These previous two relations together with Proposition \ref{proGap} then imply \eqref{eqnAPDLUprim} and \eqref{eqnAPDLUfeas}.
	\end{proof}
\end{thm}

\section{Numerical examples}

In this section, we will present some preliminary numerical results of the proposed methods. The numerical experiments are carried out on overlapped LASSO, compressive sensing, and an application on partially parallel image reconstruction. All algorithms are implemented in MATLAB 2013b on a Dell Precision T1700 computer with 3.4 GHz Intel i7 processor.

\subsection{Group LASSO with overlap}
The goal of this section is to examine the effectiveness of the proposed methods for solving UCO problems with unbounded $X$. In this experiment, our problem of interest is the group LASSO model given by \cite{jacob2009group}
\begin{align}
	\label{eqnOverlasso}
	\min_{x\in\R^n}\sum_{i=1}^{m}(\langle a_i,x\rangle - f_i)^2 +  \lambda\sum_{g\in\cG}\|x_g\|,
\end{align}
where $\{(a_i,f_i) \}_{i=1}^m\subseteq \R^n\times\R$ is a group of datasets, $x$ is the sparse feature to be extracted, and the structure of $x$ is represented by group $\cG$. In particular, $\cG\subseteq 2^{\{1,\ldots,n\}}$, and for any $g\subseteq \{1,\ldots,n\}$, $x_g$ is a vector that is constructed by components of $x$ whose indices are in $g$, i.e.,  $x_g:=(x_i)_{i\in g}$. The first term in \eqref{eqnOverlasso} describes the fidelity of data observation, and the second term is the regularization term to enforce certain group sparsity. In particular, we assume that $x$ is sparse in the group-wise fashion, i.e., for any $g\in\cG$, $x_g$ is sparse. Problem \eqref{eqnOverlasso} can be formulated as a UCO problem \eqref{eqnUCO} by defining the linear operator $K$ as $Kx = \lambda (x_{g_1}^T, x_{g_2}^T,\ldots, x_{g_l}^T)^T$, where $g_i\in\cG$ and $\cG = \{g_i\}_{i=1}^{l}$. Specially, if each $g_i$ consists $k$ elements, then \eqref{eqnOverlasso} becomes
\begin{align}
	\min_{x\in \R^n} \frac{1}{2}\|Ax - f\|^2 + \lambda\|Kx\|_{k,1},
\end{align}
where $A=(a_1,\ldots,a_m)^T$, $f=(f_1,\ldots,f_m)^T$, and $\|\cdot\|_{k,1}$ is defined by $\|u\|_{k,1}:=\sum_{i=1}^{n}\|(u^{(ki-k+1)},\ldots,u^{(ki)})^T\|$ for all $u\in\R^{kn}$, where $\|\cdot\|$ is the Euclidean norm in $\R^k$. Note that $F(\cdot):=\|\cdot\|_{k,1}$ is simple, so the solution of problem \eqref{eqnF} can be obtained directly by examining the optimality condition, which is also known as soft-thresholding. 

\begin{figure}[h]
	\begin{subfigure}[b]{.48\linewidth}
	\includegraphics[width=\linewidth]{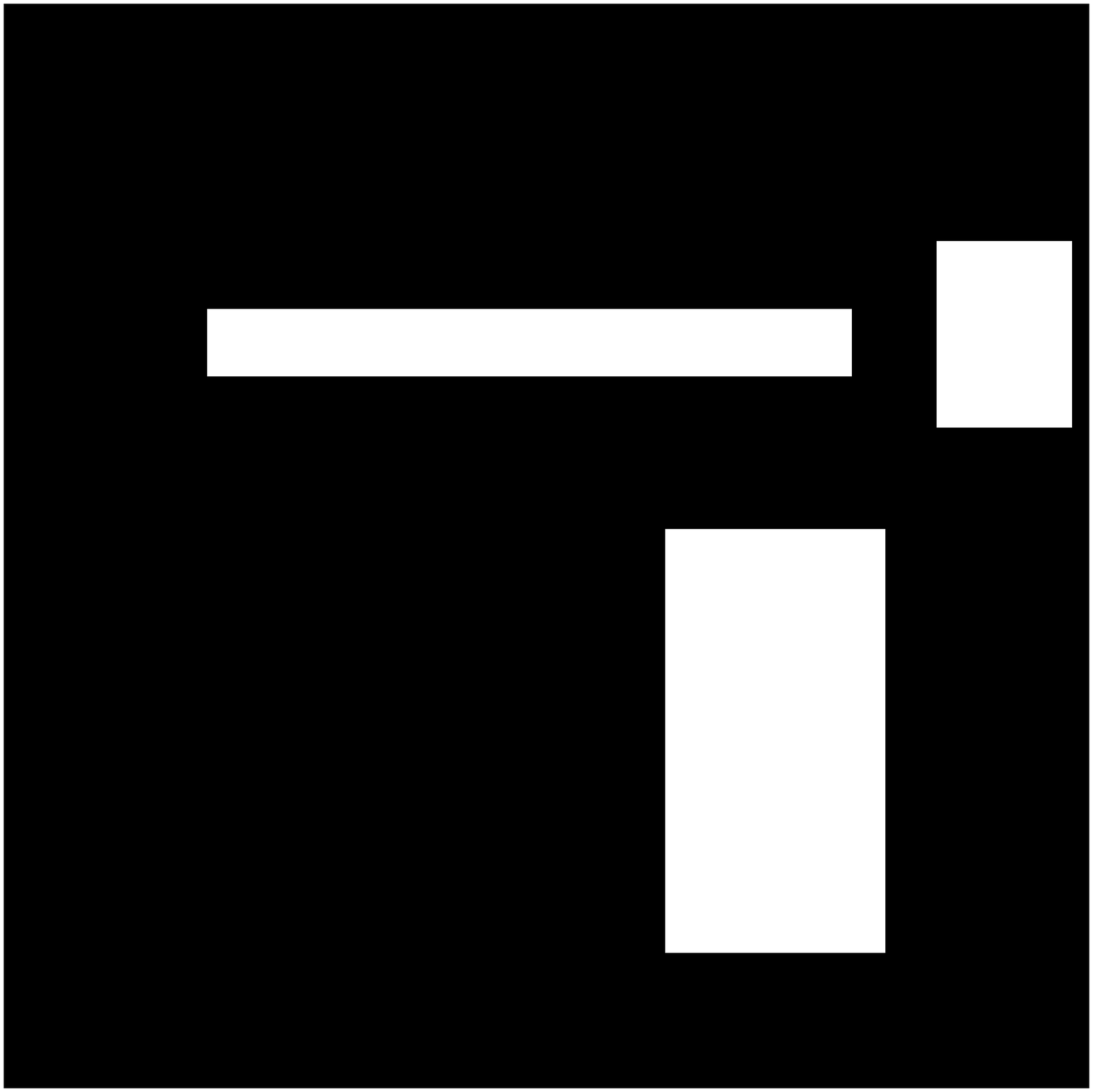}
	\end{subfigure}
	\begin{subfigure}[b]{.48\linewidth}
	\includegraphics[width=\linewidth]{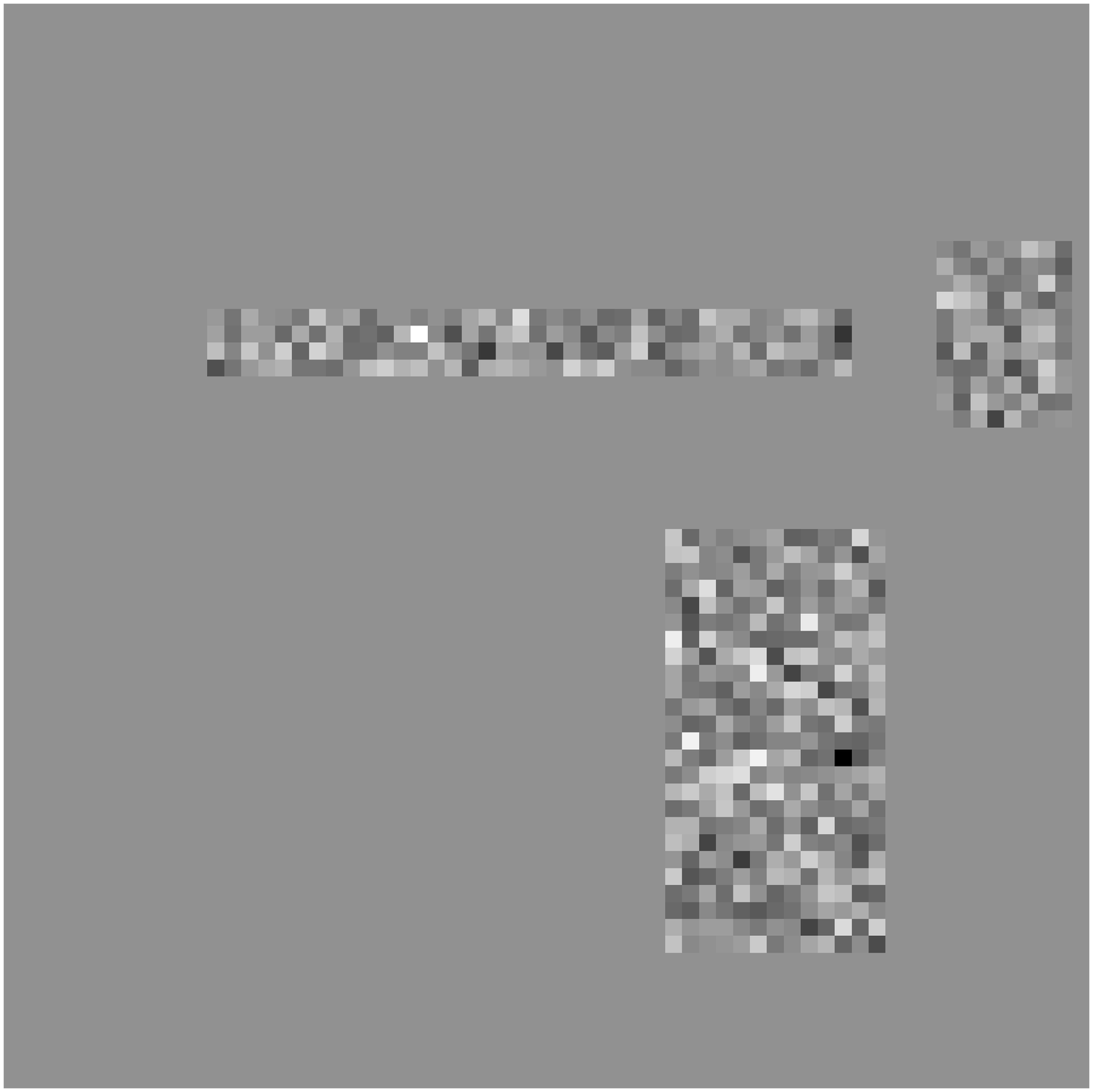}
	\end{subfigure}
	\caption{\label{figOverlassoTrue}\footnotesize True feature $x_{true}$ in the experiment of group LASSO with overlap. Left: the support of $x_{true}$. Right: the intensities of $x_{true}$.}
	\begin{subfigure}{.49\linewidth}
		\includegraphics[width=\linewidth]{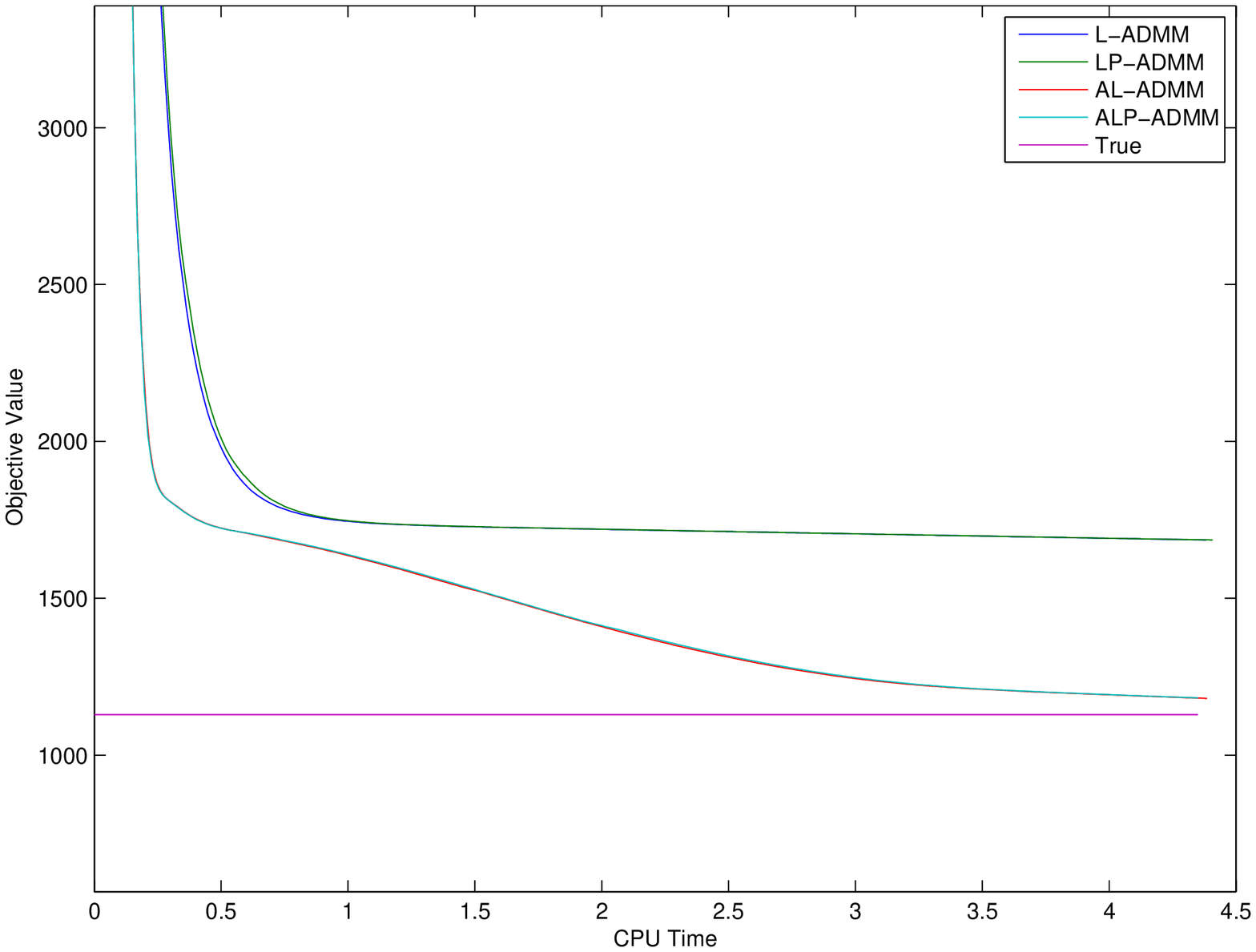}
	\end{subfigure}
	\begin{subfigure}{.49\linewidth}
		\includegraphics[width=\linewidth]{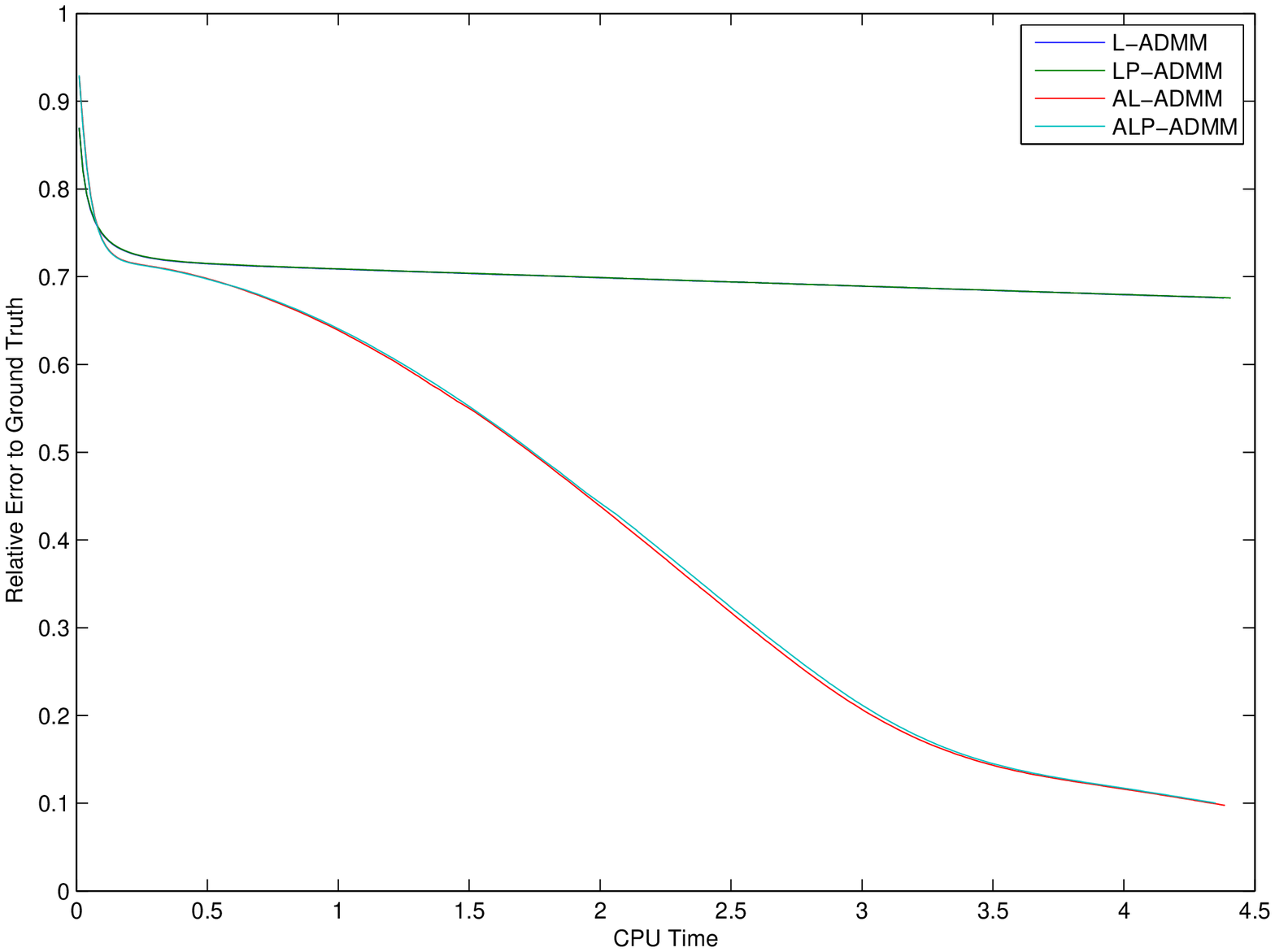}
	\end{subfigure}
	\caption{\label{figOverlasso} \footnotesize Comparisons of AL-ADMM, ALP-ADMM, L-ADMM and LP-ADMM in group LASSO with overlap. Left: the objective function values $f(x_t^{ag})$ 
		from AL-ADMM and ALP-ADMM, and $f(x_t)$ from L-ADMM and LP-ADMM vs. CPU time. The straight line at the bottom is $f(x_{true})$. Right:  the
	relative errors $r(x_t^{ag})$ from AL-ADMM and ALP-ADMM and $r(x_t)$ from L-ADMM and LP-ADMM vs. CPU time.}
\end{figure}

In this experiment, we generate the datasets $\{(a_i,f_i)\}_{i=1}^{m}$ by $f_i = \langle a_i, x_{true}\rangle + \varepsilon$, where $a_i\sim N(0,I_n)$, $\varepsilon\sim N(0,0.01)$, and the true feature $x_{true}$ is the n-vector form of a $64 \times 64$ two-dimensional signal whose support and intensities are shown in Figure \ref{figOverlassoTrue}. Within its support, the intensities of $x_{true}$ are generated independently from standard normal distribution. We set $n=4096$, $m=2048$ and choose $\cG$ to be all the $2\times 2$ blocks in the $64\times 64$ domain (so that $k=4$), and apply L-ADMM, LP-ADMM, AL-ADMM and ALP-ADMM to solve \eqref{eqnOverlasso} in which $\lambda=1$. The parameters for AL-ADMM and ALP-ADMM are chosen as in Theorem \ref{thmAADMMII}, and $N$ is set to $300$. To have a fair comparison, we use the same Lipschitz constants $L_G=\lambda_{max}(A^TA)\approx 1.6\times 10^4$, $\|K\|=2$ and $\rho=0.5$ for all algorithms without performing a backtracking. 
Both the primal objective function value
$f(\tilde x)$ and the feature extraction relative error $r(\tilde{x})$ at approximate solution $\tilde{x}\in\R$ versus CPU time are reported in Figure \ref{figOverlassoTrue}, where
\begin{align}
	\label{eqnRelErr}
	r(\tilde{x}):=\frac{\|\tilde{x} - x_{true}\|}{\|x_{true}\|}.
\end{align}

From Figure \ref{figOverlassoTrue} we can see that the performance of AL-ADMM and ALP-ADMM are almost the same, and both of them outperforms L-ADMM and LP-ADMM. This is consistent with our theoretical observations that AL-ADMM and ALP-ADMM have better rate of convergence \eqref{eqnAADMMIIprimSimple} than ADMM \eqref{eqnADMMIIprimSimple}.

\subsection{Compressive sensing}
In this subsection, we present the experimental results on the comparison of ADMM and AADMM for solving the following image reconstruction problem:
\begin{align}
	\label{eqnTVRecon}
	\min_{x\in X}\frac{1}{2}\|Ax-f\|^2 + \lambda\|D x\|_{2,1},
\end{align}
where $x$ is the $n$-vector form of a two-dimensional image to be reconstructed,
$\|Dx\|_{2,1}$ is the discrete form of the TV semi-norm,
$A$ is a given acquisition matrix (depending on the physics of the data acquisition), $f$ represents the observed
data, and $X:=\{x\in\R^n:l_*\le x^{(i)} \le u_*, \forall i=1,\ldots, n \}$. Problem \eqref{eqnTVRecon} is a special case of UCO \eqref{eqnUCO} with $\cW=\R^{2n}$, $G(x)=\|Ax-b\|^2/2$,  $F(w) = \|w\|_{2,1}$ and $K = \lambda D$. We assume that the finite difference operator $D$ satisfies the periodic boundary condition, so that the problem in \eqref{eqnAADMMxtn} with $\chi=0$ can be solved easily by utilizing the Fourier transform (see \cite{wang2008new}). 

\begin{figure}[h]
	\begin{subfigure}[b]{.49\linewidth}
		\includegraphics[width=\linewidth]{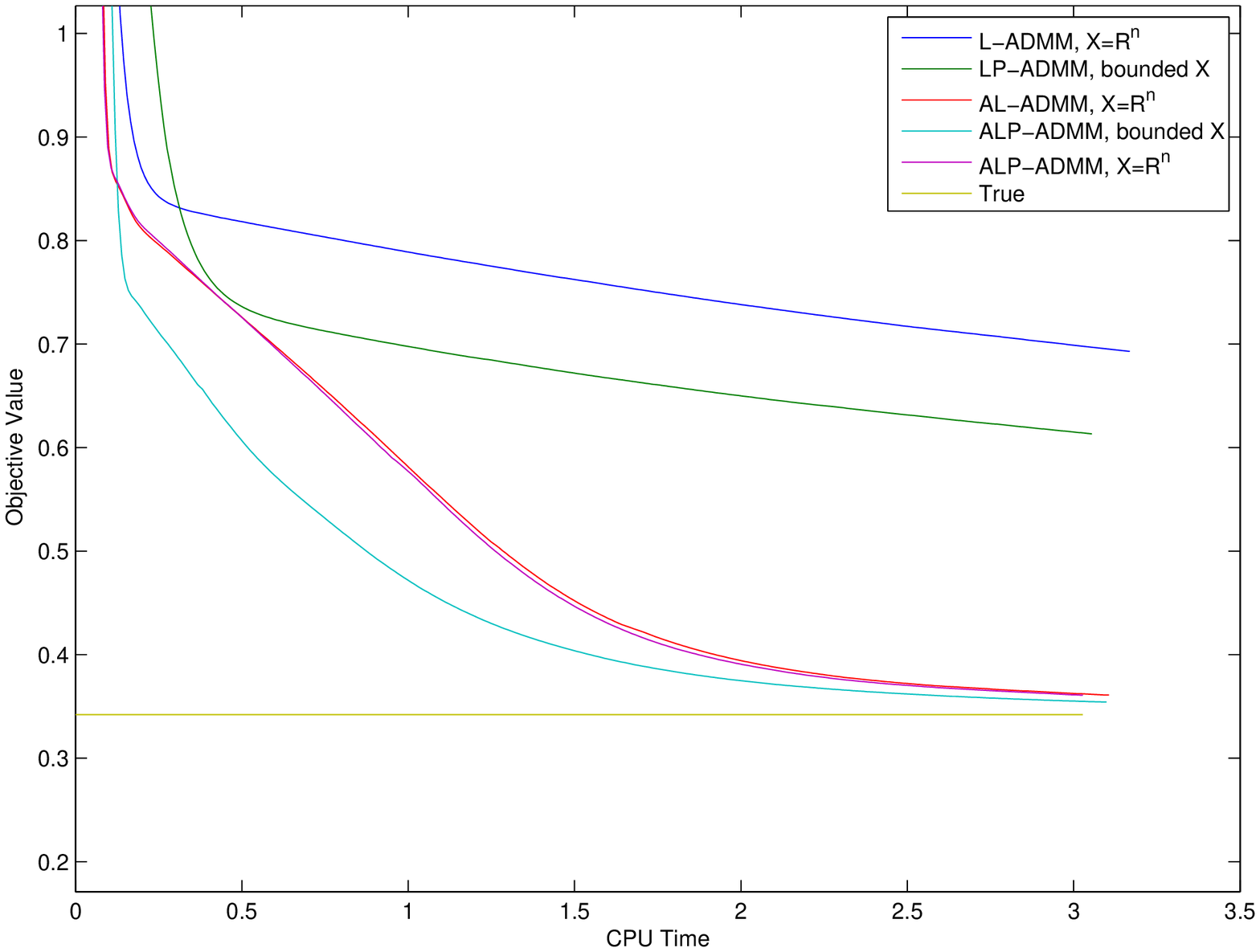}
	\end{subfigure}
	\begin{subfigure}[b]{.49\linewidth}
		\includegraphics[width=\linewidth]{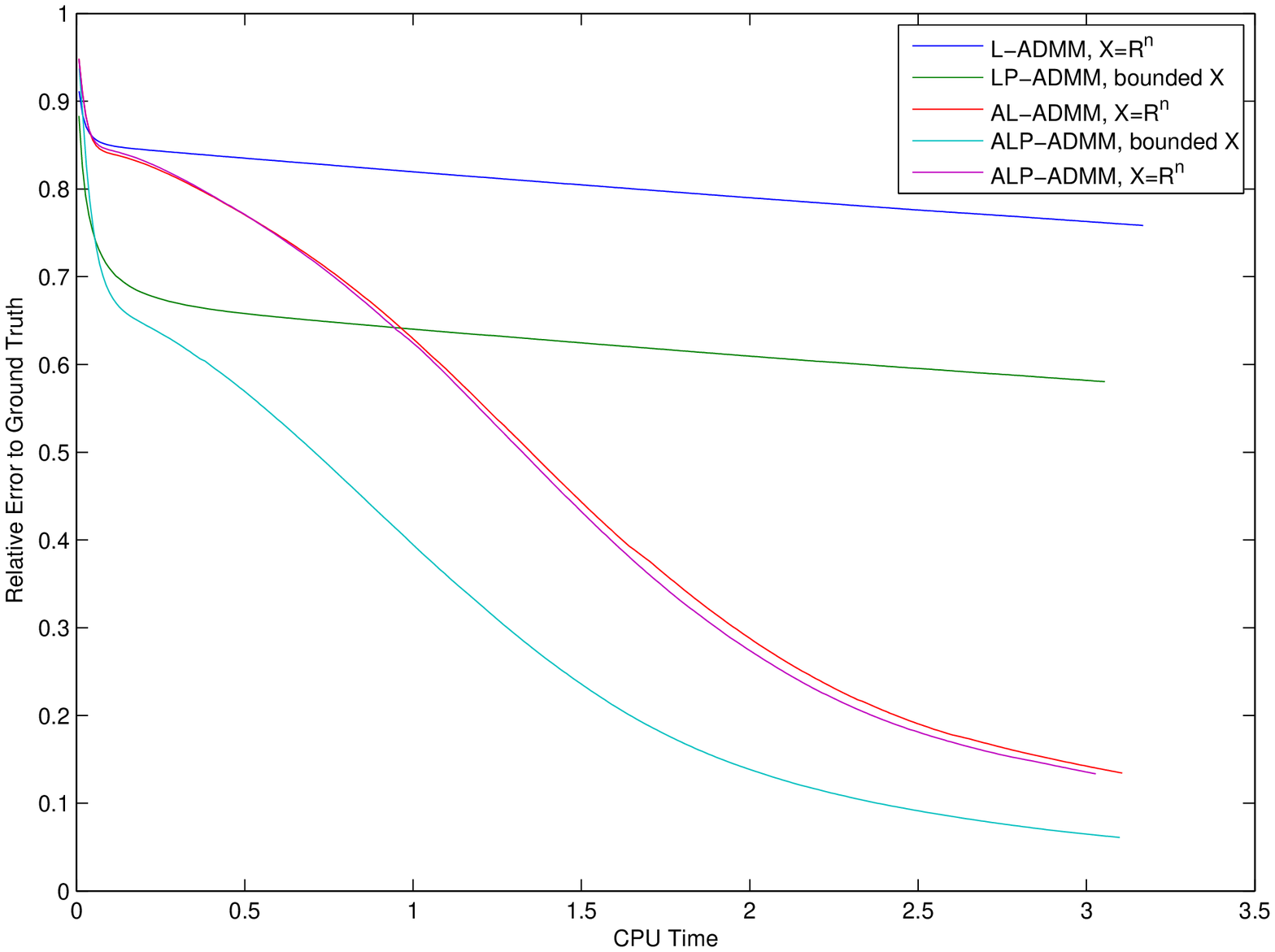}
	\end{subfigure}
	\\
	\begin{subfigure}[b]{.49\linewidth}
		\includegraphics[width=\linewidth]{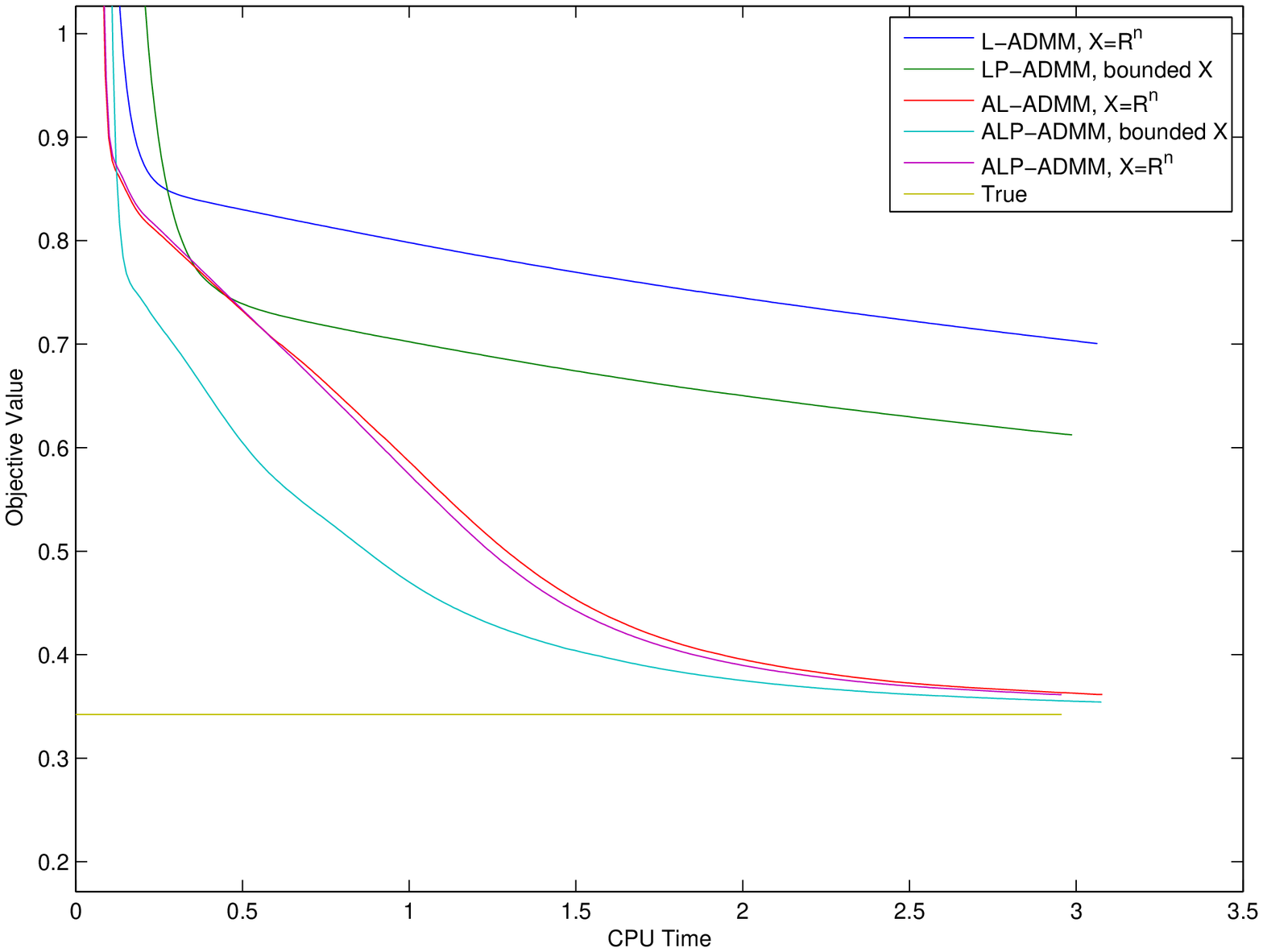}
	\end{subfigure}
	\begin{subfigure}[b]{.49\linewidth}
		\includegraphics[width=\linewidth]{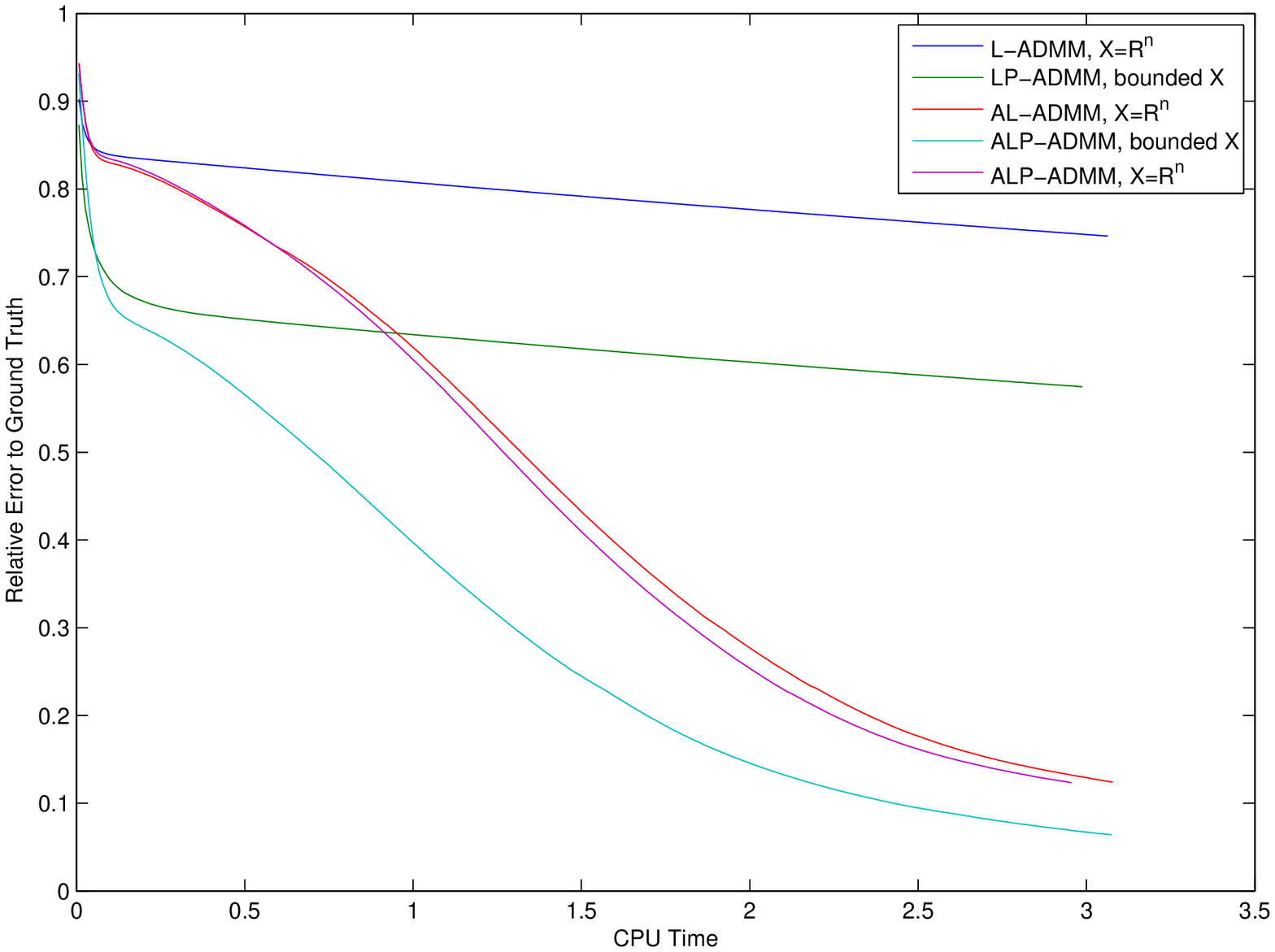}
	\end{subfigure}
	\caption{\label{figBernoulli} \footnotesize
	Comparisons of AL-ADMM, ALP-ADMM, L-ADMM and LP-ADMM in image reconstruction. The top and bottom rows, 
	respectively, show
	the performance of these algorithms on the ``Gaussian" and ``Bernoulli" instances. Left: the objective function values $f(x_t^{ag})$ 
	from AL-ADMM and ALP-ADMM, and $f(x_t)$ from L-ADMM and LP-ADMM vs. CPU time. The straight line at the bottom is $f(x_{true})$. Right:  the
relative errors $r(x_t^{ag})$ from AL-ADMM and ALP-ADMM, and $r(x_t)$ in L-ADMM and LP-ADMM vs. CPU time.}
\end{figure}

In our experiment, we consider two instances where the acquisition matrix $A \in R^{m \times n}$ is generated independently from a normal distribution $N(0, 1/\sqrt{m})$ and a Bernoulli distribution that takes equal probability for the values $1/\sqrt{m}$ and $-1/\sqrt{m}$ respectively.
Both types of acquisition matrices are widely used in compressive sensing (see,
e.g., \cite{baraniuk2008simple}).  
For a given $A$, the measurements $b$ are generated by $b = Ax_{true}+\varepsilon$, where $x_{true}$ is a $64$
by $64$ Shepp-Logan phantom \cite{shepp1974fourier} with intensities in $[0,1]$ (so $n=4096$), and
$\varepsilon\equiv N(0, 0.001 I_n)$. We choose $m=1229$ so that the compression ratio is about $30\%$, and set $\lambda=10^{-3}$ in \eqref{eqnTV}. Considering the range of intensities of $x_{true}$, we apply ALP-ADMM with parameters in Theorem \ref{thmAPD} and LP-ADMM to solve \eqref{eqnTVRecon} with  bounded feasible set $X:=\{x\in\R^n:0\le x^{(i)} \le 1, \forall i=1,\ldots, n \}$. It should be pointed that since $Y:=\dom F^*=\{y\in\R^{2n}:\|y\|_{2,\infty}:=\max_{i=1,\ldots,n}\|(y^{2i-1},y^{2i})^T\|_2\leq 1\}$, we have $D_X=D_Y=n$, which suggests that $\rho=1/\|K\|$ may be a good choice for $\rho$. We also apply L-ADMM and AL-ADMM to solve \eqref{eqnTVRecon}, with $\chi=1$ and $X=\R^n$. In this case we use the parameters in Theorem \ref{thmAADMMII} with $N=300$ for AL-ADMM. To have a fair comparison, we use the same constants $L_G=\lambda_{max}(A^TA)$ and $\|K\|=\lambda\sqrt{8}$ (see \cite{chambolle2004algorithm}) and $\rho=1/\|K\|$ for all
algorithms without performing backtracking. We report both the primal objective function value
and the reconstruction relative error \eqref{eqnRelErr} versus CPU time in Figure \ref{figBernoulli}.

It is evident from Figure \ref{figOverlasso} that AL-ADMM and ALP-ADMM outperforms L-ADMM and LP-ADMM in solving \eqref{eqnOverlasso}. This is consistent with our theoretical results in Corollaries \ref{thmADMMI}, \ref{thmAPD}, \ref{thmADMMII} and \ref{thmAADMMII}. Moreover, it is interesting to observe that ALP-ADMM with box constrained $X$ outperforms AL-ADMM with $X=\R^n$. This suggests that the knowledge of the ground truth is helpful in solving image reconstruction problems.

\subsection{Partially parallel imaging}
\begin{figure}[h]
	\begin{subfigure}[b]{.49\linewidth}
		\includegraphics[width=\linewidth]{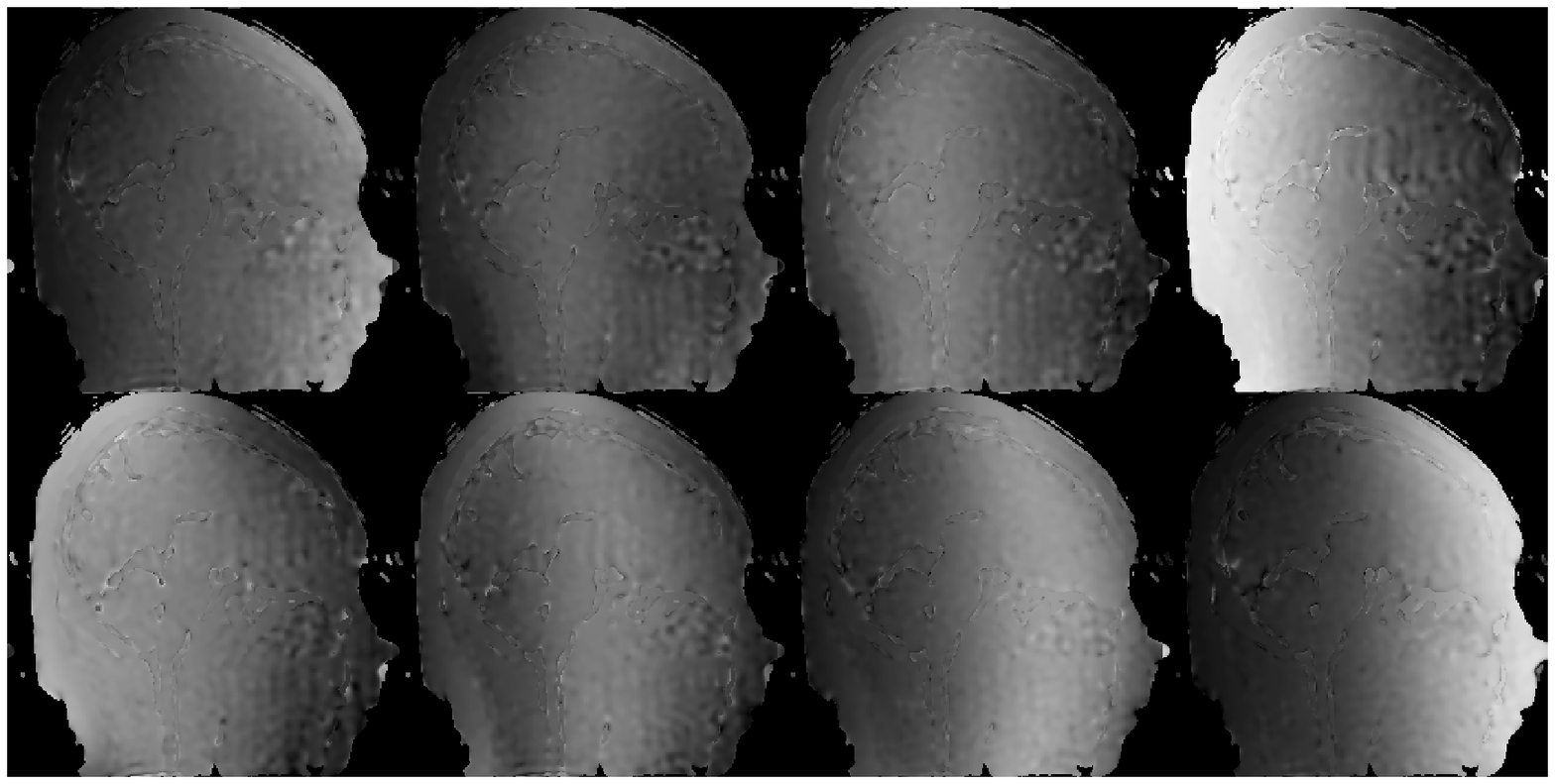}
		\caption{}
	\end{subfigure}
	\begin{subfigure}[b]{.24\linewidth}
		\includegraphics[width=\linewidth]{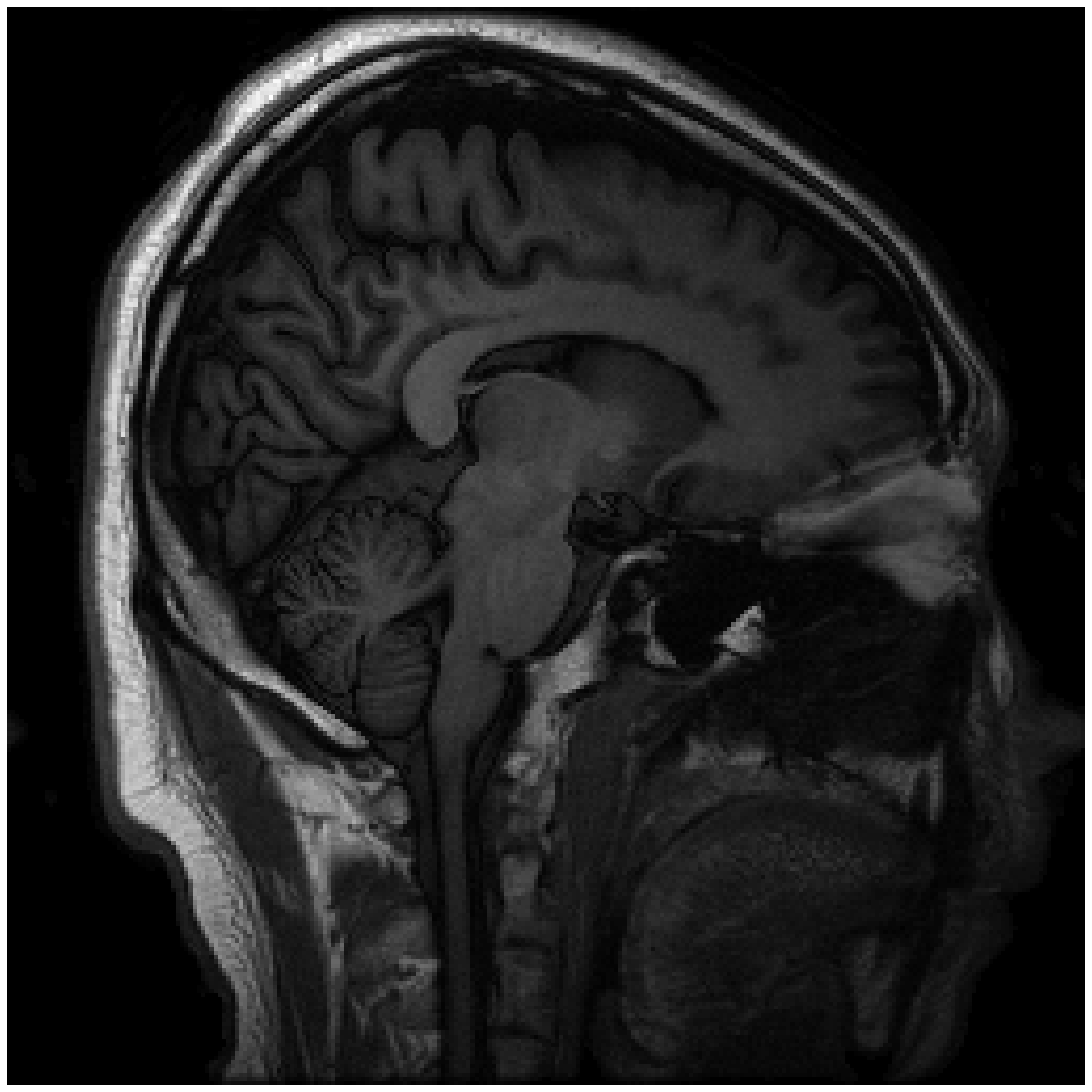}
		\caption{}
	\end{subfigure}
	\begin{subfigure}[b]{.12\linewidth}
	\begin{subfigure}[b]{\linewidth}
			\includegraphics[width=\linewidth]{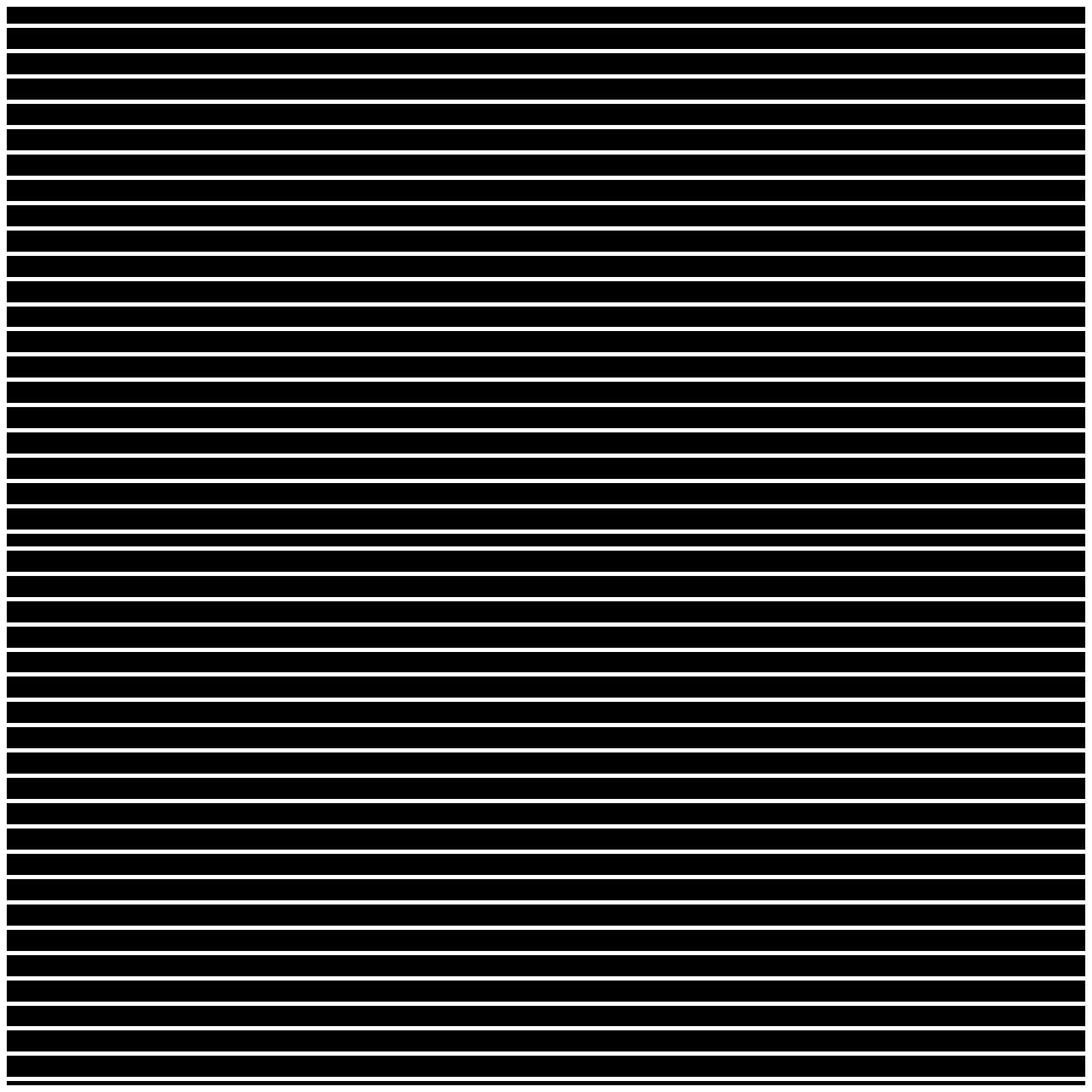}\\
			\includegraphics[width=\linewidth]{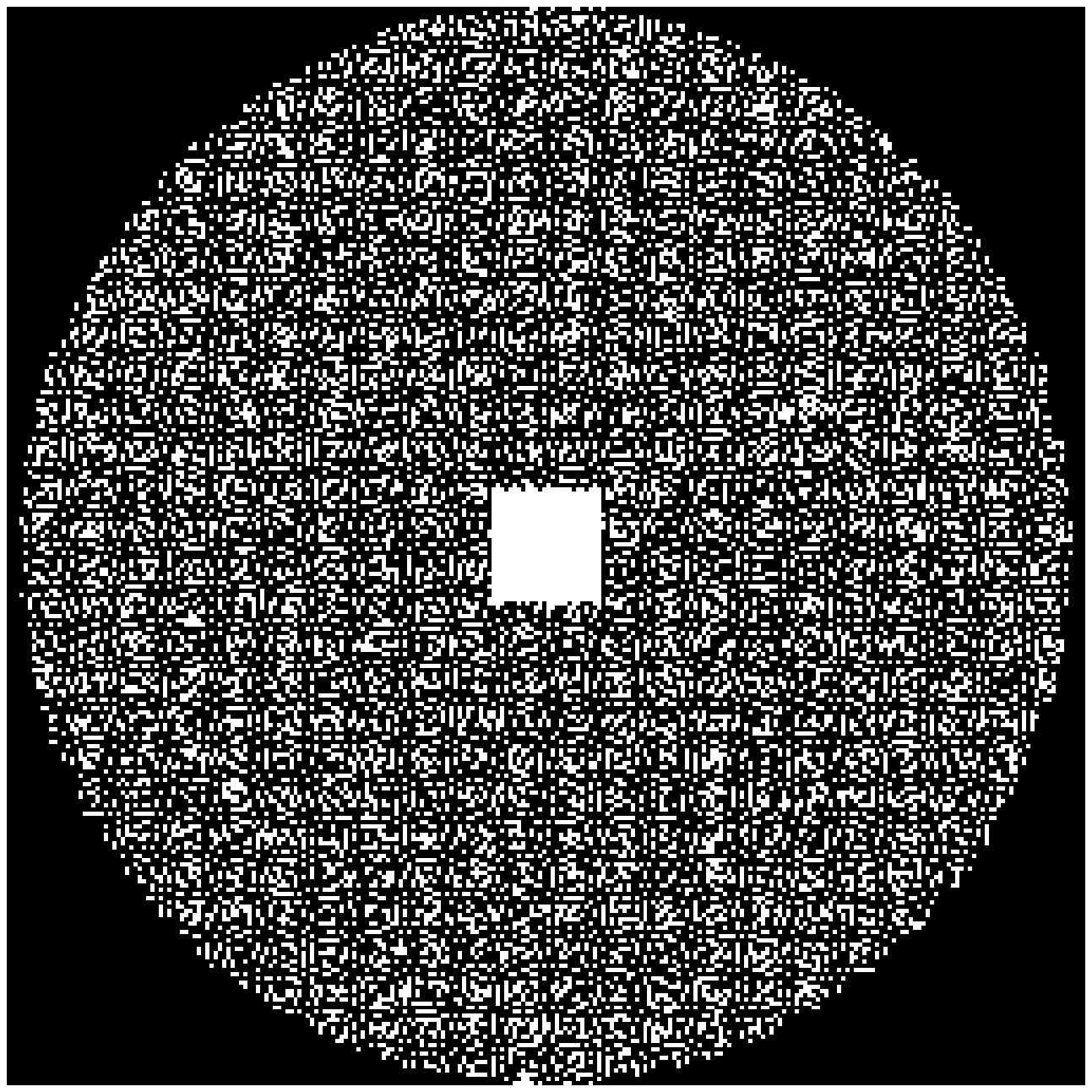}
	\end{subfigure}
	\caption{}
	\end{subfigure}
	\\
	\begin{subfigure}[b]{.49\linewidth}
		\includegraphics[width=\linewidth]{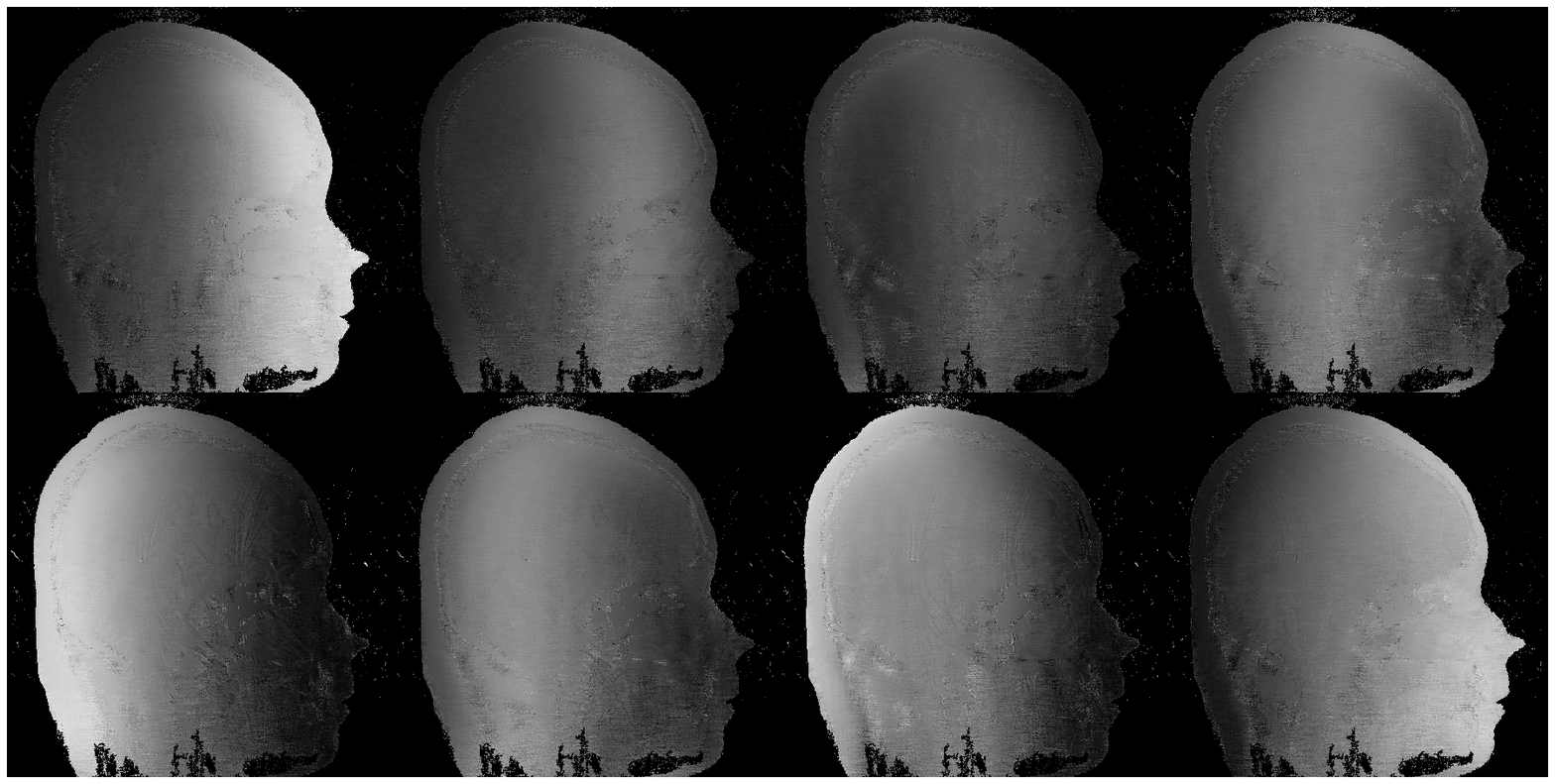}
		\caption{}
	\end{subfigure}
	\begin{subfigure}[b]{.24\linewidth}
		\includegraphics[width=\linewidth]{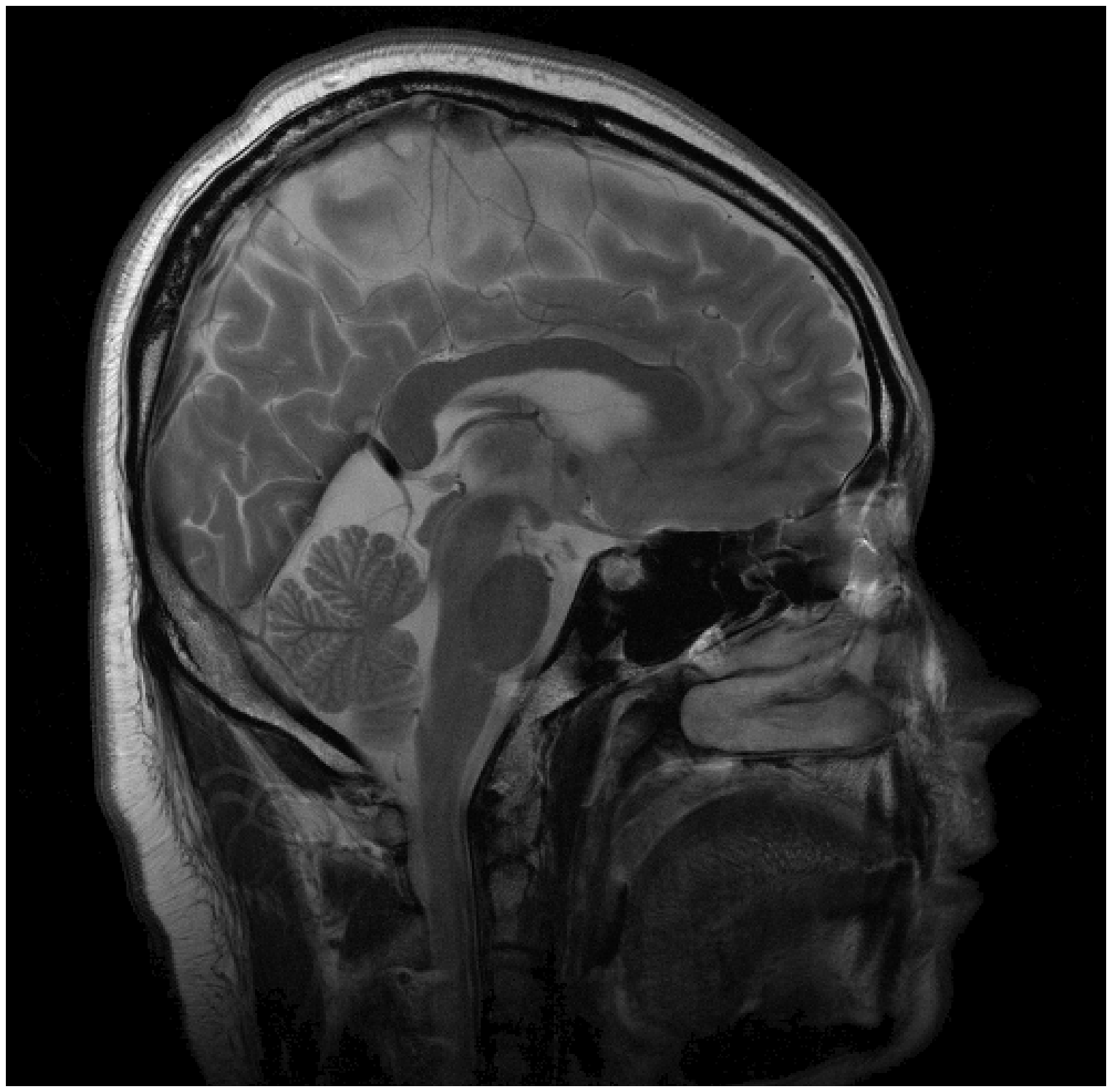}
		\caption{}
	\end{subfigure}
	\begin{subfigure}[b]{.12\linewidth}
	\begin{subfigure}[b]{\linewidth}
			\includegraphics[width=\linewidth]{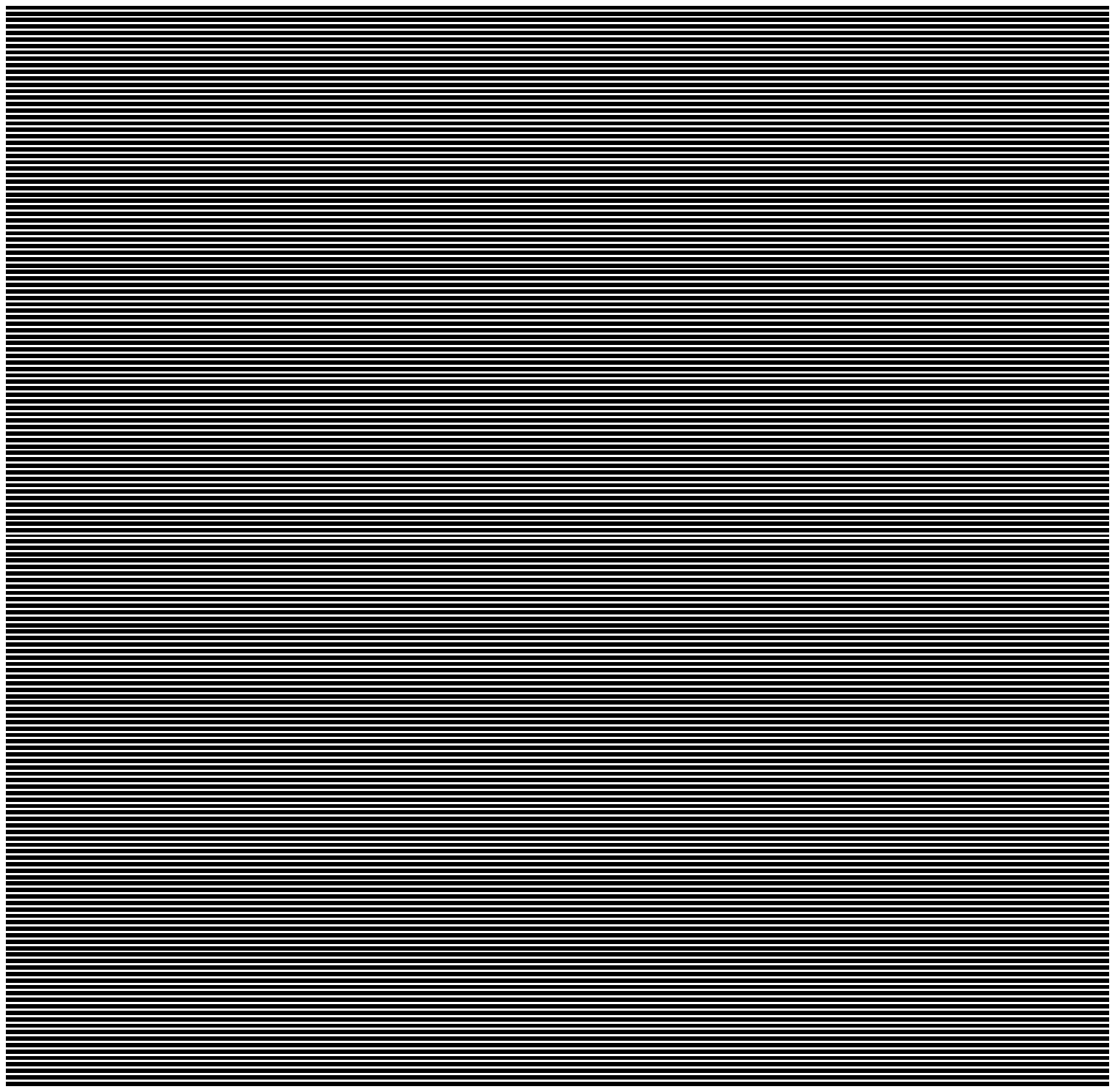}\\
			\includegraphics[width=\linewidth]{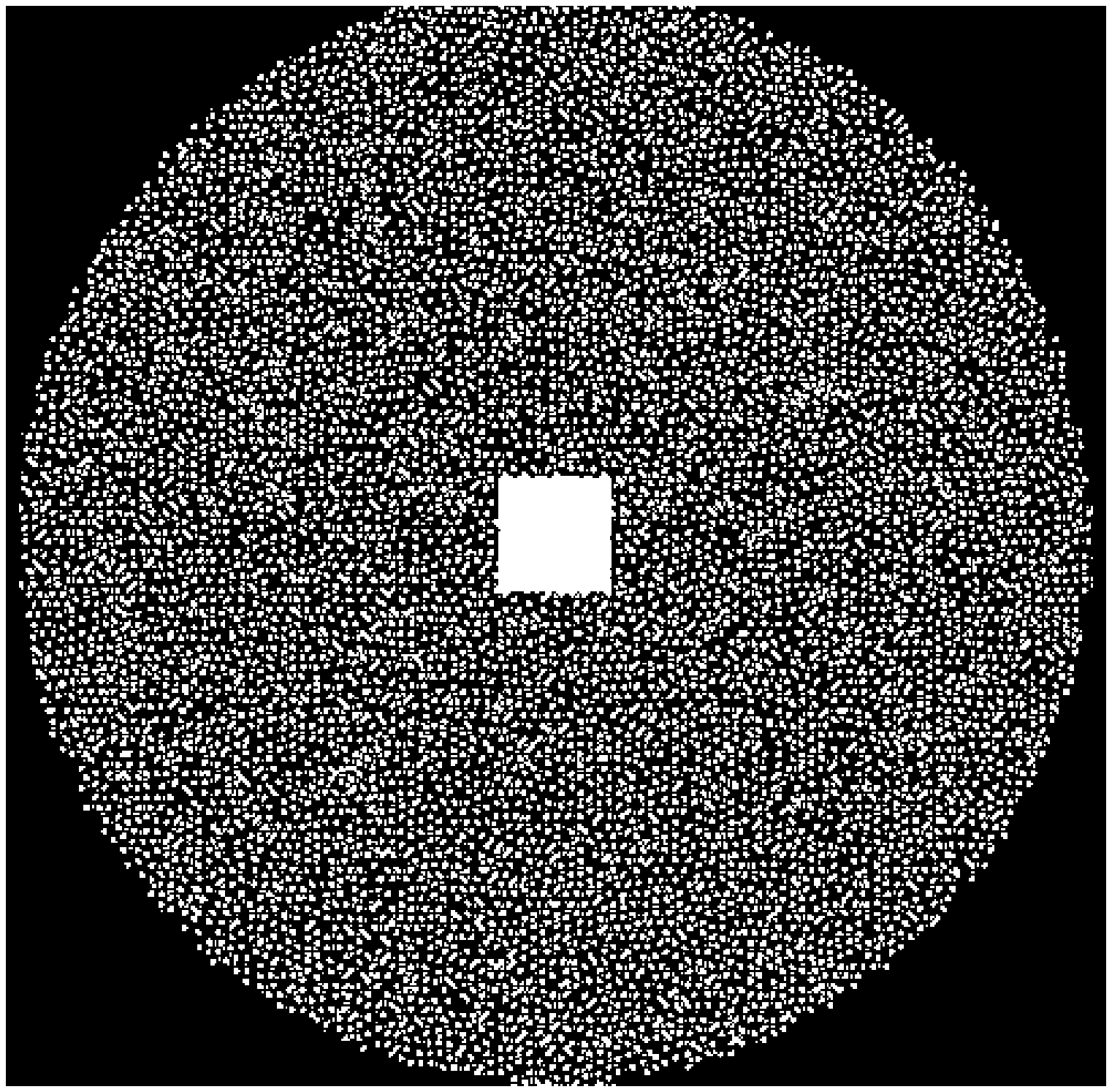}
	\end{subfigure}
		\caption{}
	\end{subfigure}
	\caption{\footnotesize\label{figPPITrue} Sensitivity map $\{\diag S_j\}_{j=1}^{8}$ (left), ground truth $x_{true}$ (middle) and mask $\diag M$ (right) in partially parallel image reconstruction. (a): The sensitivity maps in instances 1a and 1b. (b): The ground truth in instances 1a and 1b. (c): The k-space sampling trajectory in instances 1a (top) and 1b (bottom). (d): The sensitivity maps in instances 2a and 2b. (e): The ground truth in instances 2a and 2b. (f): The k-space sampling trajectory in instances 2a (top) and 2b (bottom).}
\end{figure}

	In this section, we compare the performance of AADMM with backtracking and Bregman operator splitting with variable stepsize (BOSVS) \cite{chen2013bregman}, which is a linearized ADMM method with backtracking, in reconstruction of magnetic resonance images from partially parallel imaging (PPI). In magnetic resonance PPI, a set of multi-channel k-space data is acquired simultaneously from  radiofrequency (RF) coil arrays. The imaging is accelerated by sampling a reduced number of k-space samples.
	The image reconstruction problem can be modeled as
	\begin{align}
		\label{eqnPPI}
		\min_{x\in X}\frac{1}{2}\sum_{j=1}^{n_{ch}}\|MFS_jx - f_j\|^2 + \lambda\|D x\|_{2,1},
	\end{align}	
	where $x$ is the vector form of a two-dimensional image to be reconstructed. In \eqref{eqnPPI}, $n_{ch}$ is the number of MR sensors, $F\in \C^{n\times n}$ is a 2D discrete Fourier transform matrix,  $S_j\in\C^{n\times n}$ is the sensitivity encoding map of the $j$-th sensor, and $M\in\R^{n\times n}$ describes the scanning pattern of MR sensors, and $X\subseteq \C^n$. In particular, $S_j$'s and $M$ are both diagonal matrices, and their diagonal vectors $\diag S_j\in\R^n$ and $\diag M\in\R^n$ are n-vector form of images that have the same dimension as the reconstructed image. In practice, $\diag S_j$ describes the sensitivity of the $j$-th sensor at each pixel, and $\diag M$ is a mask that takes value ones at the scanned pixels and zeros elsewhere. Figure \ref{figPPITrue} shows the two-dimensional image representations of $\{\diag S_j\}_{j=1}^{n_{ch}}$, $x_{true}$ and $\diag M$. The PPI reconstruction problems are described in more details in \cite{chen2012fast}. It should be noted that \eqref{eqnPPI} is a special case of \eqref{eqnTVRecon}, and that the percentage of nonzero elements in $\diag M$ describes the compression ratio of PPI scan. In view of the fact that $\|F\| = \sqrt{n}$, the Lipschitz constant $L_G$ of \eqref{eqnPPI} can be estimated by
	\begin{align}
		L_G = \|\sum_{j=1}^{n_{ch}}S_jF^TM^2FS_j\| \le n\|\sum_{j=1}^{n_{ch}}S_j\|^2 = n\|\sum_{j=1}^{n_{ch}}\diag S_j\|_{\infty}^2.
	\end{align}

	In this experiment, $n_{ch}=8$, and the measurements $\{f_j\}_{j=1}^{n_{ch}}$ are generated by
		\begin{align*}
			f_j = M(FS_jx_{true} + \varepsilon_j^{re}/\sqrt{2} + \varepsilon_j^{im}/\sqrt{-2}),\ j=1,\ldots,n_{ch}
		\end{align*}
	where the noises $\epsilon_j^{re}, \epsilon_j^{im}$ are independently generated from distribution ${\cal N}(0, 10^{-4}\sqrt{n}I_n)$. We generate four instances of experiments where the ground truth $x_{true}$ are the human brain image (see Figure \ref{figPPITrue}). The information of the instances is listed in Table \ref{tabPPIinfo}. In particular, instances 1a and 1b have Cartesian and pseudo-random k-space sampling trajectories respectively but share the same sensitivity map and ground truth, and so are instances 2a and 2b.
		
	\begin{table}
		\caption{\label{tabPPIinfo} Data acquisition information in partially parallel image reconstruction.}
		\centering
		\begin{tabular}{ccccc}
			\hline
			Instance & Dimension of $x_{true}$ & Sampling trajectory & Acquisition rate & $L_G$ 
			\\\hline
			1a & $n=256\times 256$ & Cartesian mask & 18\% & $3.34\times 10^{5}$
			\\\hline
			1b & $n=256\times 256$ & Pseudo random mask  & 24\% & $3.34\times 10^{5}$
			\\\hline
			2a & $n=512\times 512$ & Cartesian mask & 18\% & $1.60\times 10^{6}$
			\\\hline
			2b & $n=512\times 512$ & Pseudo random mask & 24\% & $1.60\times 10^{6}$
			\\\hline
		\end{tabular}
	\end{table}

	We first consider $X=\C^n$, and use AL-ADMM with backtracking to solve \eqref{eqnPPI}. We use the parameters in Theorem \ref{thmAADMML} with $N=400$ in all PPI experiments. We also apply the BOSVS method in \cite{chen2013bregman}\footnote{The BOSVS code is available at http://people.math.gatech.edu/\textasciitilde xye33/software/BOSVS.zip} to solve \eqref{eqnPPI} with $X=\C^n$, which is a backtracking linesearch technique for L-ADMM with Barzilai-Borwein stepsize \cite{barzilai1988two}. Furthermore, noticing that $x_{true}$ is in bounded feasible set $X:=\{x\in\C^n:|x^{(i)}|\le 1,\forall i=1,\ldots,n\}$, we also apply ALP-ADMM with backtracking to solve \eqref{eqnPPI} with aforementioned bounded feasible set $X$. We set the parameters to $\lambda = 10^{-10}n$ in \eqref{eqnPPI}, and choose $L_0 = \|F\|^2 = n$, $L_{min} = L_G/10$, $M_0 = \|K\|/10 = \lambda\sqrt{8}/10$ for Algorithm \ref{algAPDL} where $L_G$ is listed in Table \ref{tabPPIinfo}.
	
	The performance of AL-ADMM, ALP-ADMM and BOSVS is shown in Figures \ref{figPPI} and \ref{figPPIRecon}, in terms of both the primal objective function value and relative error \eqref{eqnRelErr}. It is evident that AL-ADMM and
	ALP-ADMM outperform BOSVS in terms of the decrement of both primal objective value and relative error to ground truth, especially in the case of using Cartesian sampling trajectory. Since the Cartesian sampling trajectory in our experiments collects less low-frequency data (the center part in the k-space)  and has no randomness in sampling (see Figure \ref{figPPITrue}), it makes harder to get a good reconstruction comparing with that of the pseudo-random sampling trajectory. Our experimental results indicates that in this case the AADMM is much more efficient than BOSVS in reconstruction. It is evident that AL-ADMM and ALP-ADMM outperform BOSVS in terms of the decrement of both primal objective value and relative error to ground truth. This observation is consistent with our theoretical result in Theorems \ref{thmAPDL} and \ref{thmAADMML}.

\section{Conclusion}
We present in this paper the AADMM framework by incorporating a multi-step acceleration scheme into linearized ADMM. AADMM has better rates of convergence than linearized ADMM on solving a class of convex composite optimization with linear constraints, in terms of the Lipschitz constant of the smooth component. Moreover, AADMM can handle both bounded and unbounded feasible sets, as long as a saddle point exists. For the unbounded case, the estimation for the rate of convergence depends on the distance from initial point to the set of saddle points. We also propose a backtracking scheme to improve the practical performance of AADMM. Our preliminary numerical results show that AADMM is promising for solving large-scale convex composition optimization with linear constraints.

\section*{Acknowledgment}
The authors would like to thank Invivo Philips, Gainesville, FL for providing the PPI brain scan datasets.

\bibliographystyle{./plain_first_initial}
\bibliography{./yuyuan}{}

\begin{figure}[h]
	\begin{subfigure}[b]{.49\linewidth}
		\includegraphics[width=\linewidth]{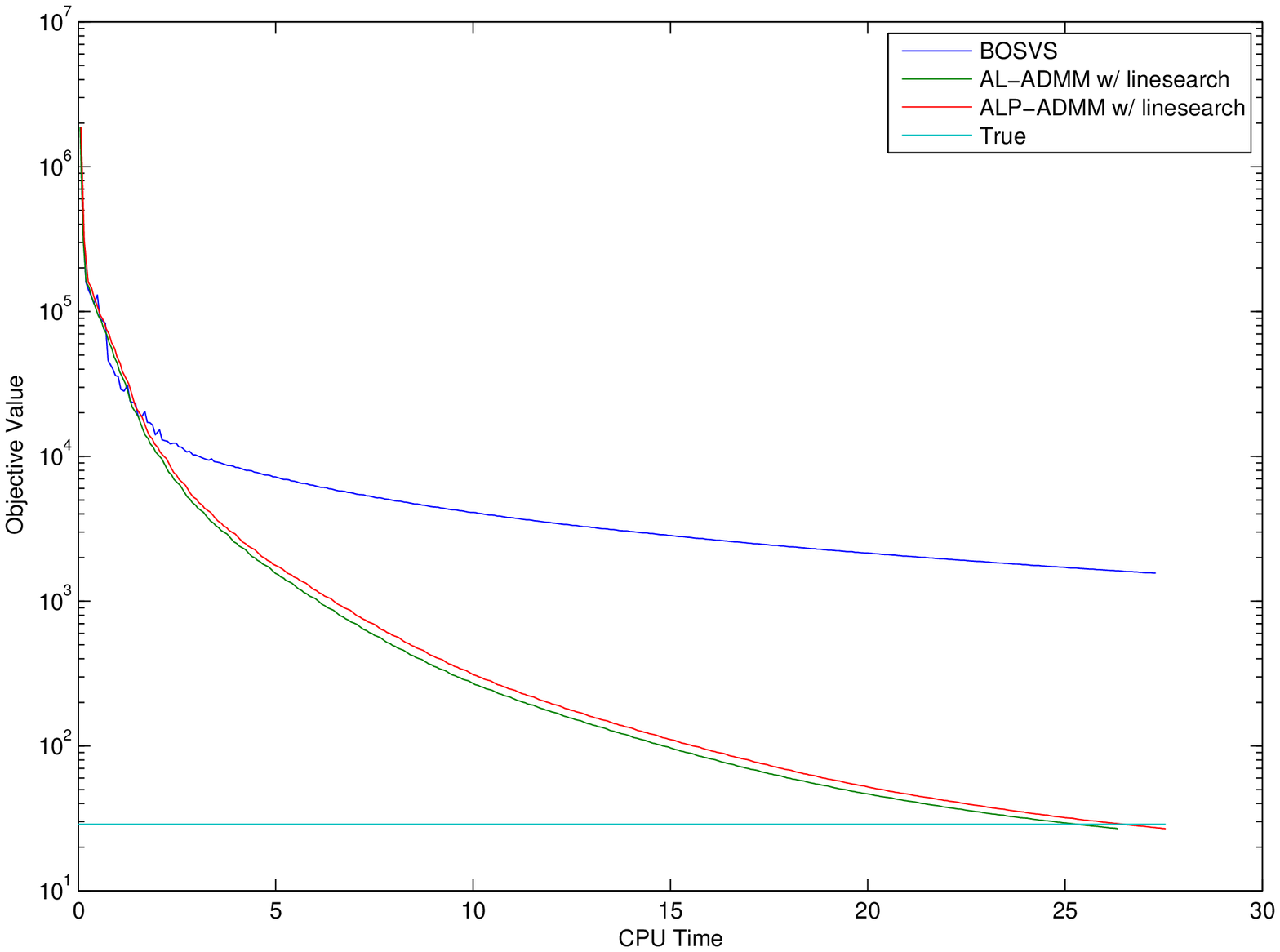}
	\end{subfigure}
	\begin{subfigure}[b]{.49\linewidth}
		\includegraphics[width=\linewidth]{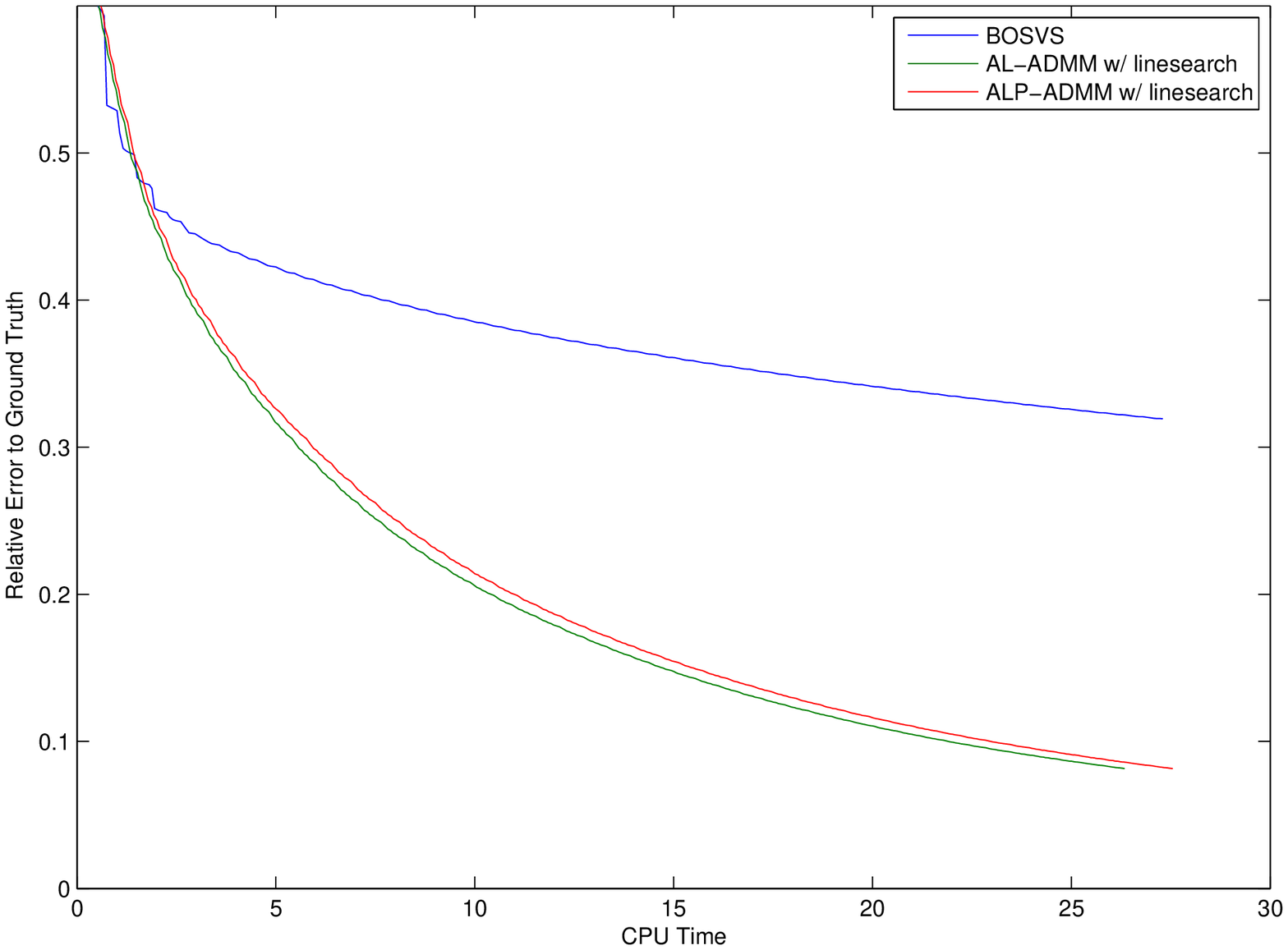}
	\end{subfigure}
	\begin{subfigure}[b]{.49\linewidth}
		\includegraphics[width=\linewidth]{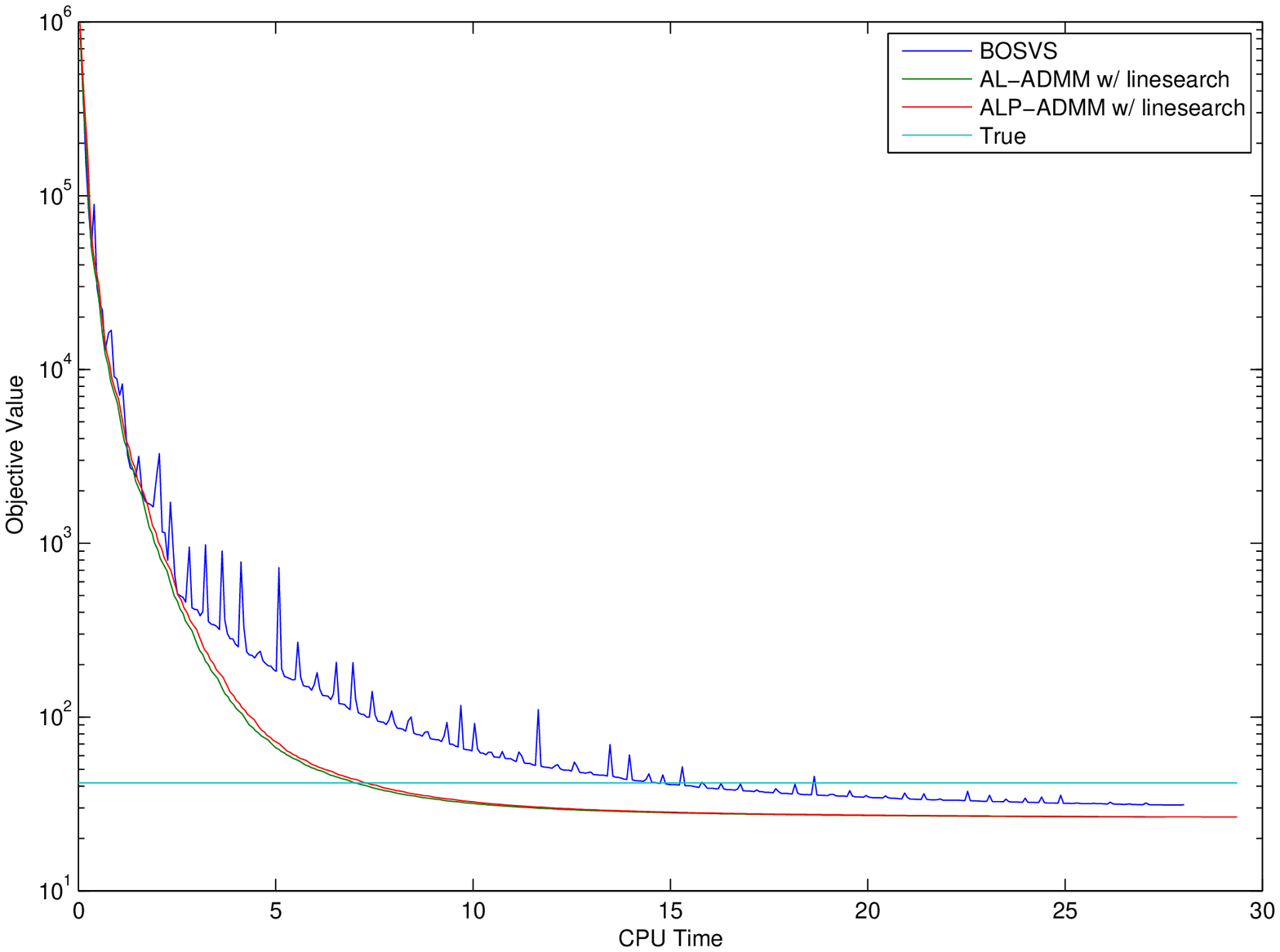}
	\end{subfigure}
	\begin{subfigure}[b]{.49\linewidth}
		\includegraphics[width=\linewidth]{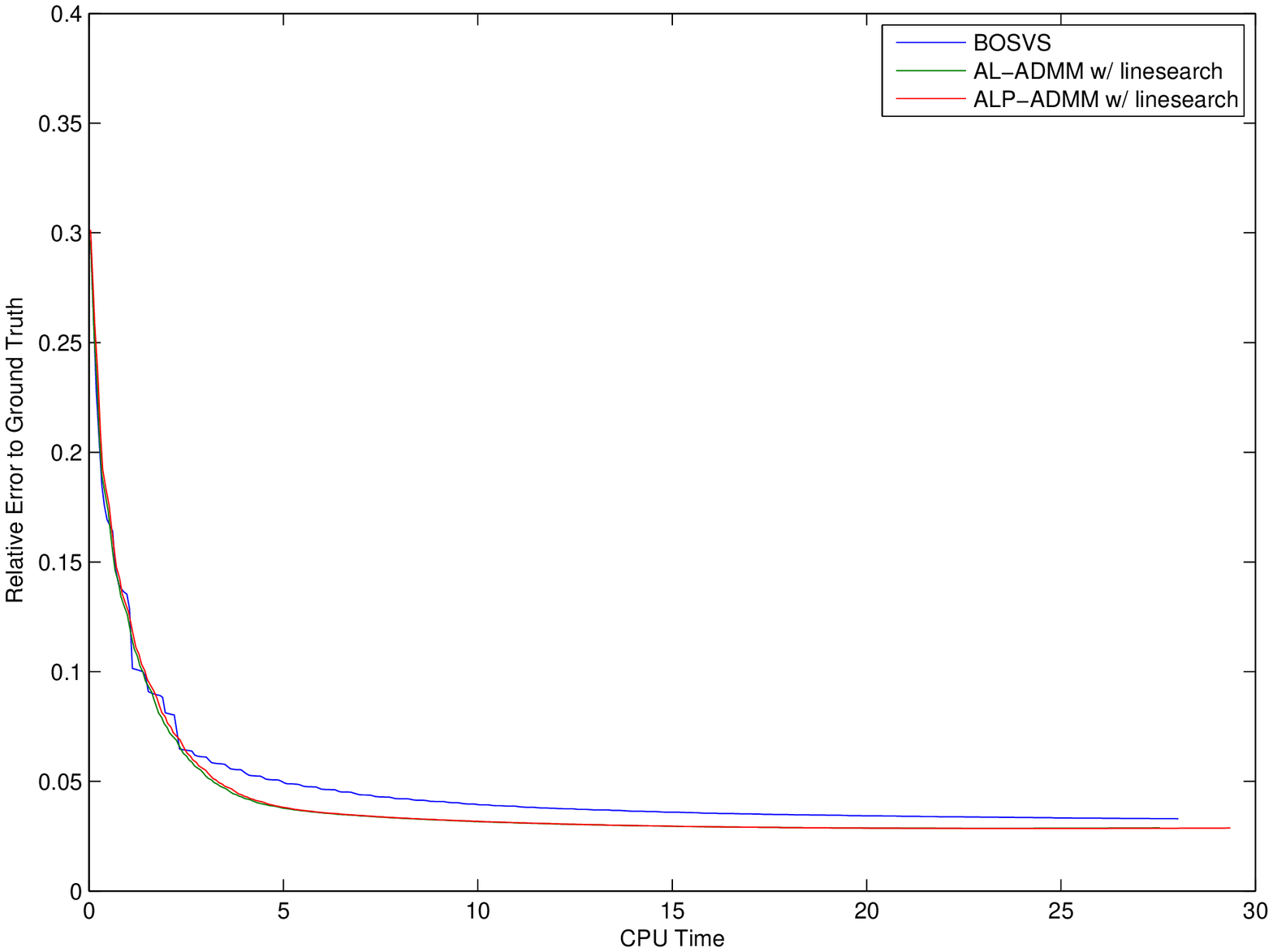}
	\end{subfigure}
	\caption{\label{figPPI} \footnotesize
		Comparisons of AL-ADMM, ALP-ADMM and BOSVS in partially parallel image reconstruction. From top to bottom: performances of algorithms in instances 1a and 1b. Left: the objective function values vs. CPU time. The straight line at the bottom is $f(x_{true})$. Right: the relative errors  vs. CPU time.}
\end{figure}
\begin{figure}
	\ContinuedFloat
	\begin{subfigure}[b]{.49\linewidth}
		\includegraphics[width=\linewidth]{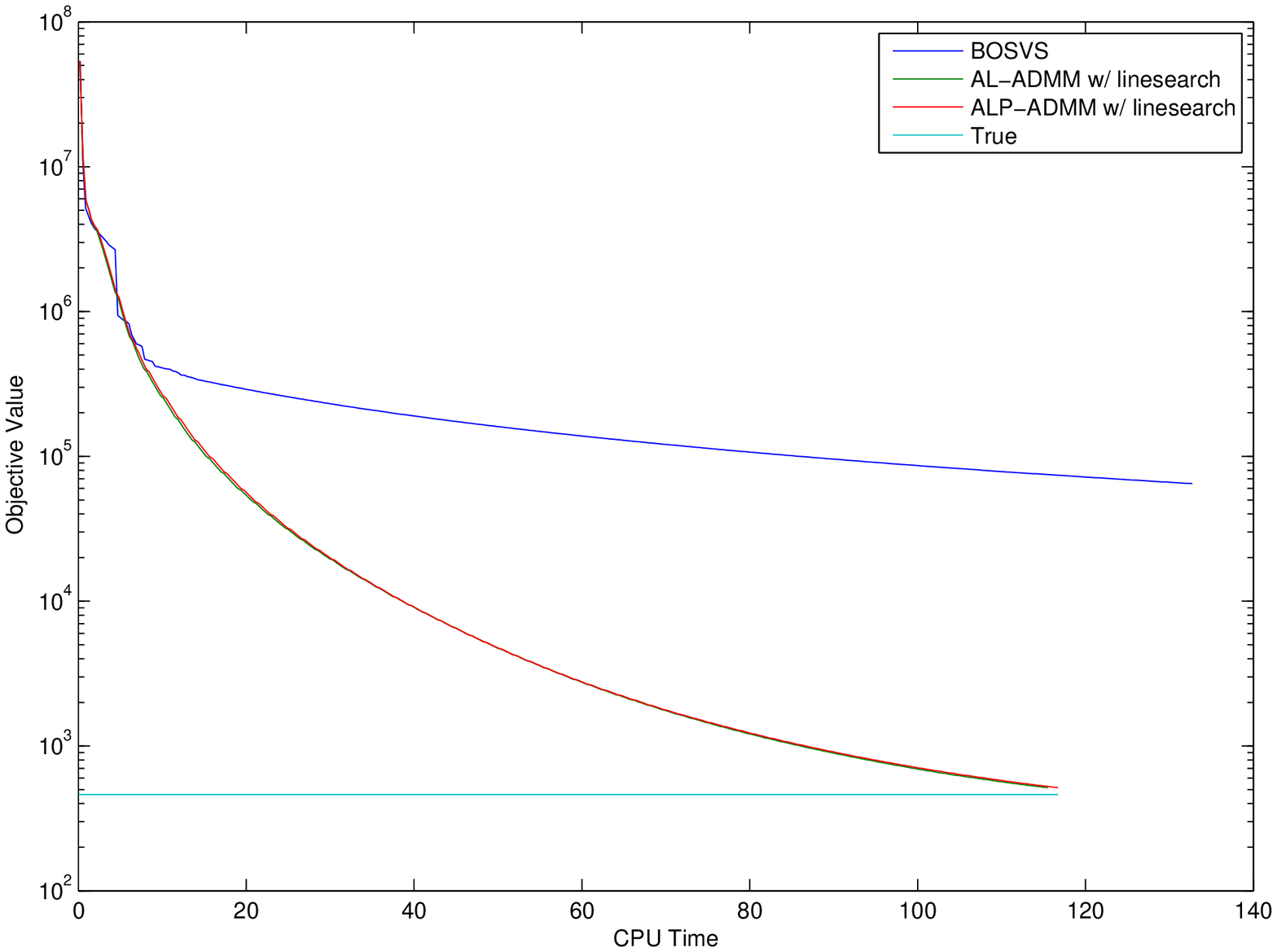}
	\end{subfigure}
	\begin{subfigure}[b]{.49\linewidth}
		\includegraphics[width=\linewidth]{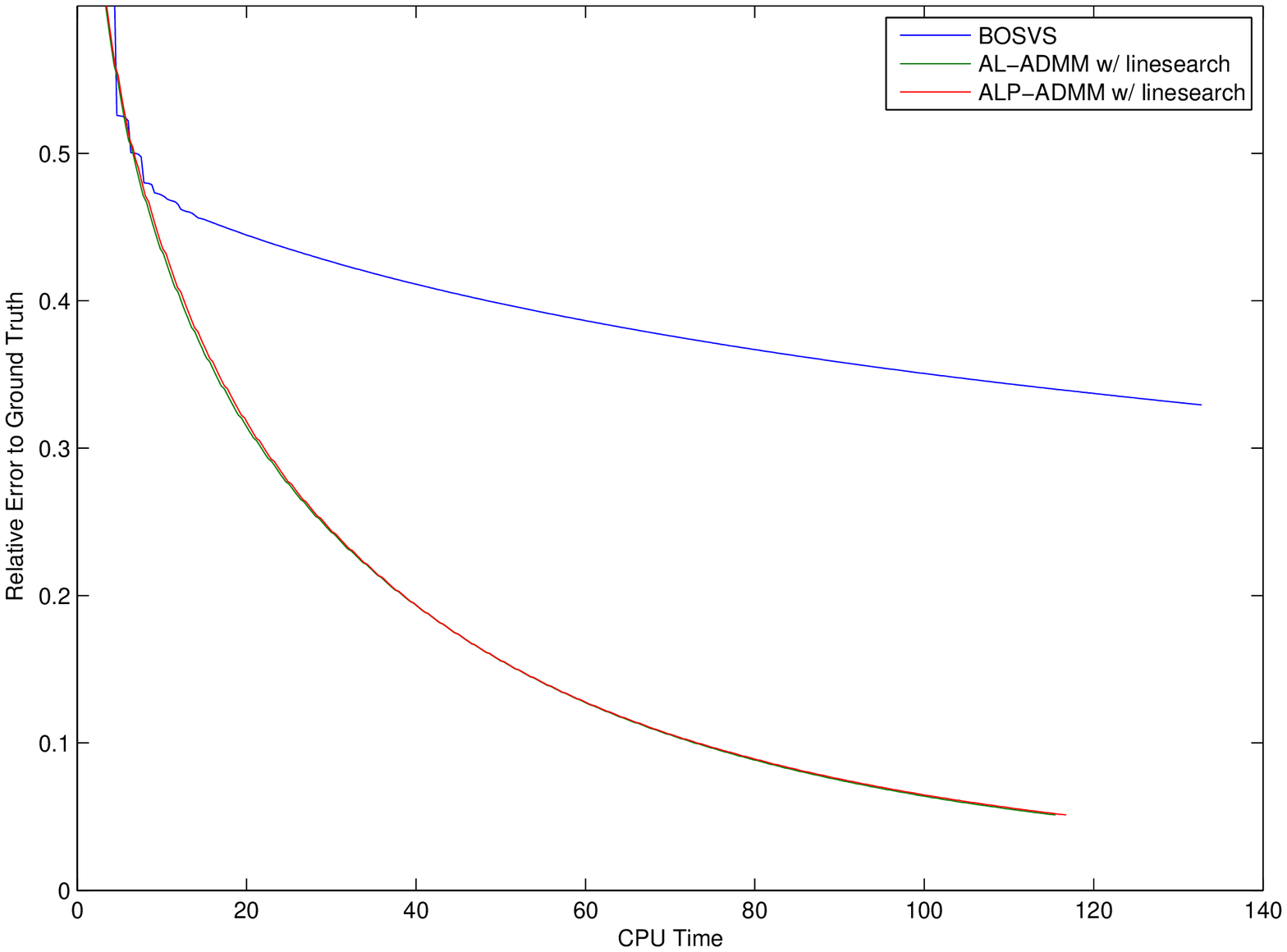}
	\end{subfigure}
	\begin{subfigure}[b]{.49\linewidth}
		\includegraphics[width=\linewidth]{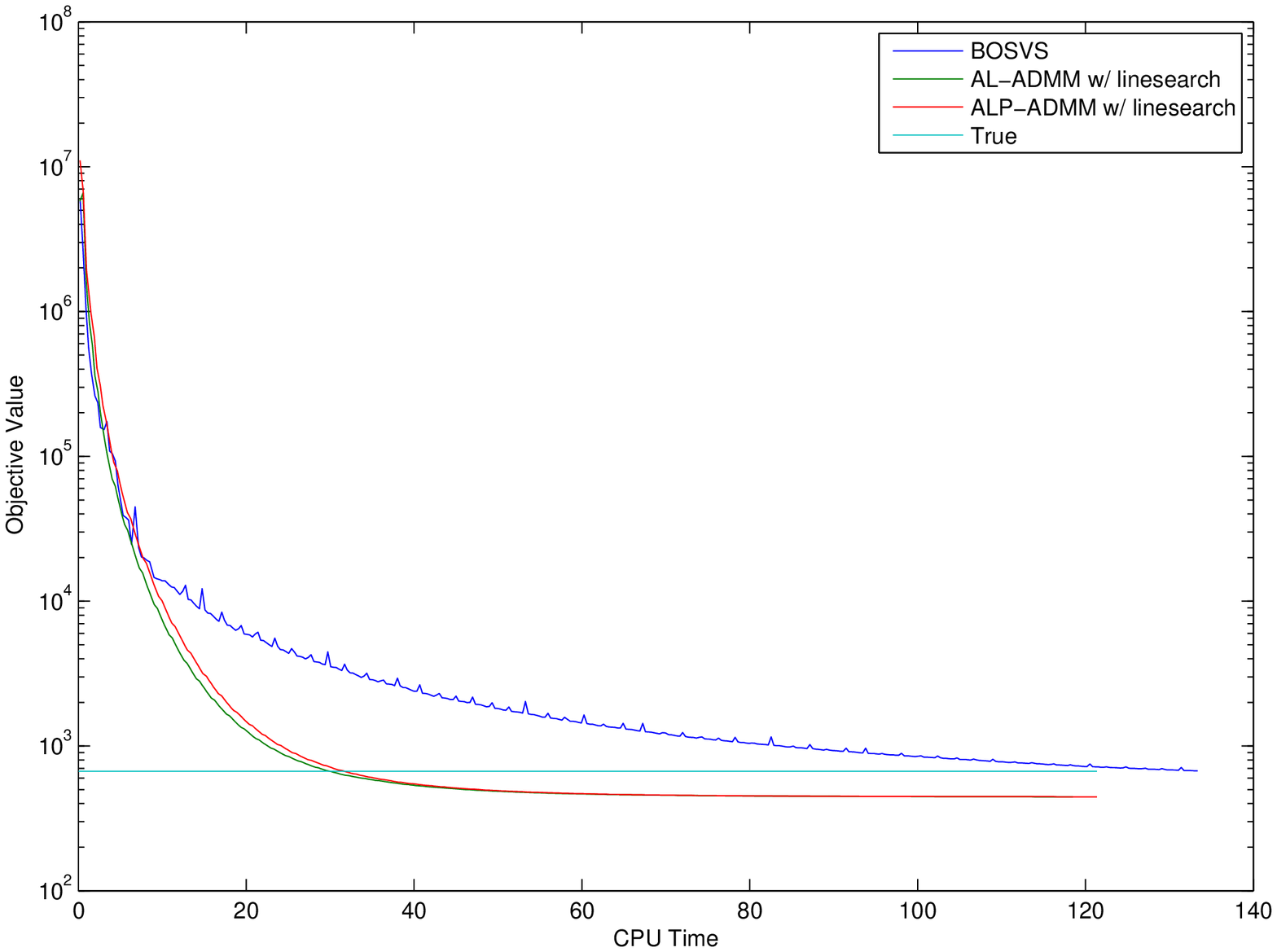}
	\end{subfigure}
	\begin{subfigure}[b]{.49\linewidth}
		\includegraphics[width=\linewidth]{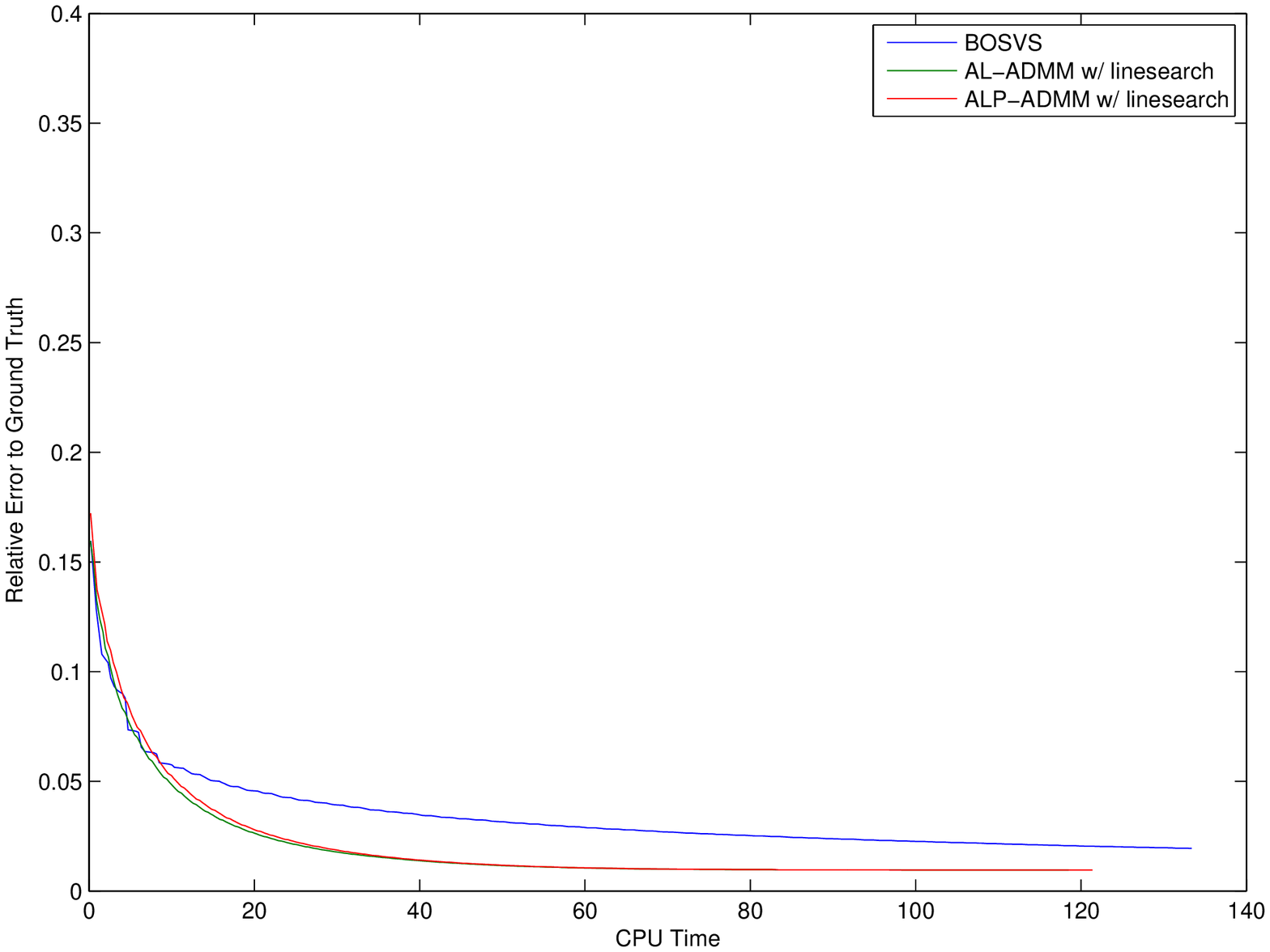}
	\end{subfigure}
	\caption{\footnotesize
	Comparisons of AL-ADMM, ALP-ADMM and BOSVS in partially parallel image reconstruction (cont'd). From top to bottom: performances of algorithms in instances 2a and 2b. Left: the objective function values vs. CPU time. The straight line at the bottom is $f(x_{true})$. Right: the relative errors  vs. CPU time.}
\end{figure}

\begin{figure}
	\begin{subfigure}[b]{.3\linewidth}
		\includegraphics[width=\linewidth]{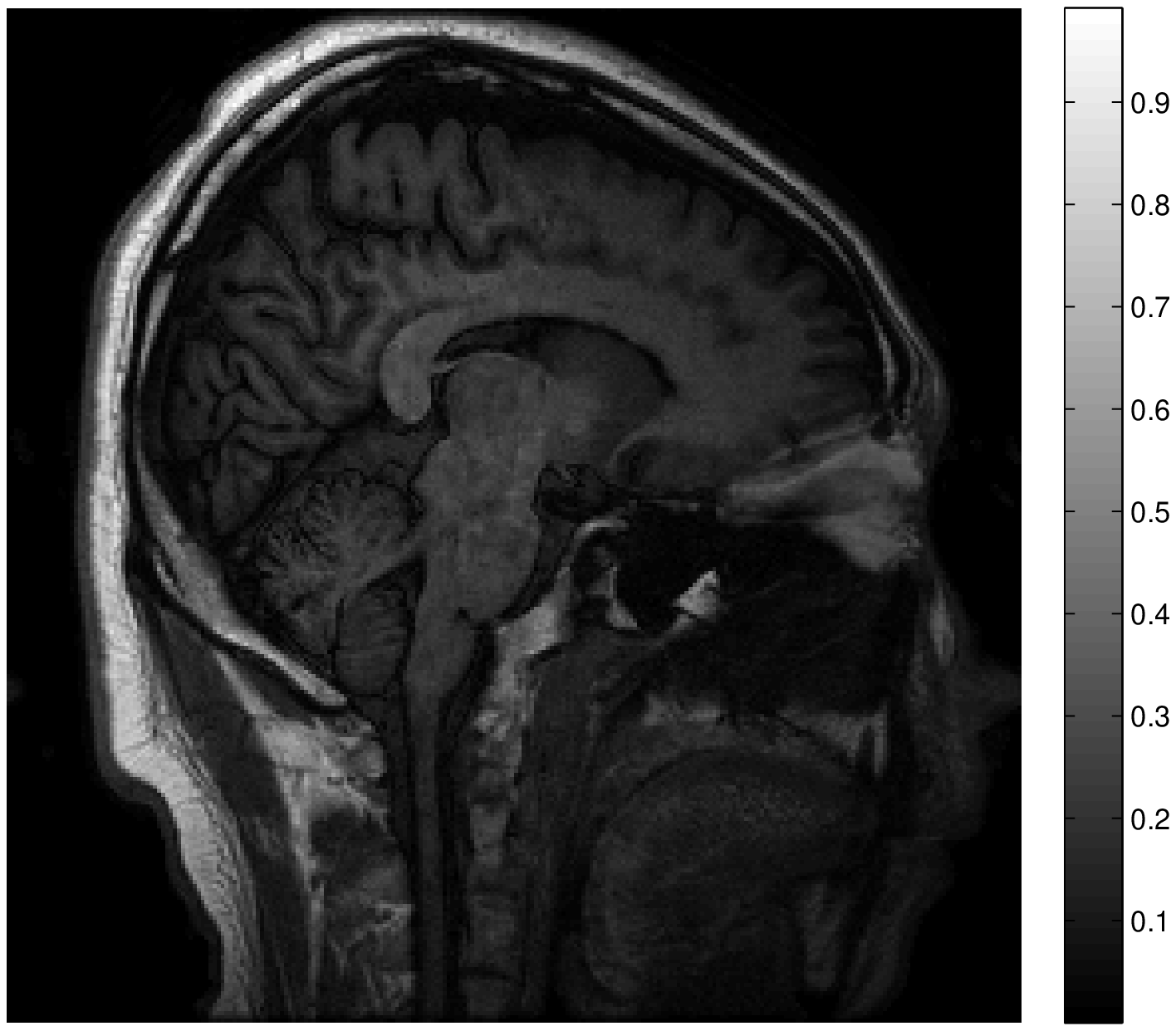}
	\end{subfigure}
	\begin{subfigure}[b]{.3\linewidth}
		\includegraphics[width=\linewidth]{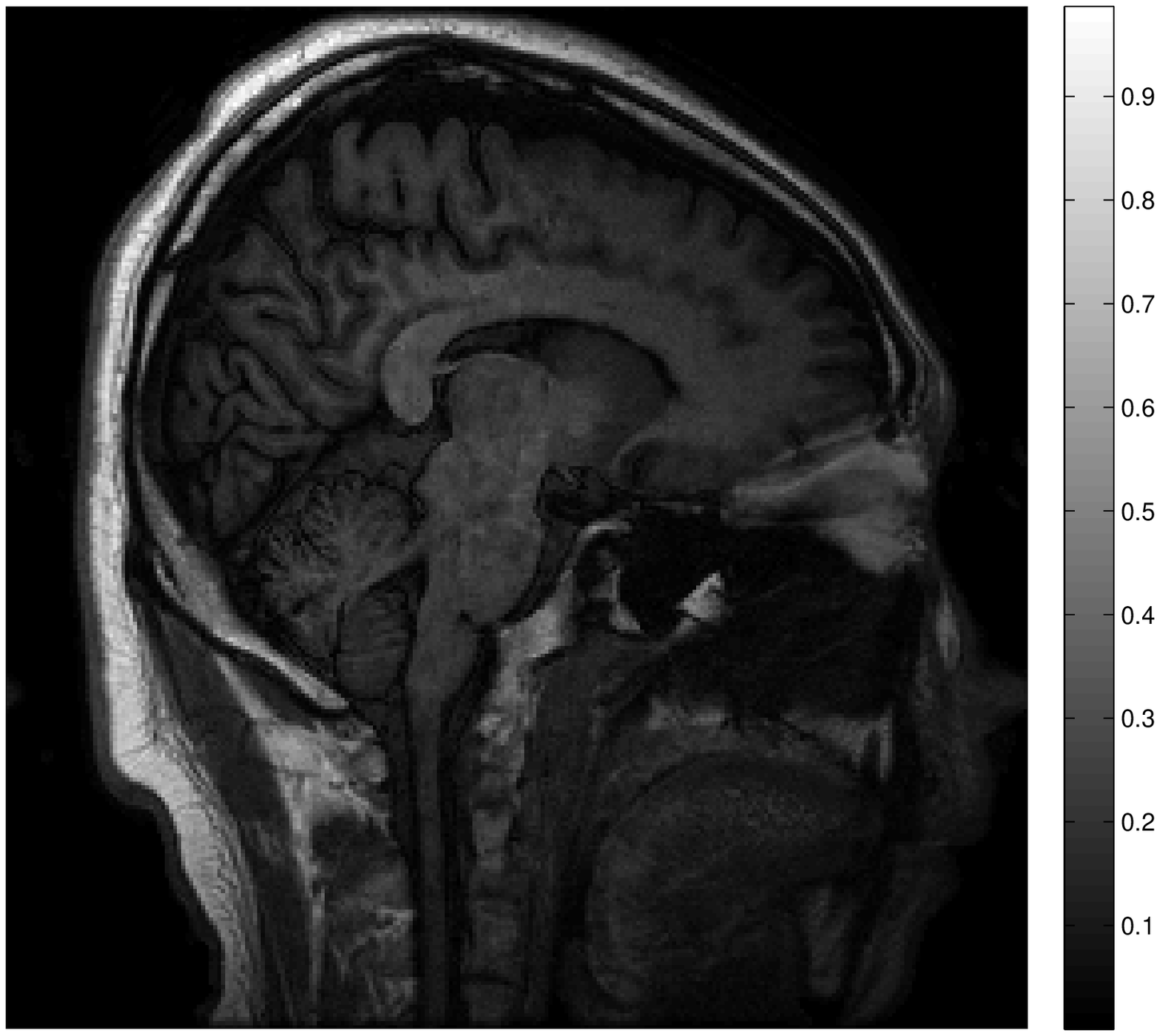}
	\end{subfigure}
	\begin{subfigure}[b]{.3\linewidth}
		\includegraphics[width=\linewidth]{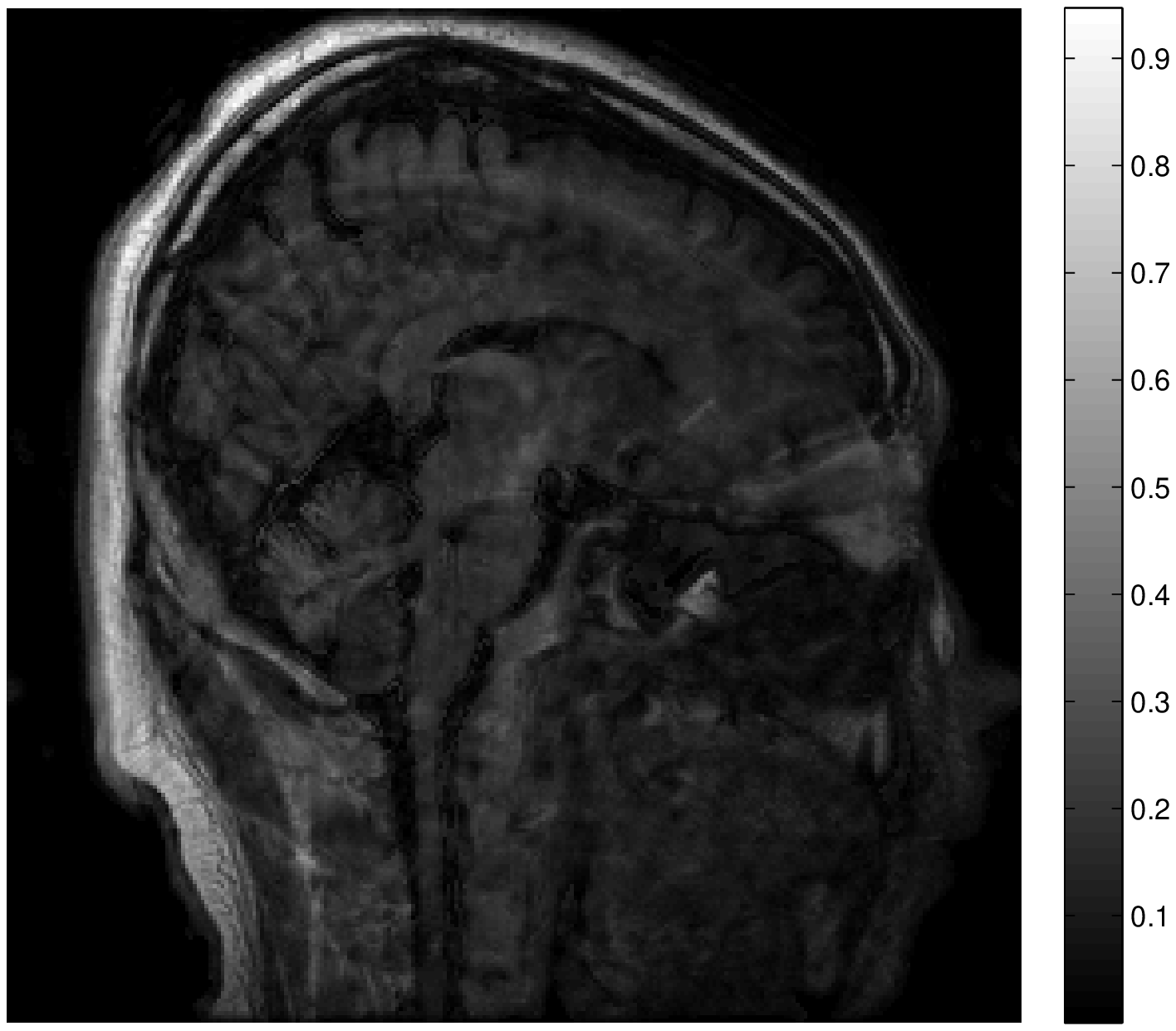}
	\end{subfigure}
	\\
	\begin{subfigure}[b]{.3\linewidth}
		\includegraphics[width=\linewidth]{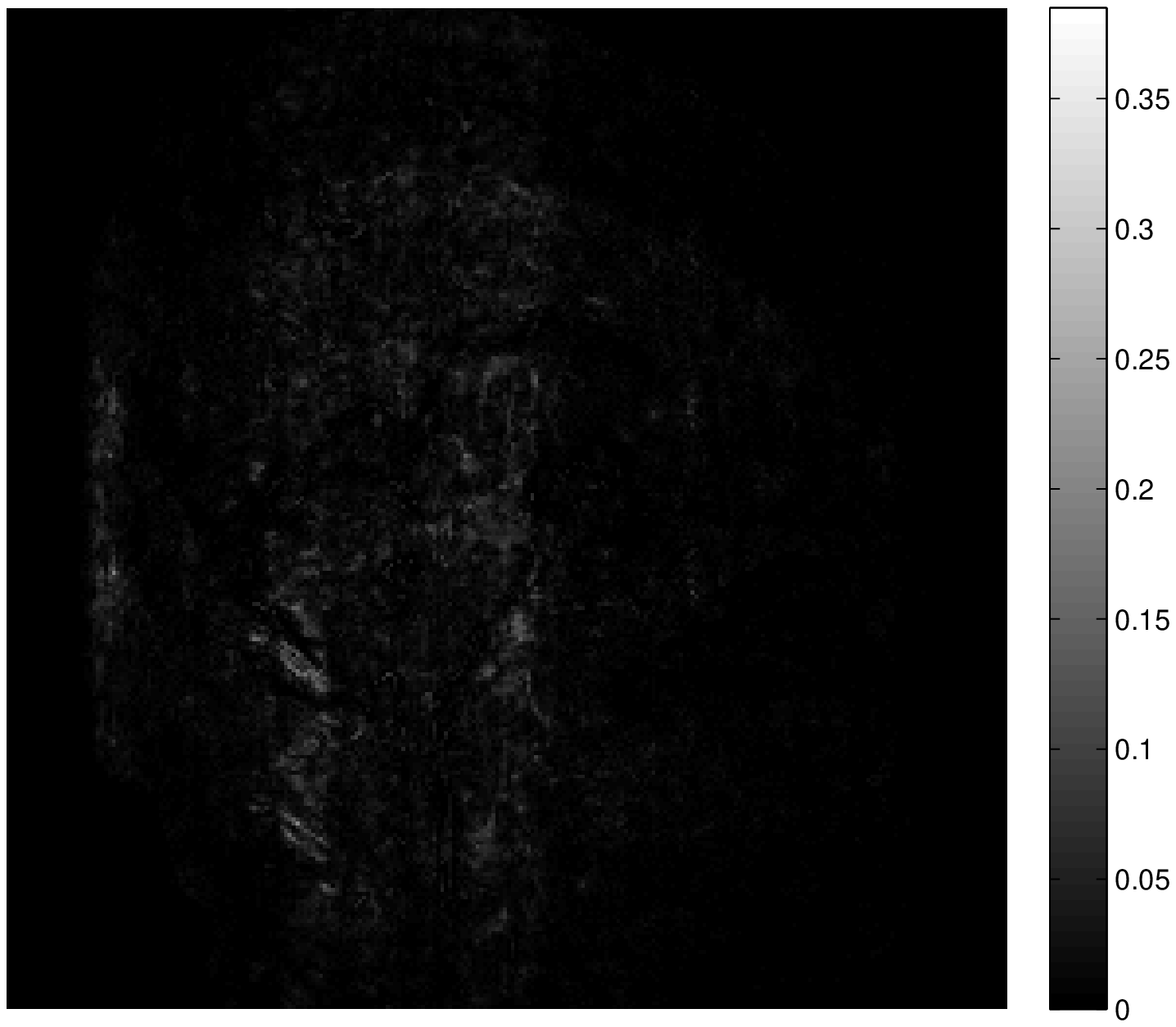}
	\end{subfigure}
	\begin{subfigure}[b]{.3\linewidth}
		\includegraphics[width=\linewidth]{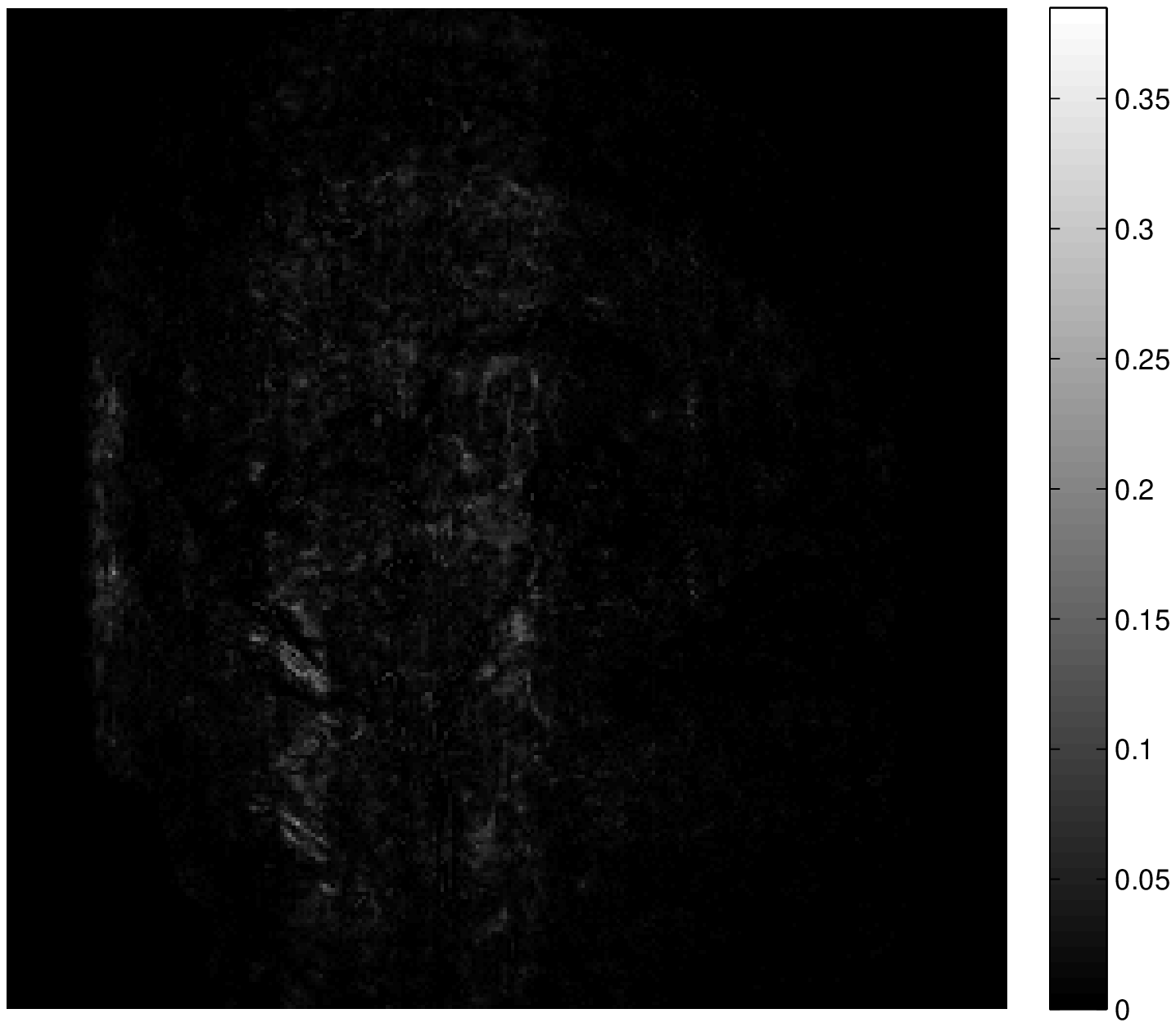}
	\end{subfigure}
	\begin{subfigure}[b]{.3\linewidth}
		\includegraphics[width=\linewidth]{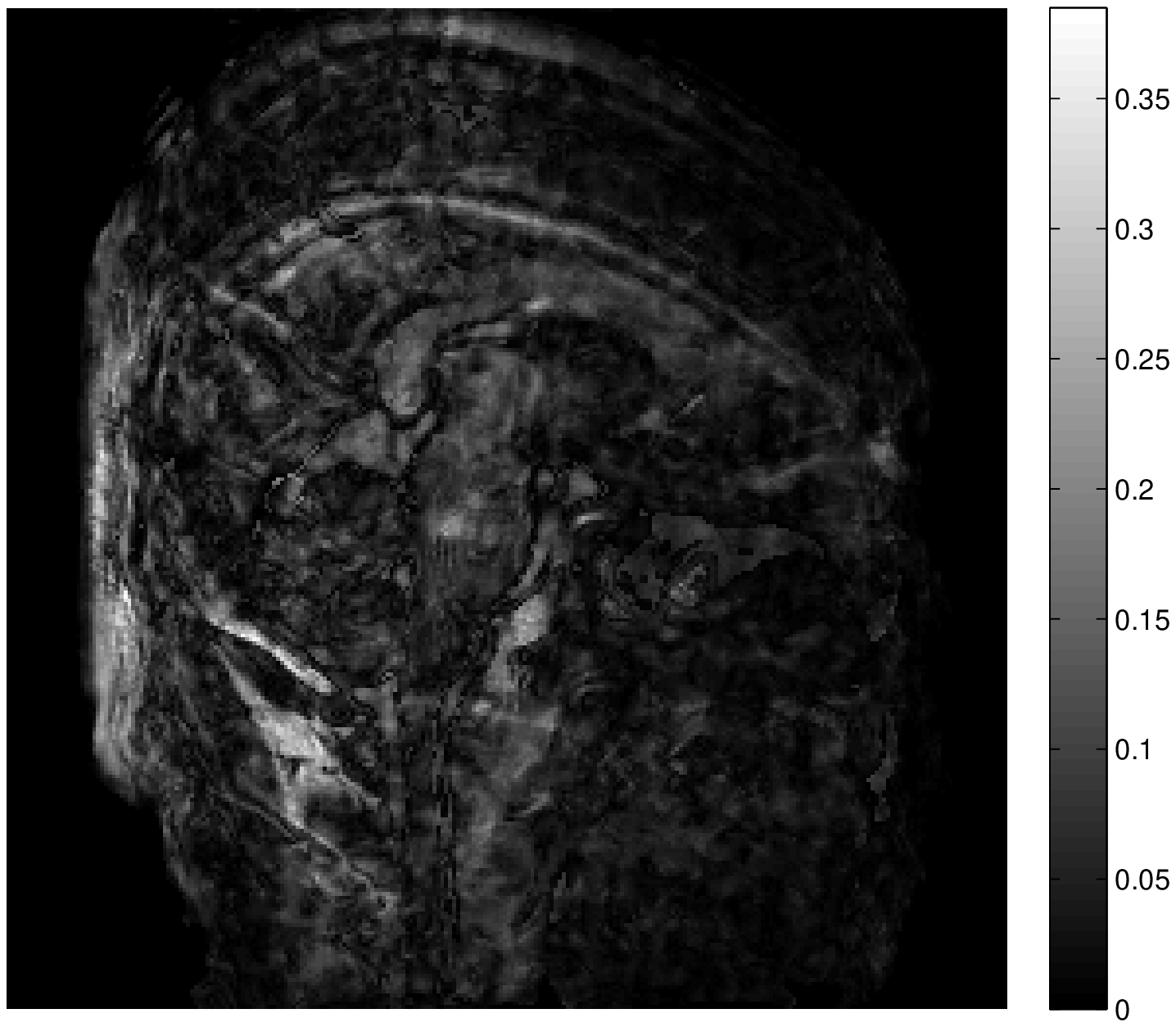}
	\end{subfigure}
	\\
	\begin{subfigure}[b]{.3\linewidth}
		\includegraphics[width=\linewidth]{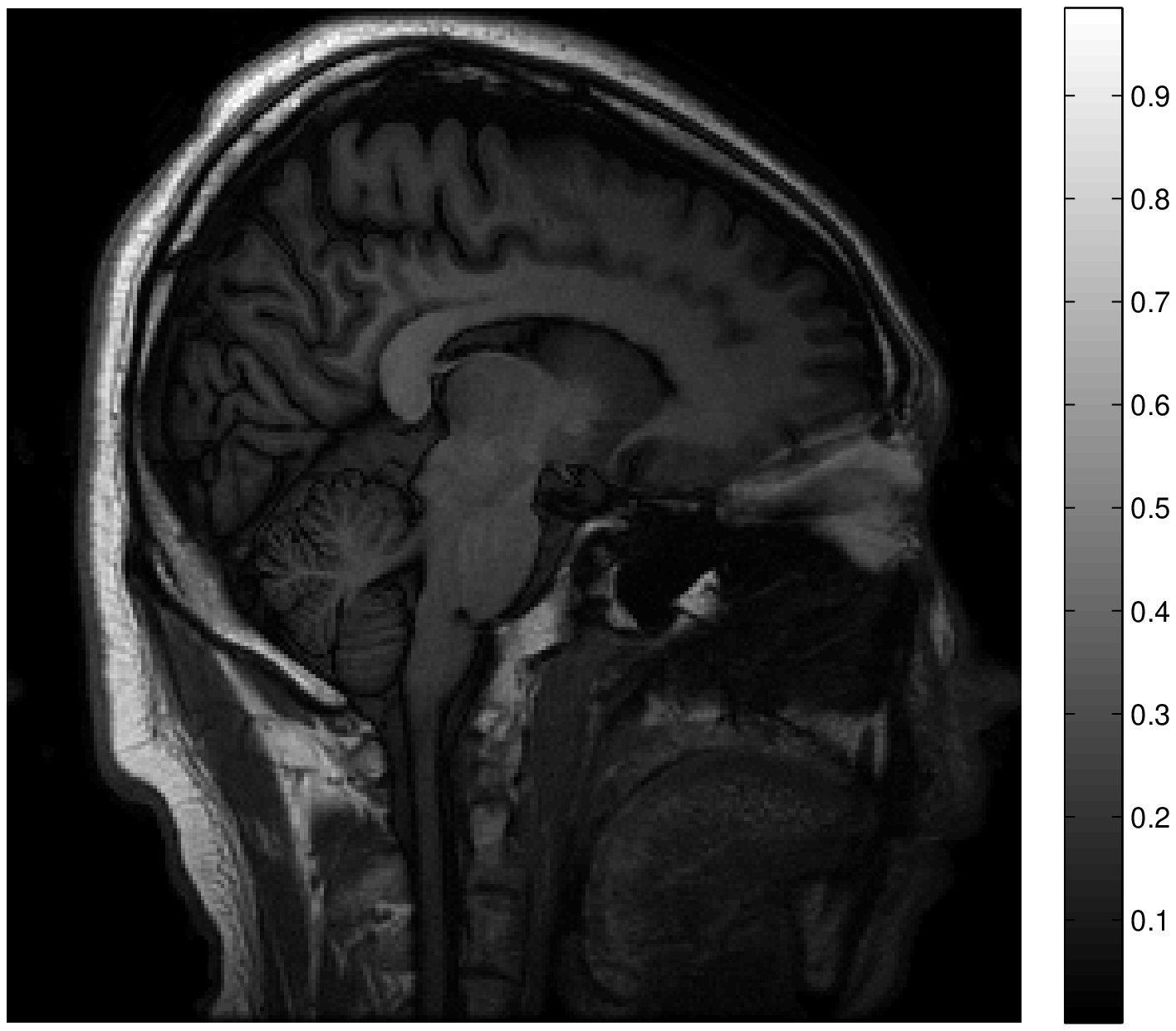}
	\end{subfigure}
	\begin{subfigure}[b]{.3\linewidth}
		\includegraphics[width=\linewidth]{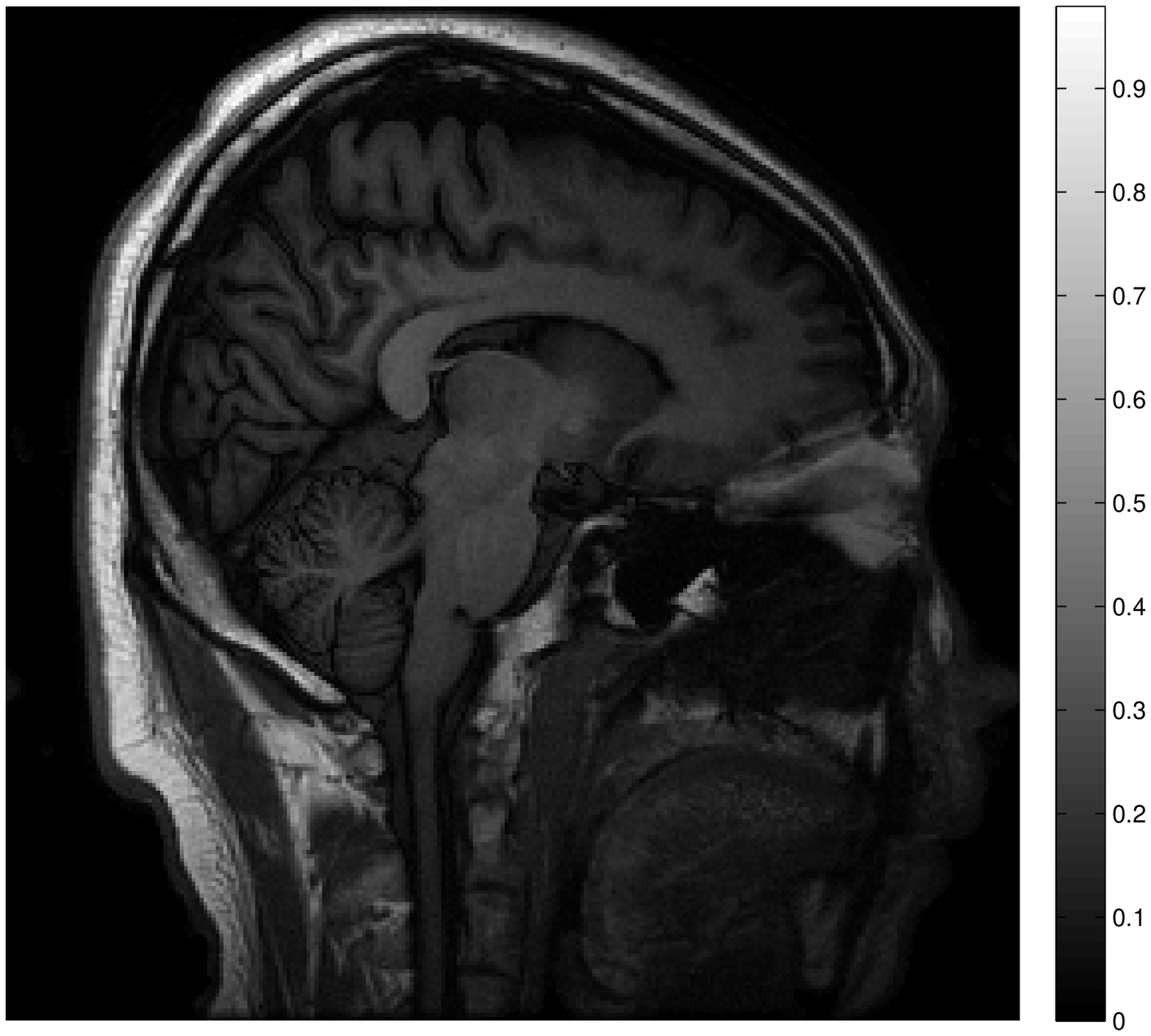}
	\end{subfigure}
	\begin{subfigure}[b]{.3\linewidth}
		\includegraphics[width=\linewidth]{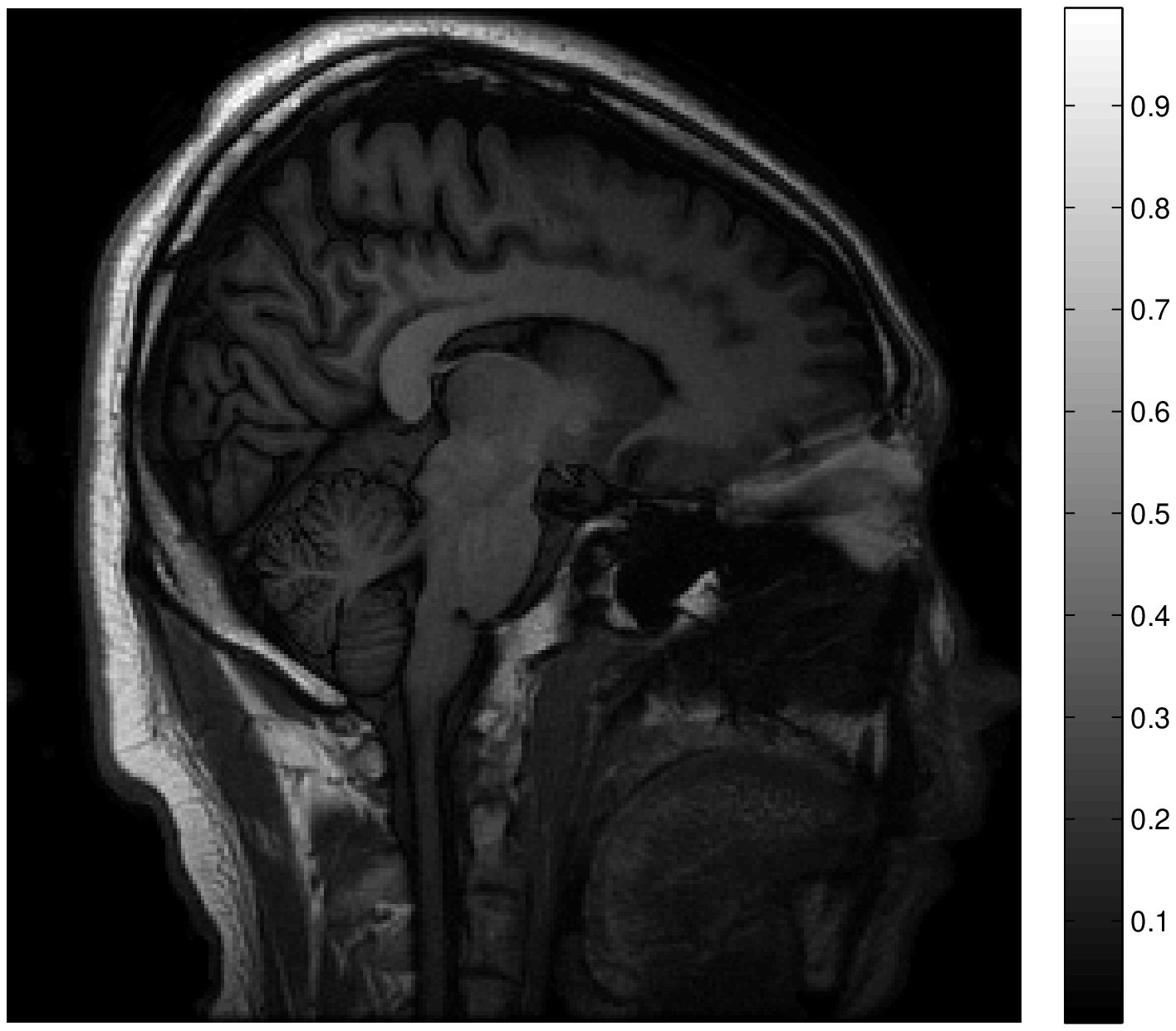}
	\end{subfigure}
	\\
	\begin{subfigure}[b]{.3\linewidth}
		\includegraphics[width=\linewidth]{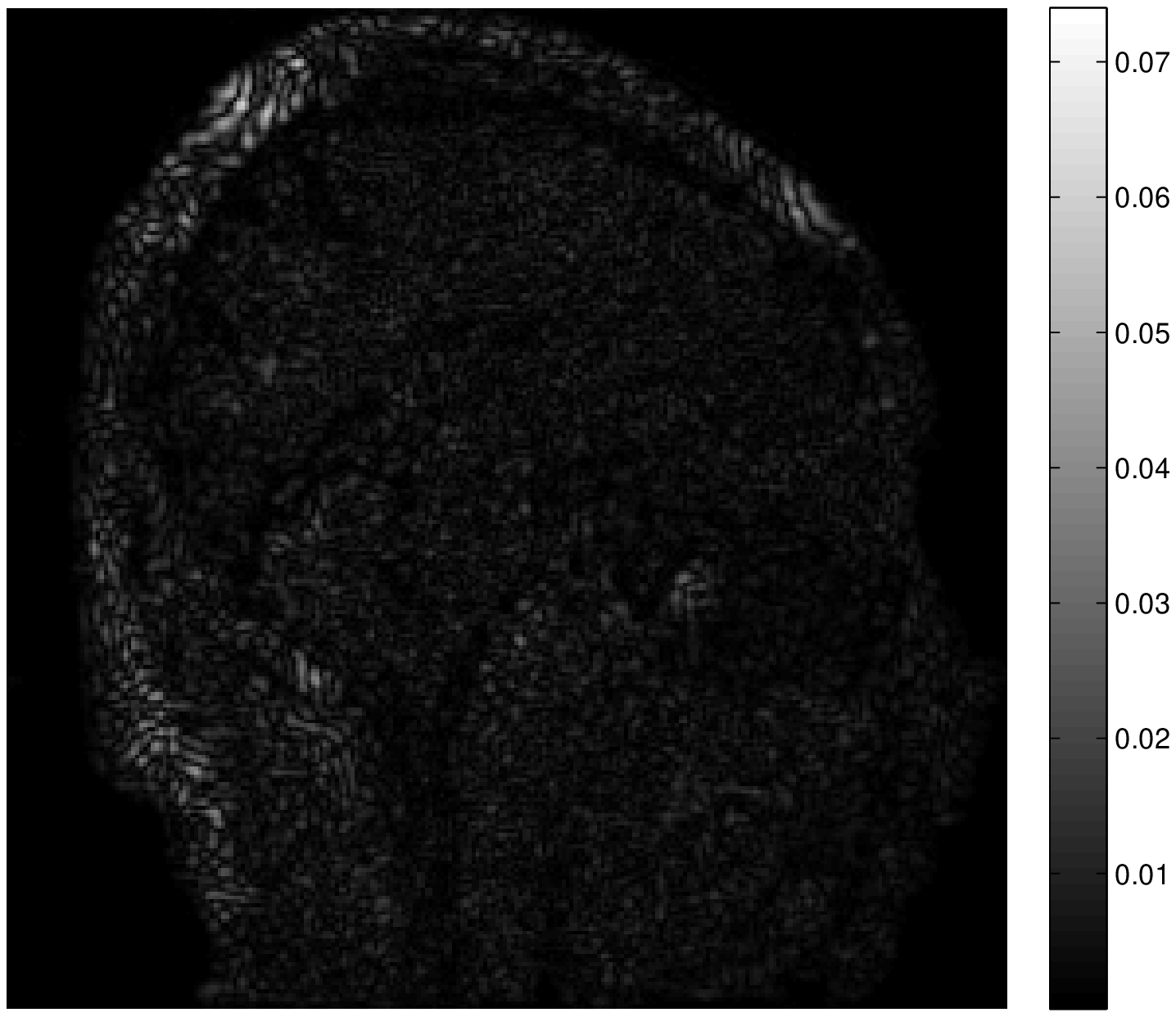}
	\end{subfigure}
	\begin{subfigure}[b]{.3\linewidth}
		\includegraphics[width=\linewidth]{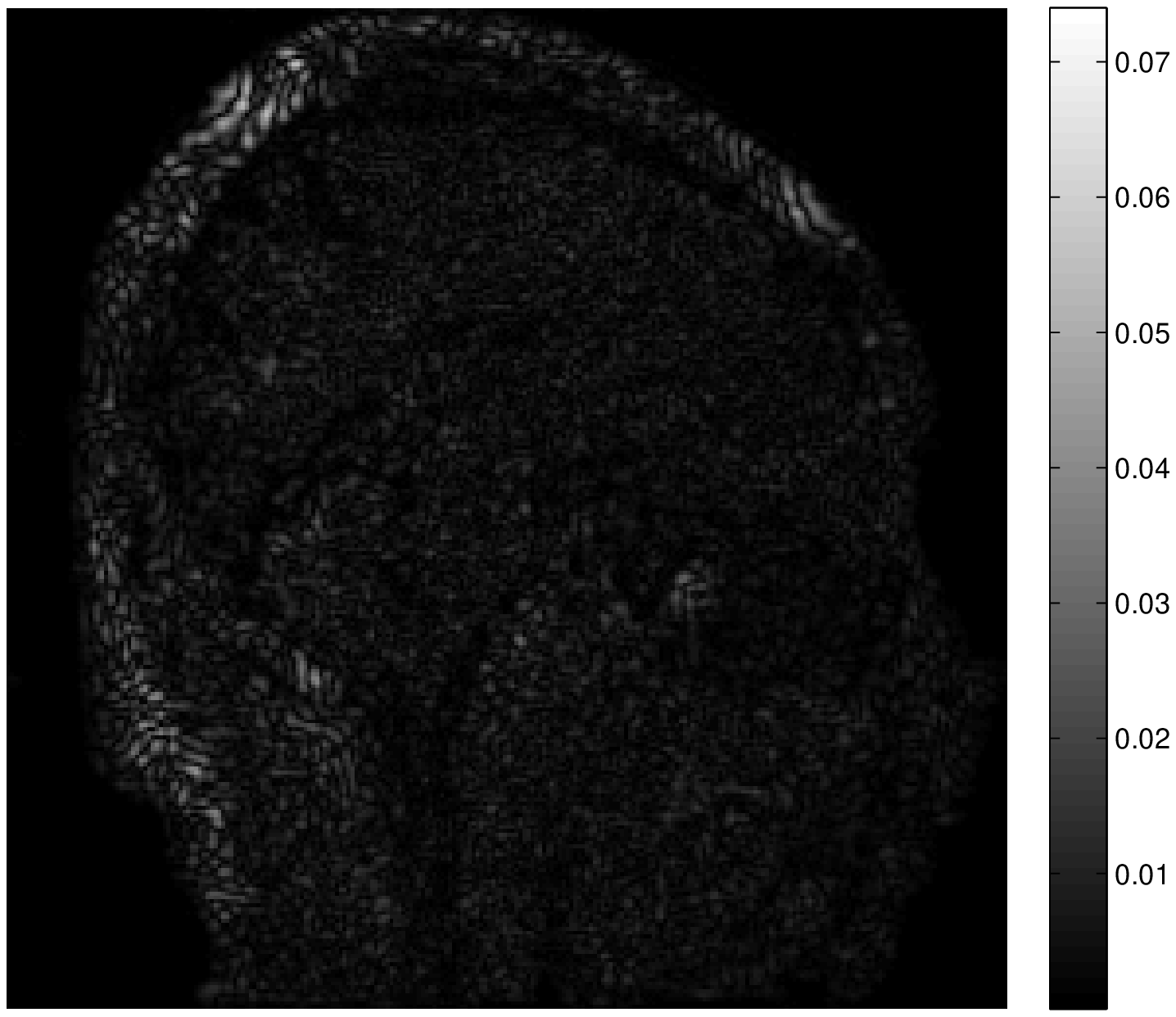}
	\end{subfigure}
	\begin{subfigure}[b]{.3\linewidth}
		\includegraphics[width=\linewidth]{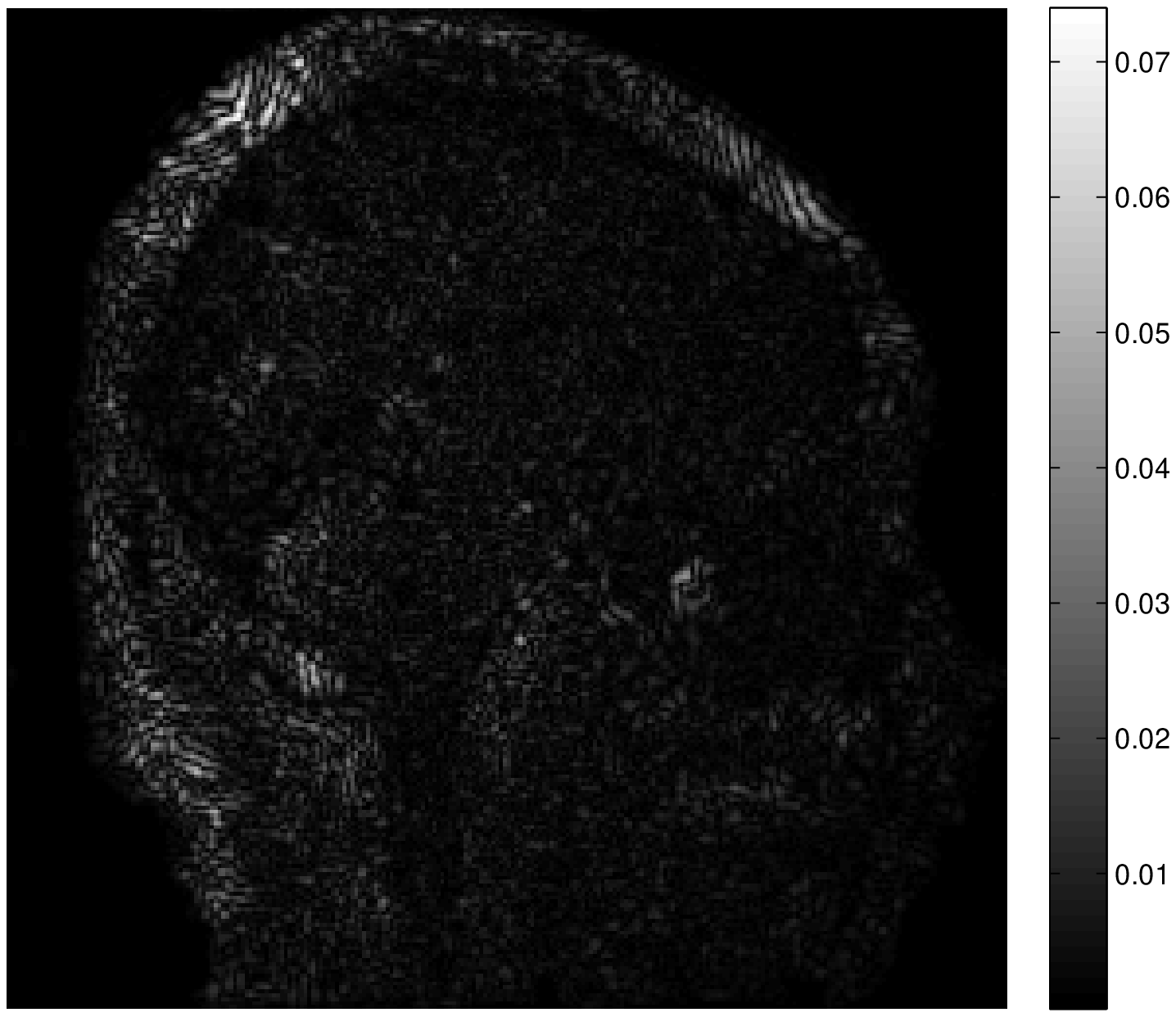}
	\end{subfigure}
	\\
	\caption{\label{figPPIRecon}\footnotesize Comparison of AL-ADMM, ALP-ADMM and BOSVS in partially parallel image reconstruction. From top to bottom: Reconstructed images and reconstruction errors in instances 1a and 1b, respectively. From left to right: AL-ADMM, ALP-ADMM and BOSVS.}
\end{figure}

\begin{figure}
	\ContinuedFloat
	\begin{subfigure}[b]{.3\linewidth}
		\includegraphics[width=\linewidth]{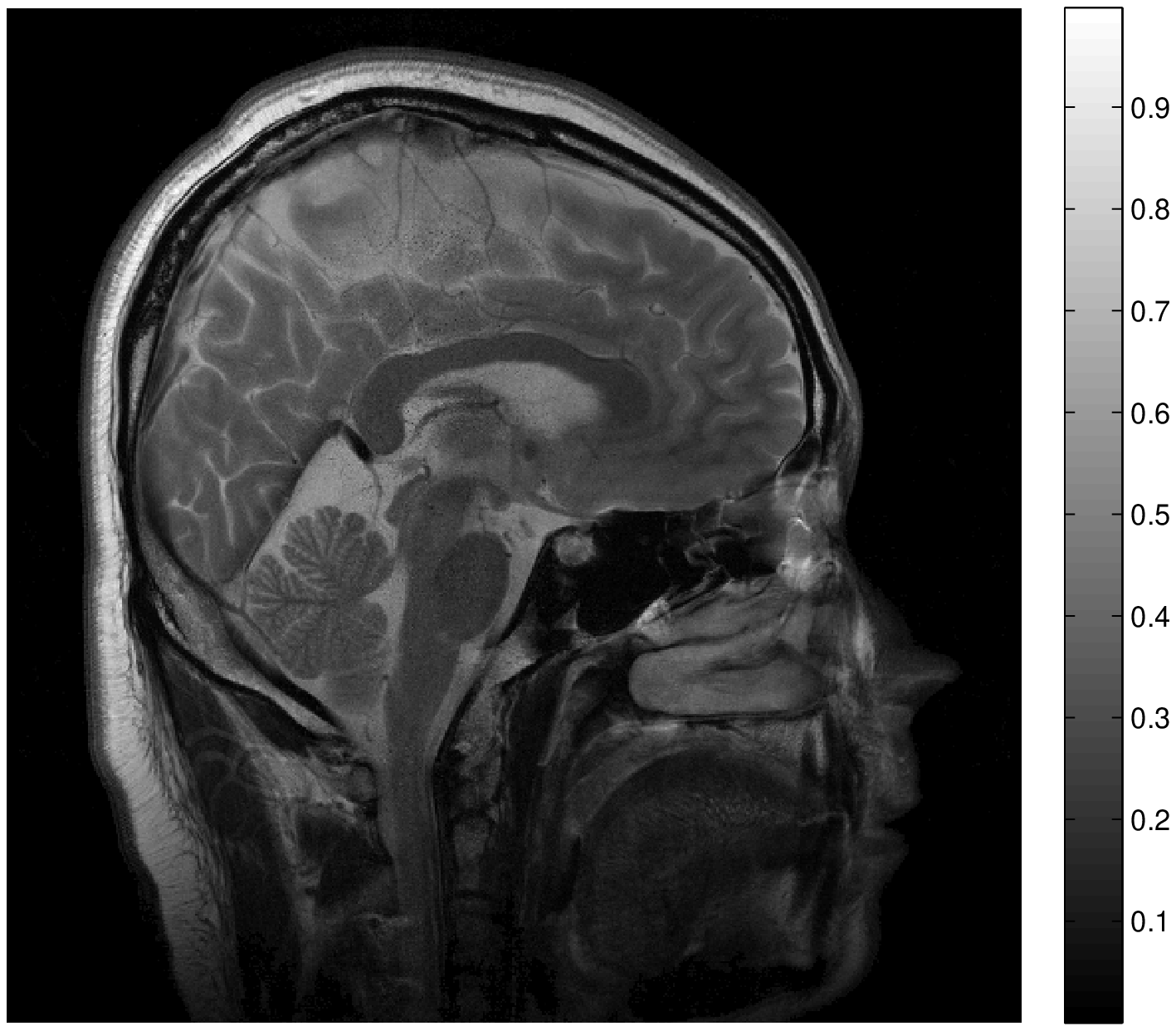}
	\end{subfigure}
	\begin{subfigure}[b]{.3\linewidth}
		\includegraphics[width=\linewidth]{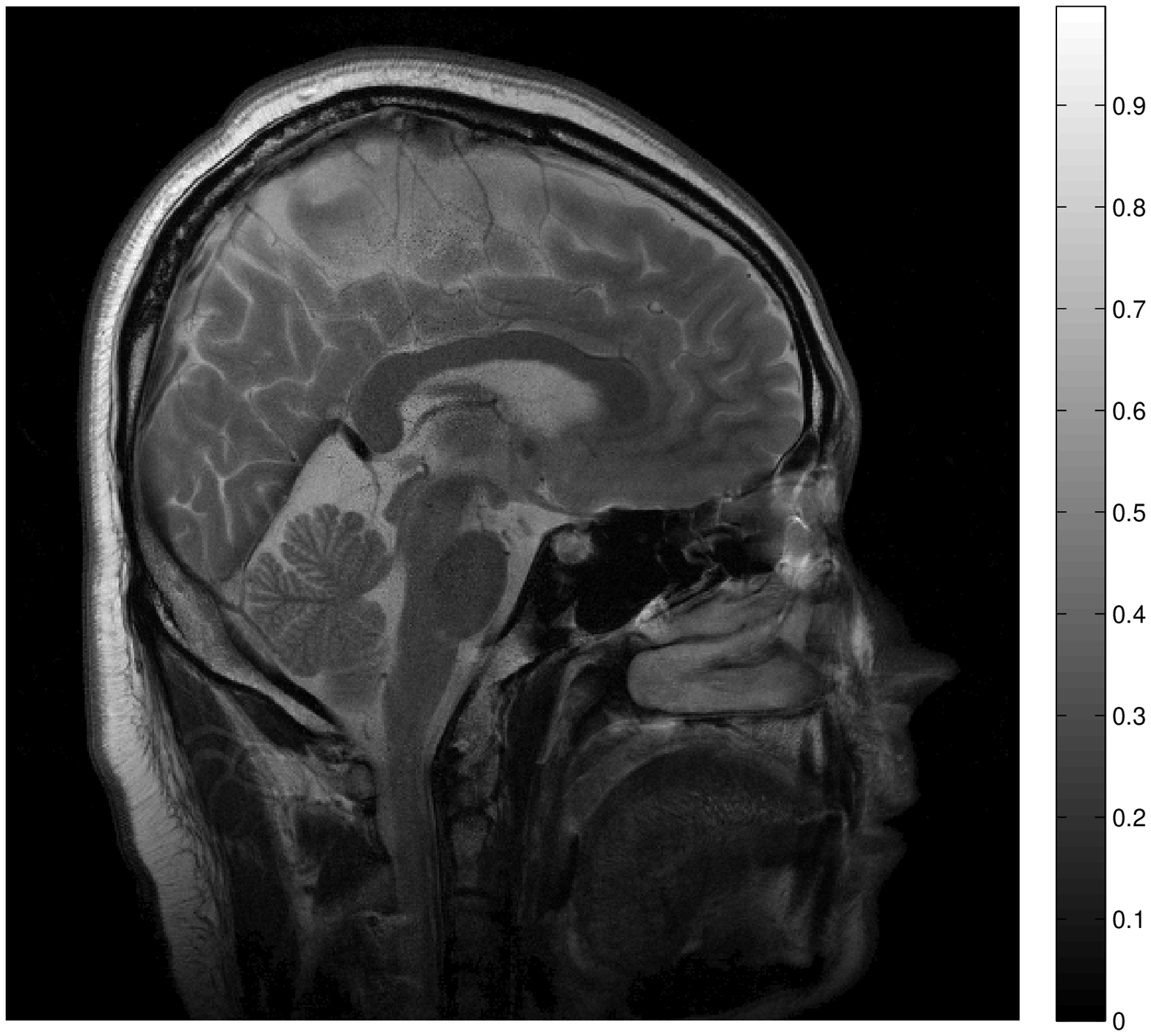}
	\end{subfigure}
	\begin{subfigure}[b]{.3\linewidth}
		\includegraphics[width=\linewidth]{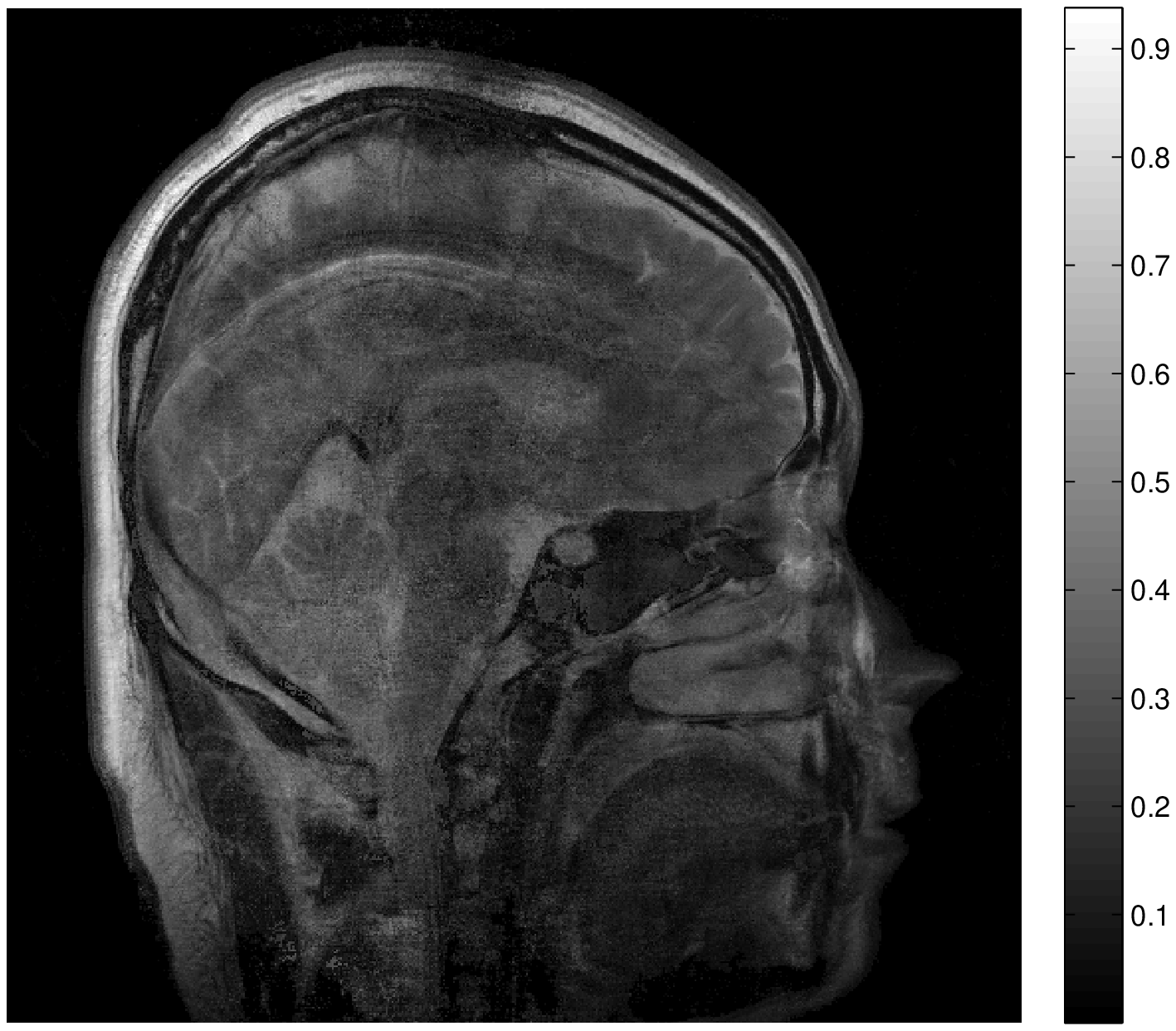}
	\end{subfigure}
	\\
	\begin{subfigure}[b]{.3\linewidth}
		\includegraphics[width=\linewidth]{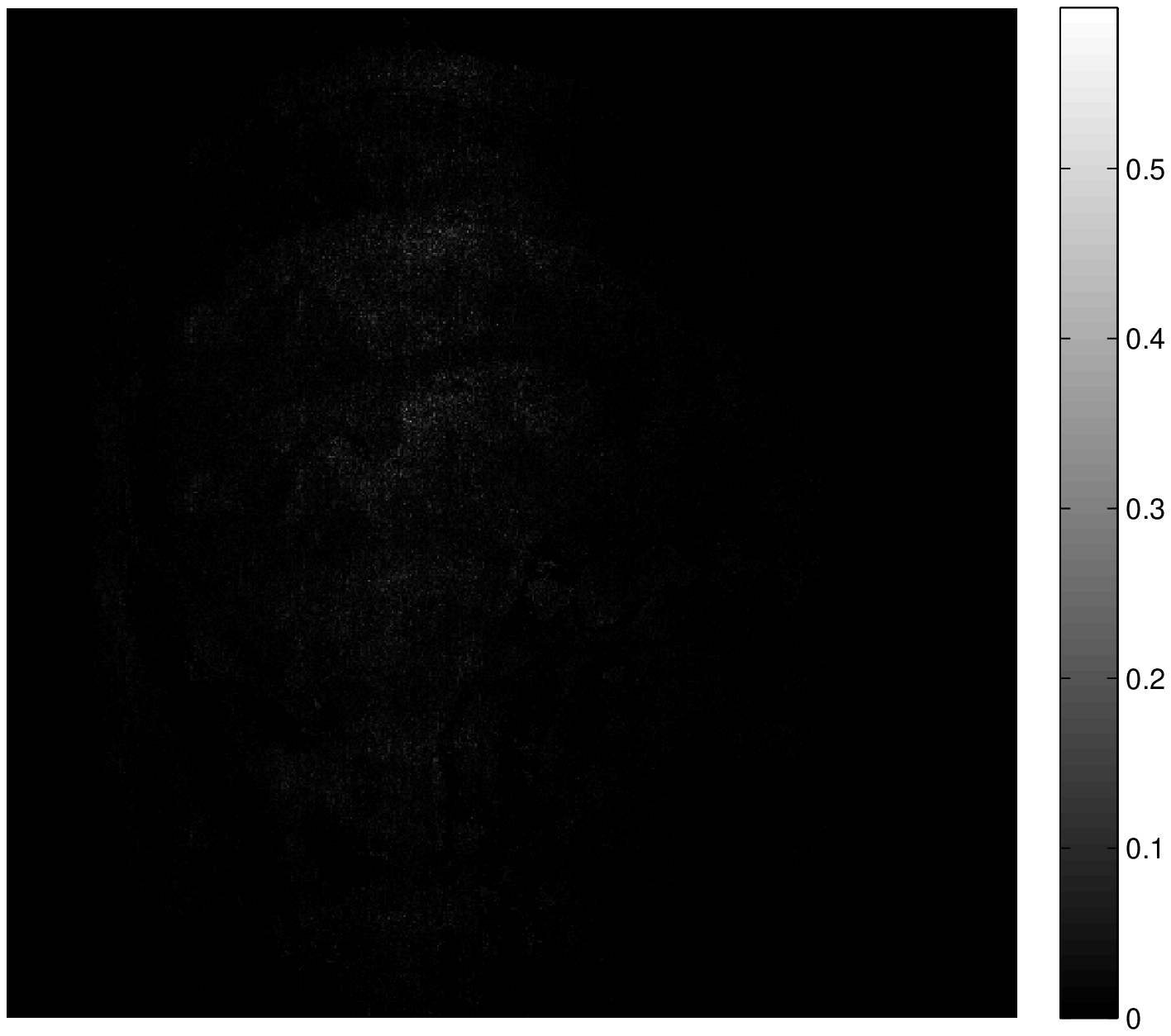}
	\end{subfigure}
	\begin{subfigure}[b]{.3\linewidth}
		\includegraphics[width=\linewidth]{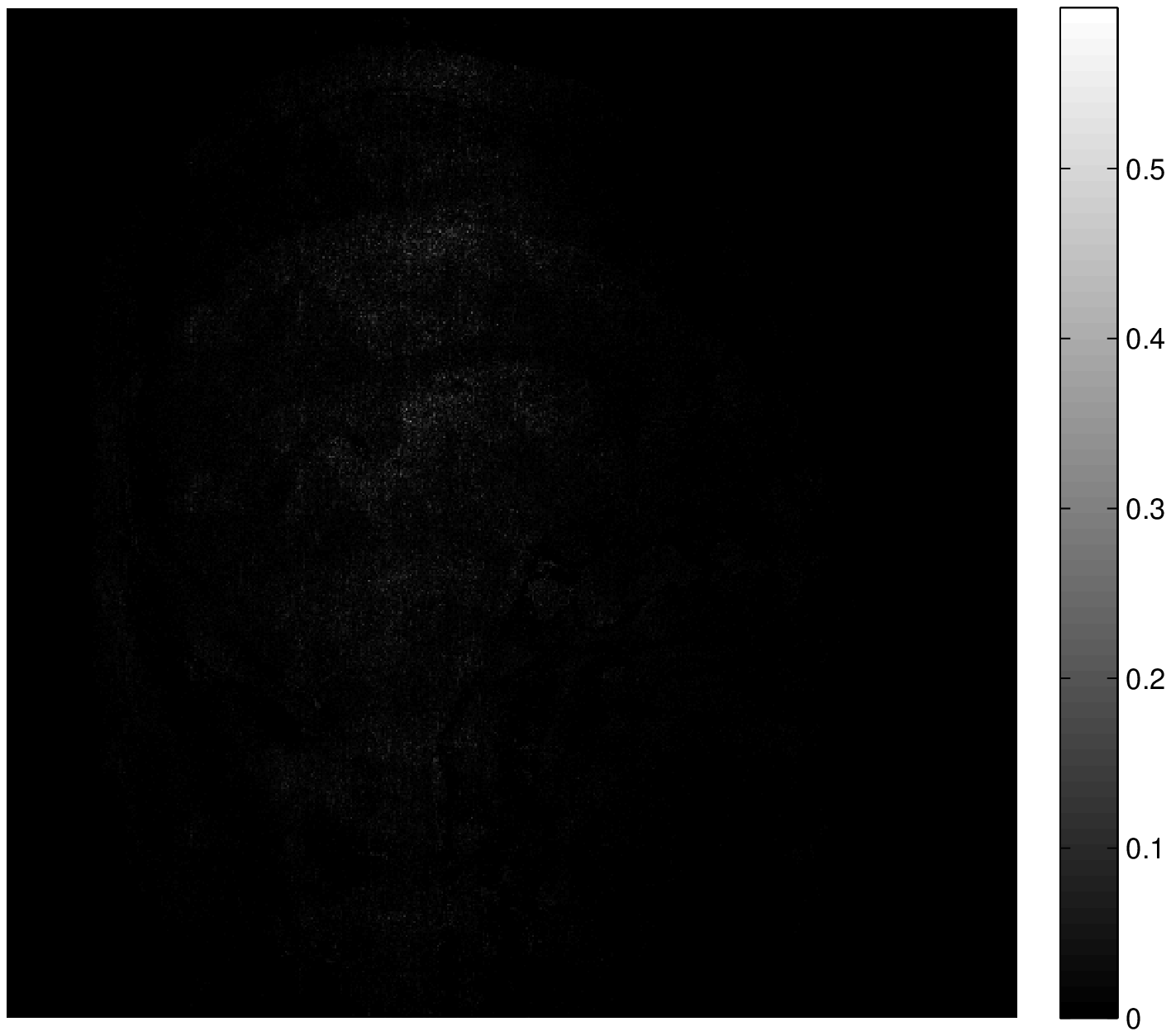}
	\end{subfigure}
	\begin{subfigure}[b]{.3\linewidth}
		\includegraphics[width=\linewidth]{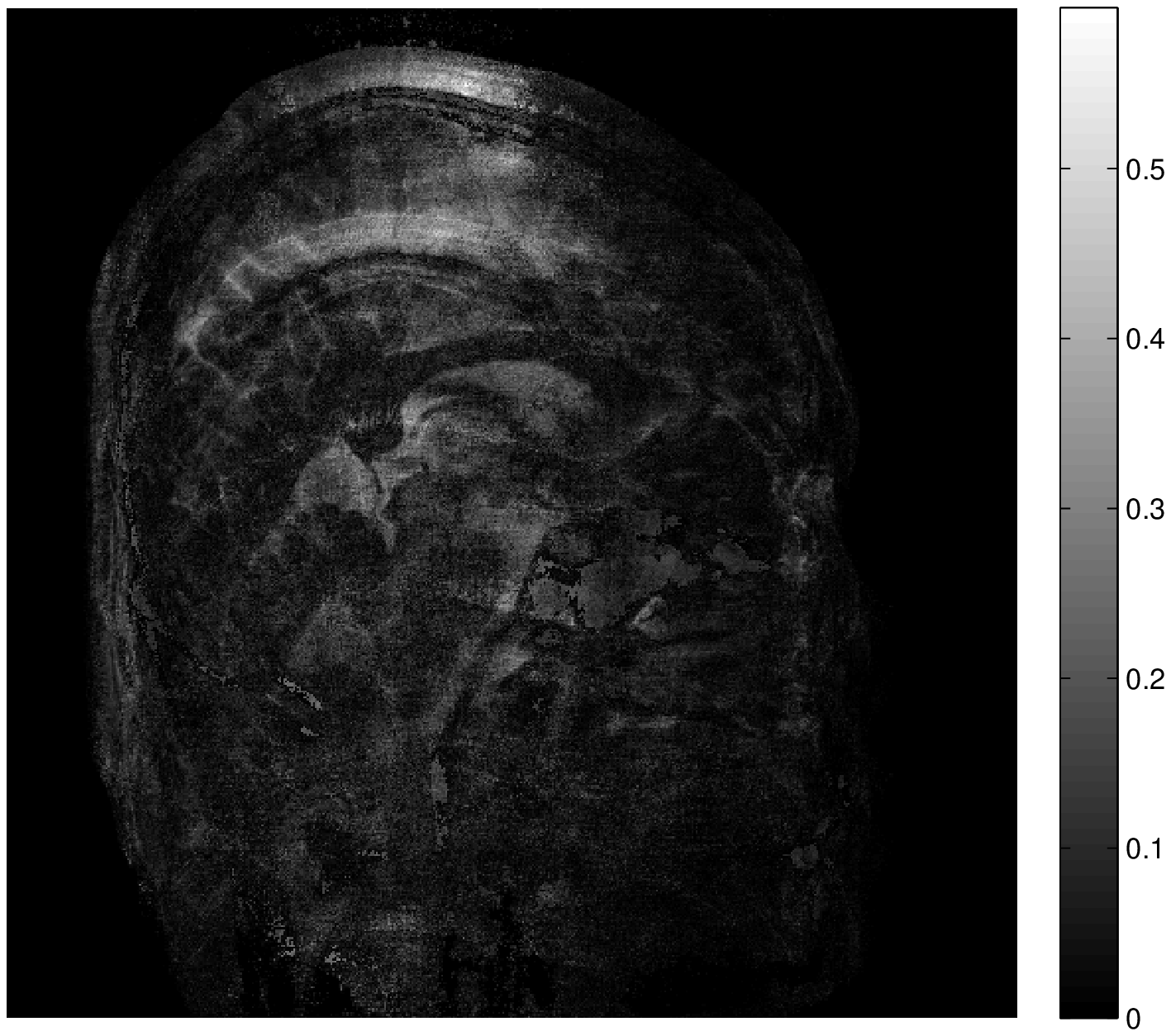}
	\end{subfigure}
	\\
	\begin{subfigure}[b]{.3\linewidth}
		\includegraphics[width=\linewidth]{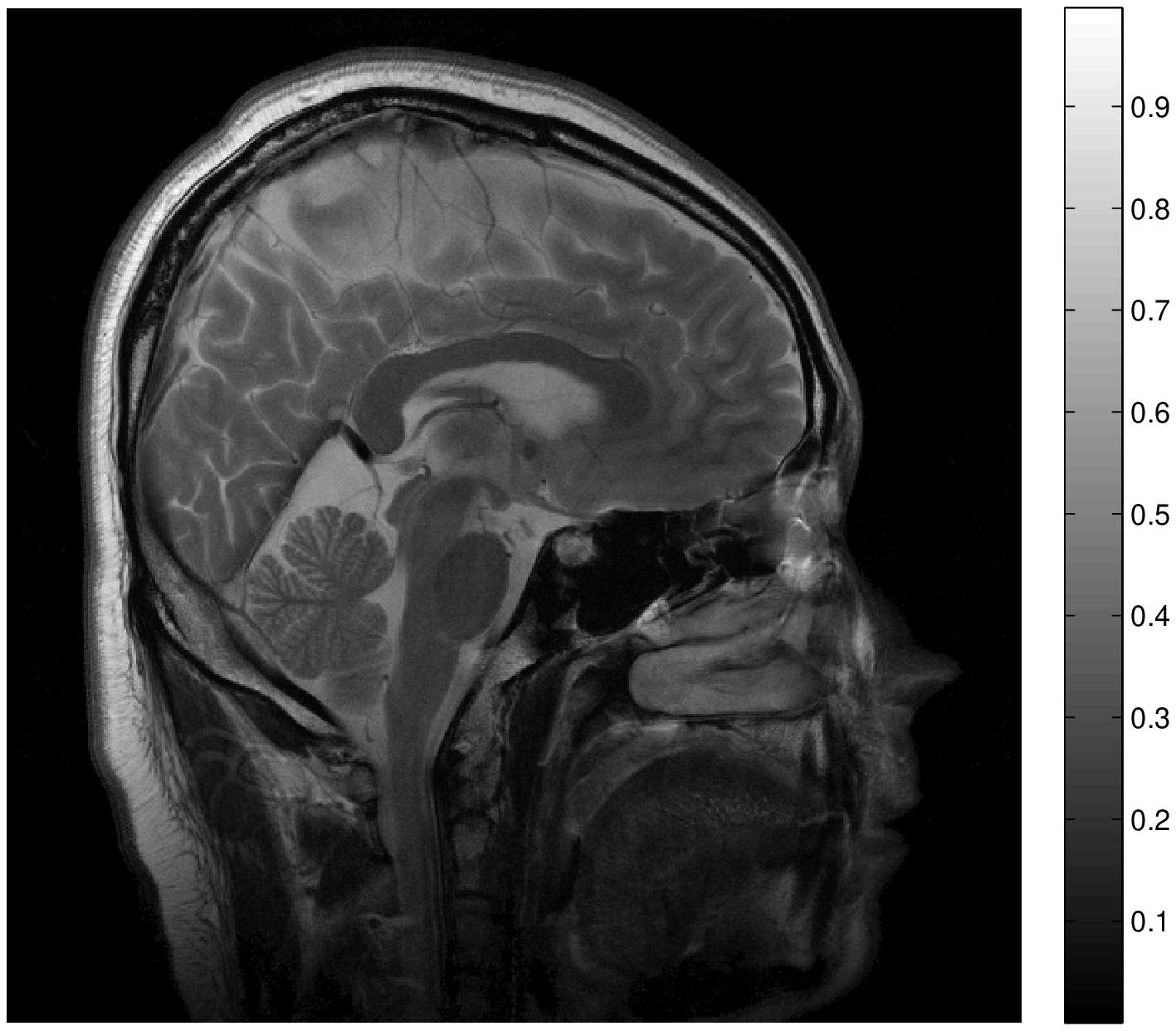}
	\end{subfigure}
	\begin{subfigure}[b]{.3\linewidth}
		\includegraphics[width=\linewidth]{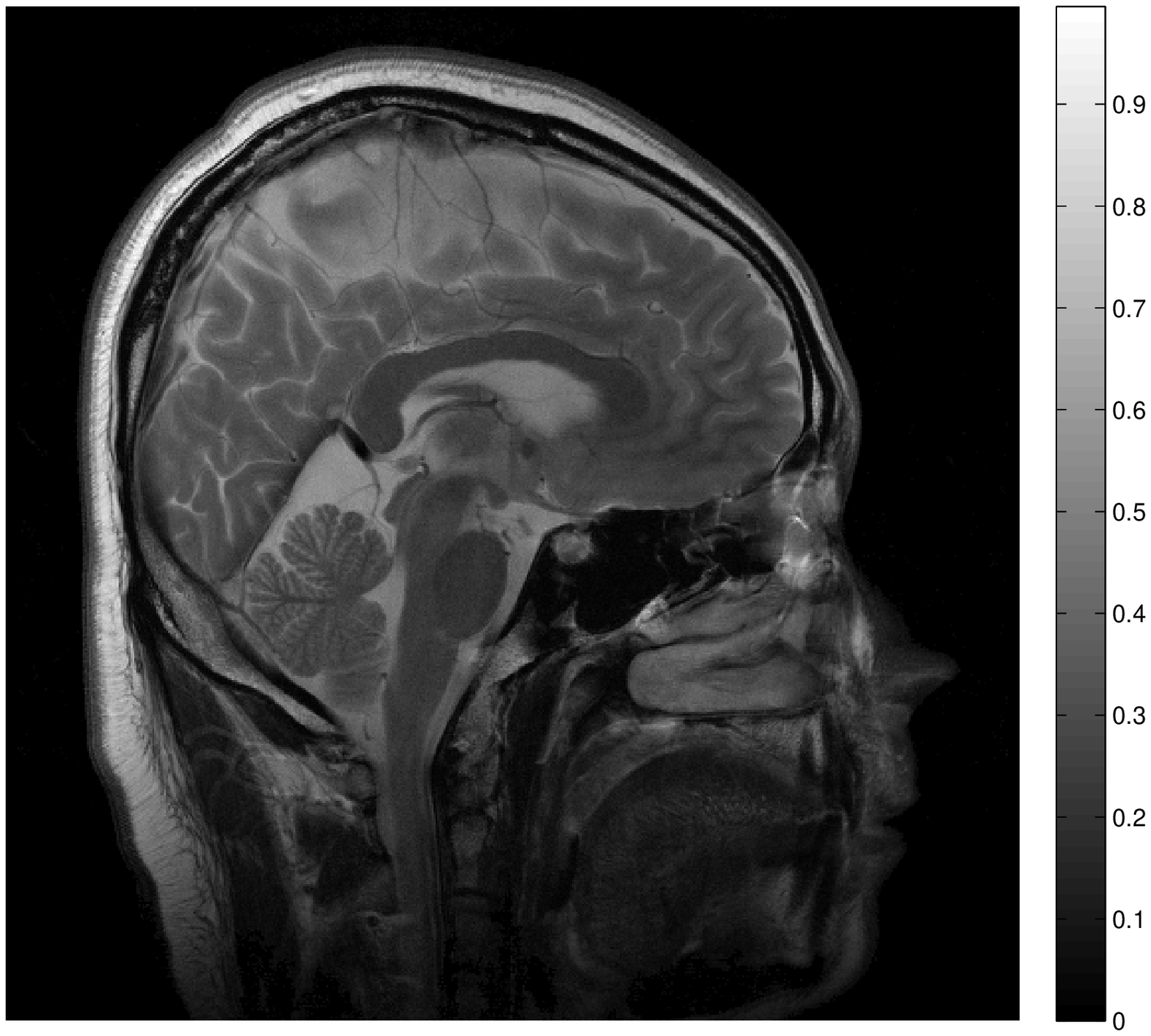}
	\end{subfigure}
	\begin{subfigure}[b]{.3\linewidth}
		\includegraphics[width=\linewidth]{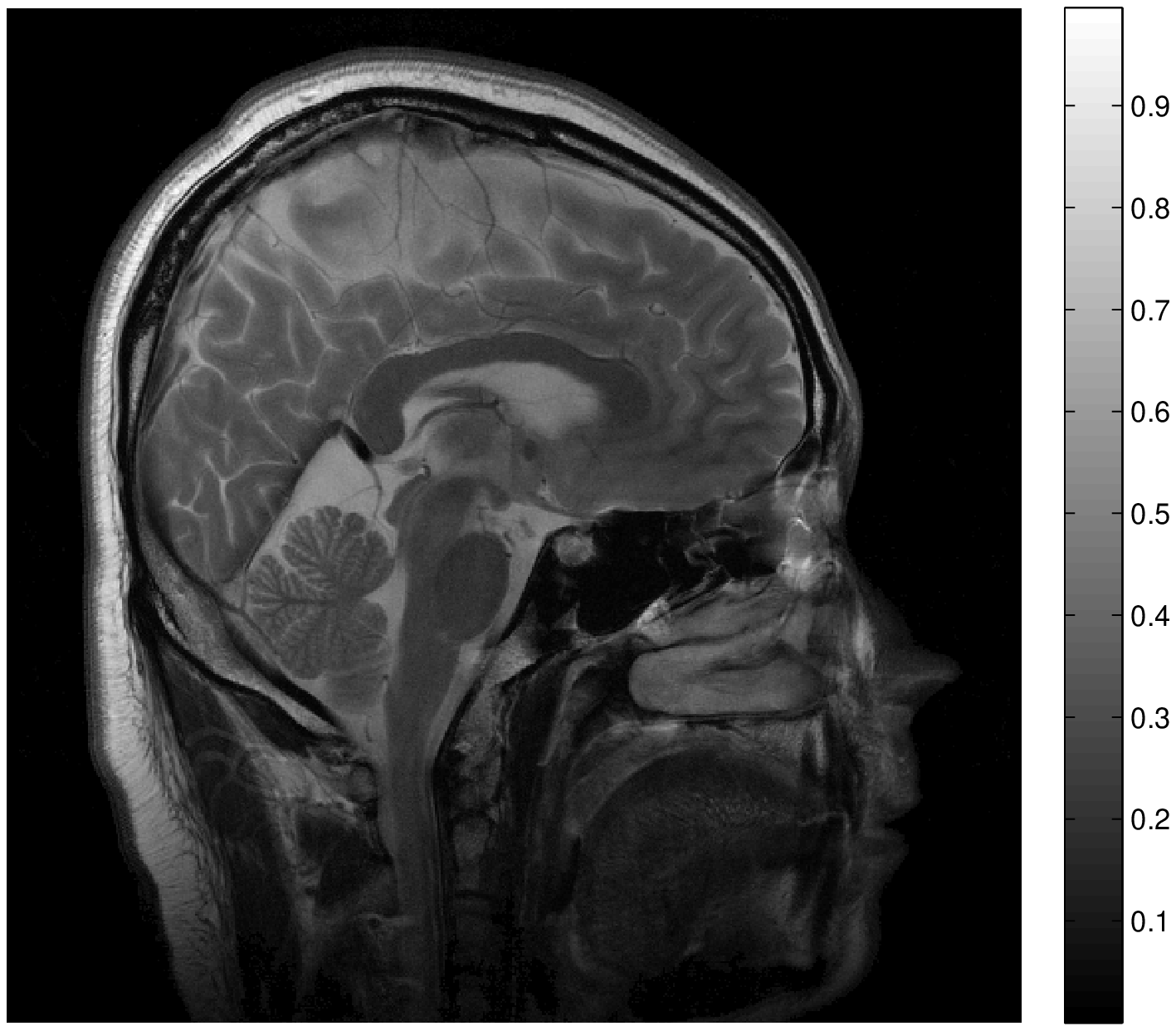}
	\end{subfigure}
	\\
	\begin{subfigure}[b]{.3\linewidth}
		\includegraphics[width=\linewidth]{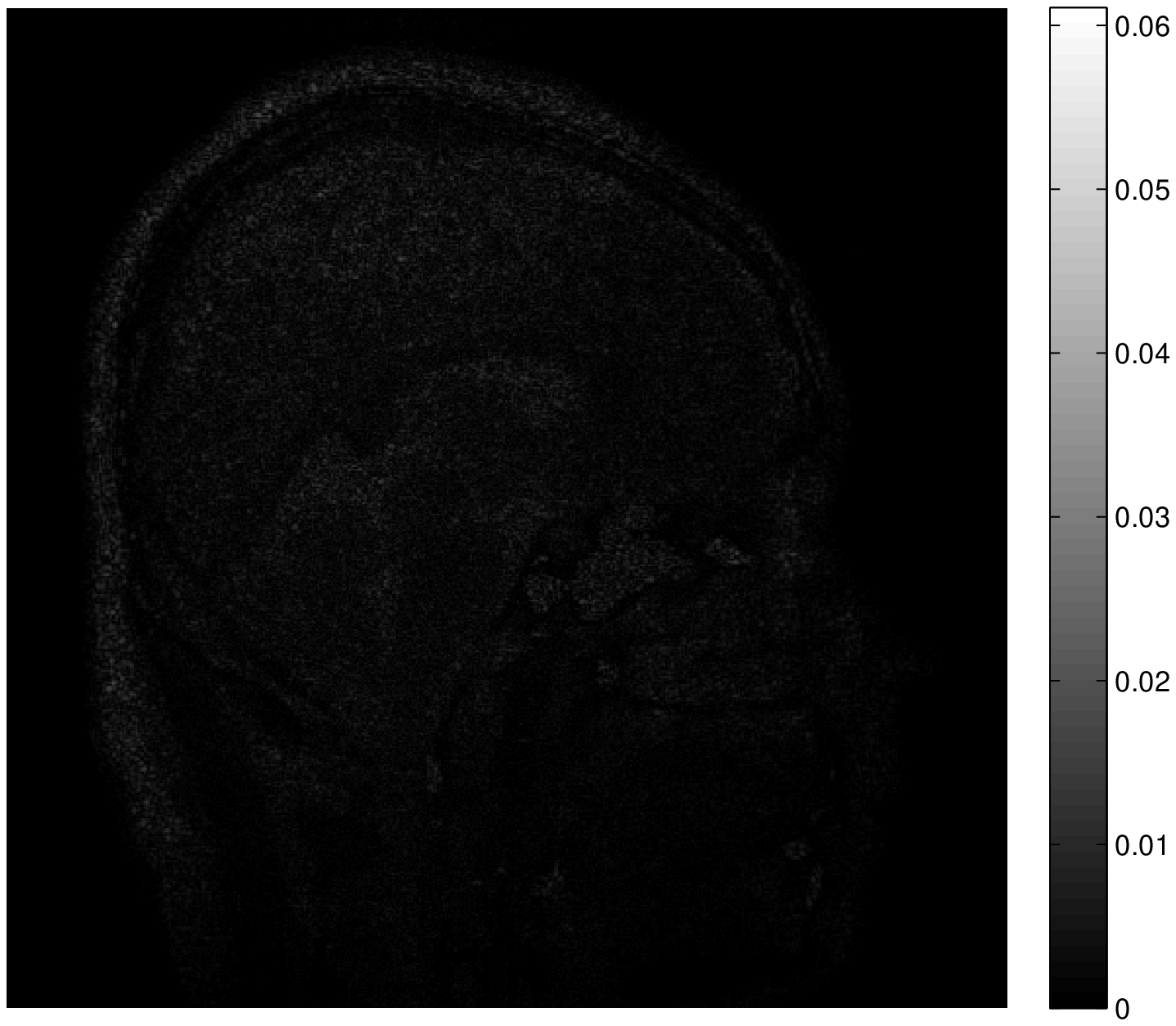}
	\end{subfigure}
	\begin{subfigure}[b]{.3\linewidth}
		\includegraphics[width=\linewidth]{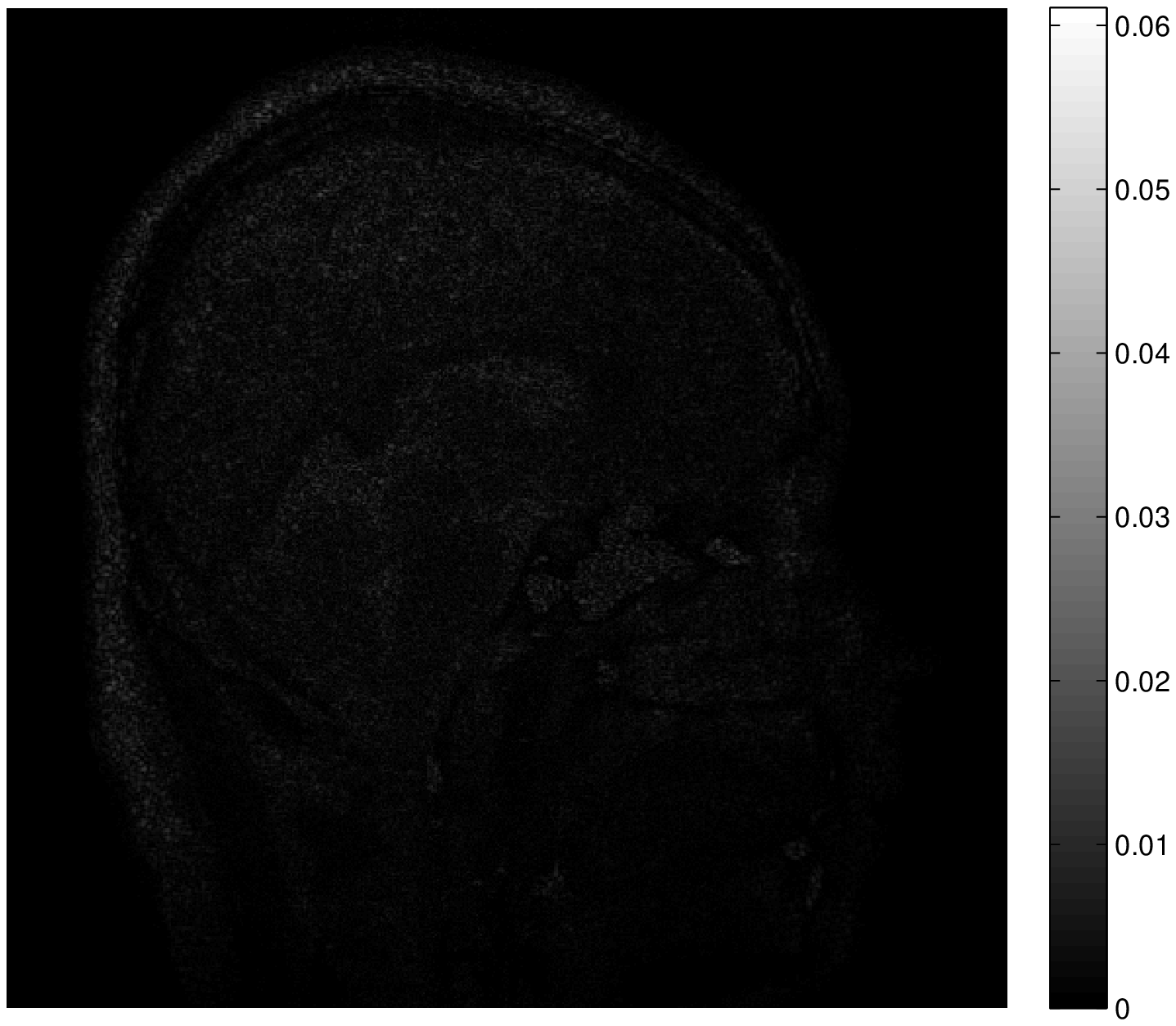}
	\end{subfigure}
	\begin{subfigure}[b]{.3\linewidth}
		\includegraphics[width=\linewidth]{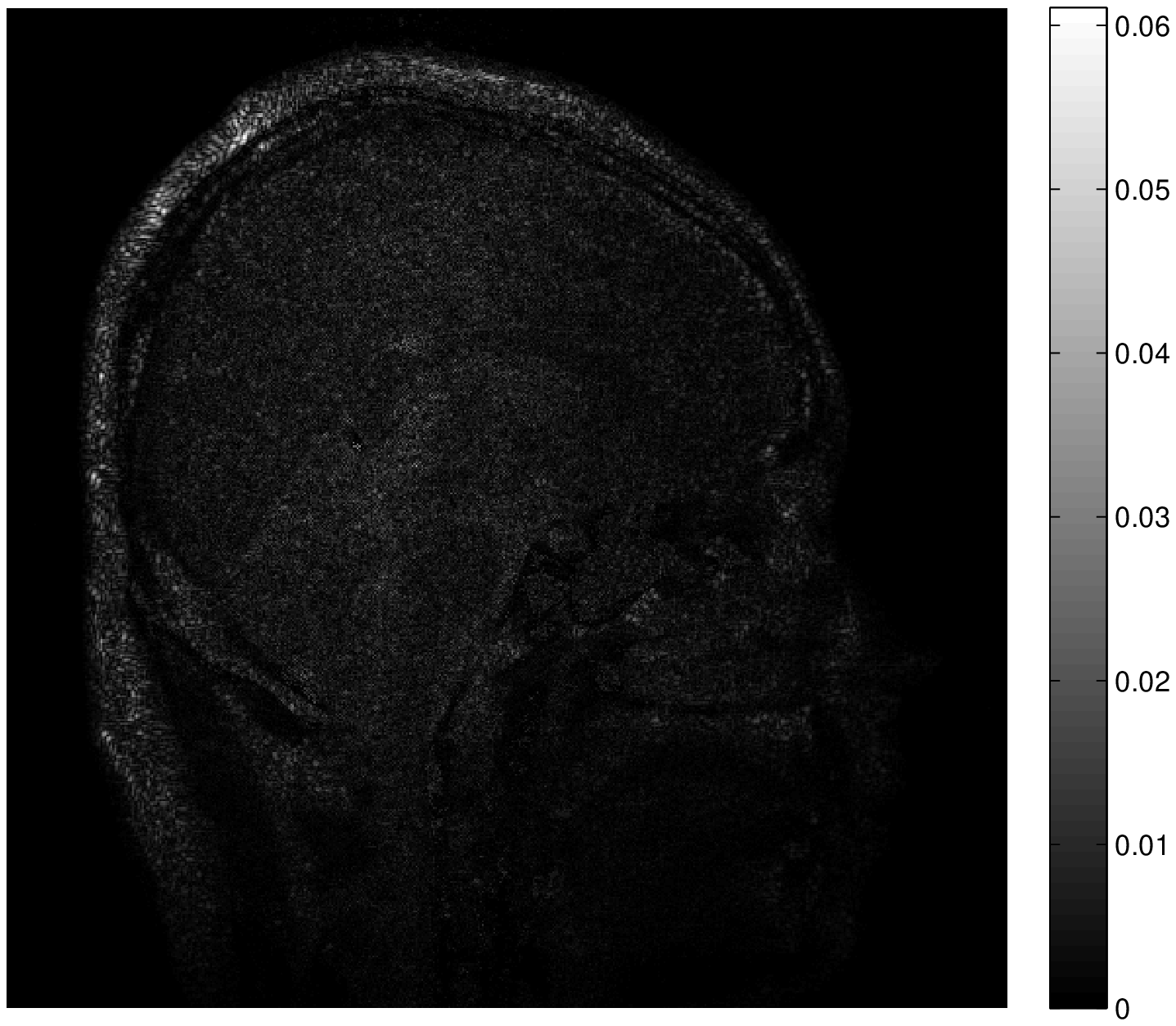}
	\end{subfigure}
	\caption{\footnotesize Comparison of AL-ADMM, ALP-ADMM and BOSVS in partially parallel image reconstruction (cont'd). From top to bottom: Reconstructed images and reconstruction errors in instances 2a and 2b, respectively. From left to right: AL-ADMM, ALP-ADMM and BOSVS.}
\end{figure}

\end{document}